\documentclass[sn-mathphys]{sn-jnl}

\jyear{2021}%

\theoremstyle{thmstyleone}%
\newtheorem{theorem}{Theorem}[section]
\newtheorem{proposition}[theorem]{Proposition}%
\newtheorem{corollary}[theorem]{Corollary}
\newtheorem{definition}[theorem]{Definition}

\theoremstyle{thmstyleone}%
\newtheorem{lemma}{Lemma}[section]

\newcommand{\calA}{{\mathcal{A}}}
\newcommand{\calB}{{\mathcal{B}}}
\newcommand{\calC}{{\mathcal{C}}}

\newcommand{\calF}{{\mathcal{F}}}

\newcommand{\calI}{{\mathcal{I}}}

\newcommand{\calK}{{\mathcal{K}}}
\newcommand{\calL}{{\mathcal{L}}}
\newcommand{\calM}{{\mathcal{M}}}
\newcommand{\calN}{{\mathcal{N}}}
\newcommand{\calO}{{\mathcal{O}}}
\newcommand{\calP}{{\mathcal{P}}}

\newcommand{\calR}{{\mathcal{R}}}
\newcommand{\calS}{{\mathcal{S}}}
\newcommand{\calT}{{\mathcal{T}}}
\newcommand{\calU}{{\mathcal{U}}}
\newcommand{\calV}{{\mathcal{V}}}

\newcommand{\N}{{\mathbb{N}}} 

\begin{document}

\title[A theory of integration for Ces\`aro limits]{A theory of integration for Ces\`aro limits}

\author*[1]{\fnm{Jonathan M.} \sur{Keith}}\email{jonathan.keith@monash.edu}

\author[1]{\fnm{Greg} \sur{Markowsky}}\email{greg.markowsky@monash.edu}

\affil[1]{\orgdiv{School of Mathematics}, \orgname{Monash University}, \orgaddress{\street{Wellington Road}, \city{Clayton}, \postcode{3800}, \state{VIC}, \country{Australia}}}


\abstract{The Ces\`aro limit - the asymptotic average of a sequence of real numbers - is an operator of fundamental importance in probability, statistics and analysis. Surprisingly, spaces of sequences with Ces\`aro limits have not previously been studied. This paper introduces spaces of such sequences, denoted $K_p(\calA)$, with the Ces\`aro limit acting as a kind of integral. 

The space $\calF$ comprised of all binary sequences with a Ces\`aro limit is studied first, along with the associated functional $\nu: \calF \rightarrow [0,1]$ mapping each such sequence to its Ces\`aro limit. It is shown that $\calF$ can be factored to produce a monotone class on which $\nu$ induces a countably additive set function. 

The space $K_p(\calA)$ is then defined, and a quotient denoted $\calK_p(\calA)$ is shown to be isometrically isomorphic, under certain conditions, to the function space $\calL_p(\N,\calA,\nu)$, where $\calA$ is a field of sets isomorphic to a subset of $\calF$, and $\nu$ is a finitely additive measure induced by the functional mentioned above. The Ces\`aro limit of an element of $K_p(\calA)$ is shown to be equal to its integral.

The complete $\calL_p(\N,\calA,\nu)$ spaces (and by implication, the $\calK_p(\calA)$ spaces isomorphic to them) are characterised, and a sufficient condition for these spaces to be separable is identified. }




\maketitle

\section{Introduction}\label{sec1}

Consider limits of the form
\[
\lim_{N \rightarrow \infty} \frac{1}{N} \sum_{n=1}^N x_n
\]
where $x_n \in \mathbb{R}$ for $n \in \N$ (here and throughout this paper, $\N$ denotes the natural numbers excluding 0). These are known as Ces\`aro limits (see \cite{bishop2014} for example) or sometimes as Ces\`aro means or Ces\`aro averages (as in \cite{crismale2017}), and arise naturally in multiple mathematical fields, including statistics, probability, functional and general analysis, and in the study of stochastic processes, particularly ergodic processes. They are important in many applications (see list of references in~\cite{bishop2014}). They are named for mathematician Ernesto Ces\`aro (1859-1906), who was not the first to consider the asymptotic properties of a sequence of averages, but used them to define a generalised limit for divergent series (see \cite{ferraro1999first}, and the references therein, for an interesting historical account).

Ces\`aro limits seem to be a kind of expectation - at any rate, ergodic theorems, for example that of Birkhoff~\cite{birkhoff1931}, establish one kind of relationship between a Ces\`aro limit and an expectation operator. Moreover, one may think of any sequence of real numbers $x = (x_1, x_2, \ldots)$ as a function $x: \N \rightarrow \mathbb{R}$. It seems reasonable then to consider a space of such functions for which Ces\`aro limits exist, and to identify conditions under which the Ces\`aro limit may be regarded as a kind of integral, or expectation operator. Such spaces could potentially be useful in the analysis of ergodic processes, especially if one can construct complete, separable, normed (or pseudo-normed) linear spaces on which an ergodic process may be regarded as a random element. This potential application is the motivation for introducing the $K_p$ function spaces analysed in this paper, although connections to ergodic processes will not be explored herein. At any rate, these spaces are of interest in their own right, and may have wider applications than this original motivation.

This paper first considers binary sequences with Ces\`aro limits. The collection of such sequences may be identified with a collection $\calF$ of subsets of $\N$, defined in Section~\ref{F_and_nu}. A set function $\nu$ that maps such subsets to their corresponding Ces\`aro limits is also defined. The basic properties of $\calF$ and $\nu$ are enumerated in that section. It turns out that $\nu$ has many of the properties of a finitely additive measure, also known as a {\em charge} (see \cite{bhaskararao1983}). However, $\calF$ is not a field, and thus $\nu$ is not a charge unless restricted to a field of sets contained in $\calF$. Section~\ref{F_and_nu} also introduces the collection of {\em null sets} $\calN$, comprised of subsets of $\N$ that induce binary sequences with zero Ces\`aro limits.

Although $\nu$ is not countably additive on $\calF$, it turns out that chains (totally ordered sets) in $\calF$ on which $\nu$ is countably additive have a number of useful properties. Two sections of the paper are devoted to exploring the properties of such chains. Section~\ref{uniform_convergence_section} characterises such chains in terms of uniform convergence to Ces\`aro limits. Section~\ref{null_modification_section} develops a construction that is here called a {\em null modification}. This construction modifies the elements of a chain of sets in $\calF$ by adding and/or removing null sets to produce a new chain on which $\nu$ is countably additive. This section makes frequent reference to Boolean algebras and their quotients: a brief review of this topic is therefore included in Section~\ref{Boolean_review}, with special attention to the Boolean quotient $\calP(\N) / \calN$. 

Section~\ref{quotient_space} considers the space $[\calF]$ - the image of $\calF$ under the quotient map $\xi : \calP(\N) \rightarrow \calP(\N) / \calN$. The collection of equivalence classes $[\calF]$ is shown to be a monotone class in the Boolean quotient $\calP(\N) / \calN$. This is a useful insight into the structure of $\calF$ because it implies that every field of sets in $\calF$ can be extended in such a way that the extension maps to a countably complete subalgebra of $[\calF]$, as a consequence of the monotone class theorem for Boolean algebras, which is reviewed in Section~\ref{Boolean_review}. This version of the monotone class theorem is an abstraction of the well known version for fields of sets, and is a slight generalisation of similar results in the literature. Consequently, a full proof is presented in the preprint version of this paper. 

The function space $K_p(\calA)$ is defined in Section~\ref{cesaro_integrals}, where $\calA \subset \calF$ is a field of sets (not necessarily a $\sigma$-field) and $p \in [1,\infty)$. The space $K_p(\calA)$ contains real-valued functions on a {\em charge space} $(\N,\calA,\nu)$. Charge spaces are generalisations of measure spaces, and are comprised of a sample space (here $\N$), a field of subsets of the sample space (here $\calA \subset \calF$), and a finitely additive measure, also known as a {\em charge}. A comprehensive introduction to the theory of charges is provided in~\cite{bhaskararao1983}. A concise summary is provided in~\cite{keith2022}, including an introduction to $L_p$ spaces and their quotients under equivalence almost everywhere (denoted $\calL_p$ spaces). That paper also extends the theory of bounded charges, presenting new characterisations of key properties that are applied in the present paper. 

Section~\ref{cesaro_integrals} explores the properties of $K_p(\calA)$ spaces and their quotients (denoted $\calK_p(\calA)$ spaces). In particular, $K_p(\calA)$ is shown to be a dense subspace of $L_p(\N,\calA,\nu)$, such that the integral of any function in $K_p(\calA)$ corresponds to the Ces\`aro limit of the corresponding sequence. In the rest of this paper, the function space $L_p(\N,\calA,\nu)$ and its quotient $\calL_p(\N,\calA,\nu)$ will be abbreviated as $L_p(\calA)$ and $\calL_p(\calA)$ respectively, since the sample space $\N$ and charge $\nu$ are assumed throughout. Under certain conditions, the quotient spaces $\calK_p(\calA)$ and $\calL_p(\calA)$ are shown to be isomorphic.

The remaining sections consider $L_p(\calA)$ spaces that are complete and separable. Section~\ref{Lp_countably_additive} characterises those $L_p(\calA)$ spaces for which $\calA$ is a $\sigma$-field and $\nu$ is countably additive; these spaces have a particularly simple form, lacking generality. Section~\ref{complete_calF} characterises more general $L_p(\calA)$ spaces that are complete, and identifies an isometric isomorphism between a complete $L_p(\calA)$ space and a conventional (ie. Lebesgue) function space. Section~\ref{separability_section} identifies a sufficient condition for an $L_p(\calA)$ space to be separable, specifically if $\calA$ is in a certain sense generated by a chain, and shows that this is equivalent to being generated (in a certain sense) by a countable sub-field of $\calA$. 

\section{Ces\`aro limits of binary sequences} \label{F_and_nu}

Let $\calP(X)$ denote the power set of an arbitrary set $X$ and consider the following definitions.

\begin{definition}
\label{calFdefn}
For any $A \in \calP(\N)$, define a {\em partial average}
\[
\nu_N(A) := \frac{1}{N} \sum_{n=1}^N I_A(n),
\] 
for each $N \in \N$, where $I_A$ is the indicator function for the set $A$. 
\end{definition}

\begin{definition} \label{upperandlower}
For any $A \in \calP(\N)$, define the {\em upper and lower Ces\`aro limits} respectively as
\[
\nu^+(A) := \limsup_{N \rightarrow \infty} \nu_N(A) \mbox{ and } \nu^-(A) := \liminf_{N \rightarrow \infty} \nu_N(A).
\]
\end{definition}

Naturally, the sets upon which $\nu^+$ and $\nu^-$ coincide are of particular interest, motivating the following definition.

\begin{definition}
Let $\calF$ be the collection of subsets $A \subseteq \N$ such that
\[
\nu(A) := \lim_{N \rightarrow \infty} \nu_N(A)
\] 
exists in the interval $[0,1]$. That is, $\calF$ is the collection of subsets $A$ of $\N$ for which the {\em Ces\`aro limit} of the binary sequence $x_n := I_A(n)$ (for each $n \in \N$) exists, and $\nu(A)$ is that limit for any $A \in \calF$.
\end{definition}

Note that 
\begin{enumerate}
\item $0 \leq \nu^-(A) \leq \nu^+(A) \leq 1$,  
\item $\nu^-(A) = \nu^+(A) \iff A \in \calF$ and 
\item $A \in \calF \implies \nu(A) = \nu^+(A)$.
\end{enumerate}

While $\nu_N$ is a measure on the power set of $\N$, being simply a scaling of counting measure on a finite set, it is evident that $\nu$ is not; to see this, note that any singleton set $\{k\}$ will have $\nu(\{k\})=0$, but $\cup_{k \in \N} \{k\} = \N$, so that $\nu$ is not countably additive. It is, however, finitely additive, and is thus a charge when restricted to fields of sets contained in $\calF$. It may be helpful to consider specific examples of sets for which $\nu$ exists. If $D_m$ is the set of multiples of an integer $m$, then $\nu(D_m) = \frac{1}{m}$. The same holds if $D^r_m$ is the set of all numbers equal to $r$ modulo $m$, and by taking unions of such sets one can obtain a set with any rational number as its charge. To obtain an irrational number $s$ as a charge is only slightly harder, and it can be achieved by the following algorithm. Start with $A_1 := \{1\}$, then for $N \geq 2$, if $\frac{1}{N} \sum_{n=1}^N I_{A_N}(n)<s$ let $A_{N+1} := A_N \cup\{N+1\}$, else let $A_{N+1} := A_N$. It is straightforward to show that if $A := \cup_{n \in \N} A_n$, then $\nu(A) = s$. 

The following method of describing sets is also useful for constructing specific examples. Let $z_n$ be a sequence of positive integers for $n \geq 2$, and let $z_1$ be a non-negative integer. Let $Z_n = \sum_{j=1}^n z_n$. Then let $A$ be defined by
\[
I_A(n):=\begin{cases}
1 & \mbox{ if $Z_{2k-1}+1 \leq n \leq Z_{2k}$ for some $k\geq 1$}  \\
0 & \mbox{if $n \leq z_1$ or $Z_{2k} + 1 \leq n \leq Z_{2k+1}$ for some $k \geq 1$}
\end{cases}.
\]
In words, $I_A(n)$ is $z_1$ zeroes, followed by $z_2$ ones, followed by $z_3$ zeroes, etc. It may be checked that $\nu^+(A) = \limsup_{N \to \infty} \nu_{Z_{2N}}(A) = \limsup_{N \to \infty} \frac{\sum_{j=1}^N z_{2j}}{Z_{2N}}$ and $\nu^-(A) = \liminf_{N \to \infty} \nu_{Z_{2N-1}}(A) = \liminf_{N \to \infty} \frac{\sum_{j=1}^{N-1} z_{2j}}{Z_{2N-1}}$. To form a simple example of a set which is not in $\calF$, let $z_n = 2^{n-1}$. Then $Z_n = 2^n-1$ and $\sum_{j=1}^N z_{2j} = \sum_{j=1}^N 2^{2j-1} = \frac{2}{3}(2^{2N} -1)$. Furthermore, 
\[
\nu^+(A) = \limsup_{N \to \infty} \frac{\frac{2}{3}(2^{2N} -1)}{2^{2N}-1} = \frac{2}{3}
\]
and \[
\nu^-(A) = \liminf_{N \to \infty} \frac{\frac{2}{3}(2^{2(N-1)} -1)}{2^{2N-1}-1} = \frac{1}{3}.
\]

This set can then be used to construct two sets $B, C \in \calF$, such that  $B \cap C \notin \calF$, thereby showing that $\calF$ is not a field. Let $B$ be the set of all even numbers, and let $C$ be defined by
\[
I_C(n):=\begin{cases}
1 & \mbox{ if $n$ is even and $\frac{n}{2} \in A$, or $n$ is odd and $\frac{n+1}{2} \notin A$ }  \\
0 & \mbox{otherwise}
\end{cases}.
\]
It is clear $B \in \calF$, with $\nu(B) = \frac{1}{2}$, and the same conclusion follows for $C$ upon noting that exactly one of $\{2k-1,2k\}$ lies in $C$ for every $k \in \N$. However, $B \cap C = 2A$, the set of the doubles of elements of $A$, and therefore $B \cap C \notin \calF$. 

For the purpose of intuition, it is profitable to think of sets in $\calF$ in this manner, as defined by the concatenation of alternating strings of zeroes and ones of variable length. The following result describes the lengths of these strings allowable for a set to be in $\calF$. Consider any $A \subseteq \N$. For each $N \in \N$, define $P_A(N)$ to be the smallest integer $k > 0$ such that $I_A(N+k) = 1$ or define $P_A(N) = \infty$ if no such integer exists. Similarly, define $Q_A(N)$ to be the smallest integer $k > 0$ such that $I_A(N+k) = 0$ or define $Q_A(N) = \infty$ if no such integer exists.

\begin{proposition} \label{useful_lemma}
Consider $A \subseteq \N$.
\begin{enumerate}
\item If $A$ contains $F \in \calF$ such that $\nu(F) > 0$, then $\nu^-(A) > 0$ and $P_A$ is $o(N)$.
\item If $A$ is contained in $F \in \calF$ such that $\nu(F) < 1$, then $\nu^+(A) < 1$ and $Q_A$ is $o(N)$.
\end{enumerate}
\end{proposition}

\noindent {\bf Proof: } 
Suppose $A$ contains $F \in \calF$ such that $\nu(F) > 0$. Then $\nu^-(A) \geq \nu^-(F) = \nu(F) > 0$. Moreover $F$ contains infinitely many integers, so $P_A(N) \leq P_F(N) < \infty$ for all $N \in \N$. If $F$ excludes only finitely many integers then $P_A(N) = P_F(N) = 1$ for large enough $N$, making $P_A(N)$ trivially $o(N)$, so assume $F$ excludes infinitely many integers. There are thus infinitely many positive integers $N_1 < N_2 < \ldots$ such that $I_F(N_i + 1) \neq I_F(N_i )$. Note
\[
\nu(F) = \nu_{N_1}(F) + \sum_{i = 1}^{\infty} (\nu_{N_{i+1}}(F) - \nu_{N_i}(F))
\]
is a series consisting of alternating positive and negative terms corresponding respectively to runs of ones and zeros in the sequence $( I_F(1), I_F(2), \ldots )$. Since the series converges, terms corresponding to runs of zeros must decrease in magnitude to 0.

Now consider any $N$ with $P_F(N) = k+1 > 0$.  Then 
\[
\nu_{N+k}(F) - \nu_N(F) = \nu_N(F)\frac{N}{N+k} - \nu_N(F) = - \nu_N(F) \left( \frac{k}{N+k} \right).
\]
In particular, if $N_i$ corresponds to the end of a run of ones, then the subsequent run of zeros contributes a term
\[
- \nu_{N_i}(F) \left( \frac{P_F(N_i) - 1}{N+P_F(N_i) - 1} \right)
\]
to the above series, and since $\nu_{N_i}(F) \rightarrow \nu(F) > 0$, these terms can only go to 0 if 
\[
\frac{P_F(N_i) - 1}{N+P_F(N_i) - 1} \rightarrow 0
\]
implying $P_A(N) \leq P_F(N)$ is $o(N)$. The second part of the lemma follows by applying the first part to $A^c$. 
\qed


It can be shown that the term $o(N)$ in Proposition~\ref{useful_lemma} cannot be replaced by $o(N^{1-\epsilon})$ for any $\epsilon > 0$, as follows. For positive integer $q$ let $z_n = n^q$, and form a set $A$ by the method described earlier in this section. Then, by comparing the sum with an integral, the easy estimates $\frac{N^{q+1}}{q+1} \leq Z_{N} \leq \frac{(N+1)^{q+1}}{q+1}$ are obtained, and since $\sum_{j=1}^N z_{2j} = 2^q \sum_{j=1}^N z_{j}$ it follows also that $\frac{2^q N^{q+1}}{q+1} \leq \sum_{j=1}^N z_{2j} \leq \frac{2^q (N+1)^{q+1}}{q+1}$. Thus, $\nu^+(A) = \limsup_{N \to \infty} \frac{\sum_{j=1}^N z_{2j}}{Z_{2N}} \leq \limsup_{N \to \infty} \frac{2^q (N+1)^{q+1}}{(2N)^{q+1}} = \frac{1}{2}$, and $\nu^-(A) = \liminf_{N \to \infty} \frac{\sum_{j=1}^{N-1} z_{2j}}{Z_{2N-1}} \geq \liminf_{N \to \infty} \frac{2^q (N-1)^{q+1}}{(2N-1)^{q+1}} = \frac{1}{2}$, so $A \in \calF$. However, $P_A(Z_{2N}) = (2N+1)^q + 1$, and $Z_{2N} \leq \frac{(2N+1)^{q+1}}{q+1}$, hence $\frac{P_A(Z_{2N})^{(q+1)/q}}{Z_{2N}}$ is bounded below by a positive constant for any $q > 0$, giving the result.

Sets that $\nu$ maps to zero play an important role in this paper.

\begin{definition} \label{null_sets}
The {\em null sets} in $\calF$ are the elements of $\calN := \{ A \in \calF : \nu(A) = 0 \}$. 
\end{definition}

Null sets are easy to find. The set of square numbers is a null set, as is the set of cubes, etc. The set of powers of 2, or of any other base, is a null set. The set of primes is shown to be a null set by the prime number theorem. 

The next proposition, which describes key properties of $\calF$ and the Ces\'aro limits of its elements, is of fundamental importance in subsequent sections. 

\begin{proposition}
\label{calFproperties}
The collection $\calF$ and set functions $\nu^+$, $\nu^-$ and $\nu$ have the following properties.
\begin{enumerate}
\item $\emptyset, \N \in \calF$, with $\nu(\emptyset) = 0$ and $\nu(\N) = 1$.
\item For all $A \in \calF$, $A^c \in \calF$ with $\nu(A^c) = 1 - \nu(A)$.
\item For all $A, B \in \calF$, $A \cup B \in \calF$ if and only if $A \cap B \in \calF$, and if either is true then $\nu(A \cup B) = \nu(A) + \nu(B) - \nu(A \cap B)$.
\item If $A, B \in \calP(\N)$, then $\nu^+(A \cup B) \leq \nu^+(A) + \nu^+(B)$. If $A$, $B$, and $A \cup B$ are all in $\calF$, then $\nu(A \cup B) \leq \nu(A) + \nu(B)$.
\item For $A, B \in \calP(\N)$ such that $A \subseteq B$, 
\begin{enumerate}
\item $\nu^+(A) \leq \nu^+(B)$,  
\item $\nu^+(B \setminus A) \geq \nu^+(B) - \nu^+(A)$, 
\item $\nu^-(A) \leq \nu^-(B)$, and 
\item $\nu^-(B \setminus A) \leq \nu^-(B) - \nu^-(A)$. 
\end{enumerate}
If in addition $A, B \in \calF$, then 
\begin{enumerate}
\item $\nu(A) \leq \nu(B)$, 
\item $B \setminus A \in \calF$,
\item $\nu(B \setminus A) = \nu(B) - \nu(A)$, and
\item $\nu(A) = \nu(B) \iff B \setminus A \in \calN$.
\end{enumerate}
\item If $A \in \calN$, then any $B \subseteq A$ satisfies $B \in \calN$. Consequently, for any $C \in \calP(\N)$, 
\begin{enumerate}
\item $A \cap C, A \setminus C \in \calN$, and
\item $\nu^+(A \cup C) = \nu^+(C \setminus A) = \nu^+(C)$.
\end{enumerate}
If $C \in \calF$, then $A \cup C, C \setminus A \in \calF$.
\item For pairwise disjoint sets $A_1, \ldots, A_K \in \calF$, $A_1 \cup \ldots \cup A_K \in \calF$ with 
\[
\nu(A_1 \cup \ldots \cup A_K) = \sum_{k=1}^K \nu(A_k).
\] 
\item Consider $\calC \subseteq \calP(\N)$. Then 
\[
\nu^-(\bigcup \calC) \geq \sup\{ \nu^-(A) : A \in \calC \}.
\]
In particular, if $\{ A_k \}_{k=1}^{\infty} \subset \calF$ are pairwise disjoint, then
\[
\nu^-(\cup_{k=1}^{\infty} A_k) \geq \sum_{k=1}^{\infty} \nu(A_k).
\]
\item Consider a chain $\calC \subset \calP(\N)$ such that $\nu_N(A) \leq \nu^+(A)$ for all $N \in \N$ and $A \in \calC$. Then 
\begin{enumerate}
\item $\nu^+(\bigcup \calC) = \sup_{A \in \calC} \nu^+(A)$, and
\item if $\calC \subset \calF$, then $\bigcup \calC \in \calF$ with $\nu(\bigcup \calC) = \sup_{A \in \calC} \nu(A)$.
\end{enumerate}
In particular, if $\{ A_k \}_{k=1}^{\infty} \subset \calF$ are pairwise disjoint with $\nu_N(A_k) \leq \nu(A_k)$ for all $k, N \in \N$, then $\cup_{k=1}^{\infty} A_k \in \calF$ and $\nu(\cup_{k=1}^{\infty} A_k) = \sum_{k=1}^{\infty} \nu(A_k)$.
\if2 {\tcr Very interesting fact. I think it can be generalized somewhat, for instance is all but finitely many satisfy $\nu_N(A_k) \leq \nu(A_k)$, or if $\nu_N(A_k) \leq \nu(A_k)+ \epsilon_k$ where $\sum_{k=1}^\infty \epsilon_k < \infty$. But I suspect that the statement given probably covers the important cases, and thus probably best not to complicate it.} {\color{blue} (HHH One option is to split it into two properties: this simpler version and a more general one.)} \fi 
\item Consider $\calC \subseteq \calF$. If $\sup \{ \nu(A) : A \in \calC \} = 1$, then $\bigcup \calC \in \calF$ and $\nu(\bigcup \calC) = 1$. In particular, if $\{ A_k \}_{k=1}^{\infty} \subset \calF$ are pairwise disjoint with $\sum_{k=1}^{\infty} \nu(A_k) = 1$, then $\cup_{k=1}^{\infty} A_k \in \calF$ and $\nu(\cup_{k=1}^{\infty} A_k) = 1$.
\end{enumerate}
\end{proposition}

\noindent {\bf Proof: } 
Property 1 is trivial. For Property 2, consider any $A \in \calF$ and note 
\[
\lim_{N \rightarrow \infty} \nu_N(A^c) = \lim_{N \rightarrow \infty} \frac{1}{N} \sum_{n=1}^N (1 - I_{A}(n)) 
= 1 - \lim_{N \rightarrow \infty} \nu_N(A)
\]
hence $A^c \in \calF$ and $\nu(A^c) = 1 - \nu(A)$.

For 3, consider $A,B \in \calF$ and note $I_{A \cup B} = I_A + I_B - I_{A \cap B}$, hence
\[
\lim_{N \rightarrow \infty} \nu_N(A \cup B) = \nu(A) + \nu(B) - \lim_{N \rightarrow \infty} \nu_N(A \cap B)
\]
if either limit exists, and the statement follows immediately.

For 4, note $I_{A \cup B} \leq I_A + I_B$, hence $\nu_N(A \cup B) \leq \nu_N(A) + \nu_N(B)$ for all $N \in \N$. The first statement then follows by taking the $\limsup$ and the second by taking limits as $N \rightarrow \infty$. (The second statement alternatively follows from 3).

For 5, note $I_A \leq I_B$, hence $\nu^+(A) \leq \nu^+(B)$ and $\nu^-(A) \leq \nu^-(B)$. Also $I_B = I_A + I_{B \setminus A}$, hence 
\[
\limsup_{N \rightarrow \infty} \nu_N(B) \leq \limsup_{N \rightarrow \infty} \nu_N(A) + \limsup_{N \rightarrow \infty} \nu_N(B \setminus A).
\]
That is, $\nu^+(B \setminus A) \geq \nu^+(B) - \nu^+(A)$. That $\nu^-(B \setminus A) \leq \nu^-(B) - \nu^-(A)$ is shown similarly. If $A, B \in \calF$, then $\nu(A) = \nu^+(A) \leq \nu^+(B) = \nu(B)$ and
\[
\lim_{N \rightarrow \infty} \nu_N(B \setminus A) = \lim_{N \rightarrow \infty} \nu_N(B) - \lim_{N \rightarrow \infty} \nu_N(A) = \nu(B) - \nu(A) 
\]
thus $B \setminus A \in \calF$ and $\nu(B \setminus A) = \nu(B) - \nu(A)$. Hence $\nu(B) - \nu(A) = 0 \iff \nu(B \setminus A) = 0$. 

For 6, note for every $N \in \N$,
\[
0 \leq \nu_N(B) \leq \nu_N(A)
\]
and $B \in \calN$ follows by letting $N \rightarrow \infty$. The consequences 6(a) follow because $A \cap C \subseteq A$ and $A \setminus C \subseteq A$. For 6b, note
\[
\nu^+(C) = \nu^+(C) - \nu^+(A) \leq \nu^+(C \setminus A) \leq \nu^+(A \cup C) \leq \nu^+(A) + \nu^+(C) = \nu^+(C).
\]
If $C \in \calF$, $A \cup C \in \calF$ by 3, and $C \setminus A \in \calF$ by 5.

For 7, use Property 3 with $A_1 \cap A_2 = \emptyset$, to conclude that $A_1 \cup A_2 \in \calF$ and $\nu(A_1 \cup A_2) = \nu(A_1) + \nu(A_2)$. The property then follows by induction.

For 8, by Property~5, $\nu^-(\bigcup \calC) \geq \nu^-(A)$ for all $A \in \calC$. Hence $\nu^-(\bigcup \calC) \geq \sup\{ \nu^-(A) : A \in \calC \}$. The second part of 8 follows by defining $B_j := \cup_{k=1}^j A_k$, so that by Property~7, $B_j \in \calF$ with $\nu(B_j) = \sum_{k=1}^j \nu(A_k)$ for each $j \in \N$. Then apply the first part of 8 to $\calC := \{ B_j \}_{j=1}^{\infty}$.

For 9, note for any $N \in \N$, $\nu_N(\bigcup \calC) = \sup_{A \in \calC} \nu_N(A) \leq \sup_{A \in \calC} \nu^+(A)$, hence $\nu^+(\bigcup \calC) \leq \sup_{A \in \calC} \nu^+(A)$. Moreover, for any $\epsilon > 0$, one can choose $C \in \calC$ such that $\nu^+(\bigcup \calC) \geq \nu^+(C) > \sup_{A \in \calC} \nu^+(A) - \epsilon$. Letting $\epsilon \rightarrow 0$ gives $\nu^+(\bigcup \calC) = \sup_{A \in \calC} \nu^+(A)$. If $\calC \subset \calF$, then by 8,
\[
\nu^-(\bigcup \calC) \geq \sup_{A \in \calC} \nu(A) = \nu^+(\bigcup \calC),
\]
hence $\bigcup \calC \in \calF$ with $\nu(\bigcup \calC) = \sup_{A \in \calC} \nu(A)$. The last part of 9 follows by setting $\calC := \{ B_j \}_{j=1}^{\infty}$ as defined in the proof of 8. Then 9b gives $\cup_{k=1}^{\infty} A_k = \cup_{k=1}^{\infty} B_k \in \calF$ with $\nu(\cup_{k=1}^{\infty} A_k) = \sup_k \nu(B_k) = \sum_{k=1}^{\infty} \nu(A_k)$.

For 10, note that $1 \geq \nu^+(\bigcup \calC) \geq \nu^-(\bigcup \calC) \geq \sup\{ \nu(A) : A \in \calC \} = 1$, using Property~8. Hence $\nu^+(\bigcup \calC) = \nu^-(\bigcup \calC)$ and the first part follows. The second part of 10 follows by applying the first part to $\calC = \{ B_j \}_{j=1}^{\infty}$ defined above in the proof of 8.
\qed 

Property 9 is particularly important in what follows, so it may be helpful to discuss an example. Let $\calO$ denote the set of odd numbers which are at least 3, and for $k \in \N \cup \{0\}$ set $D_k = \{2^k m: m \in \calO\}$. Then it may be checked that $\nu_N(D_k) \leq \nu(D_k)$ for all $N$ and $\nu(D_k) = \frac{1}{2^{k+1}}$. The set $D := \cup_{k \in \N} D_k$ is all of $\N$ with a null set (the powers of $2$) removed, and thus $\nu(D) = 1 = \sum_{k = 0}^ \infty \nu(D_k)$, so that Property 9 holds. If $1$ were included in $\calO$ the same conclusion would hold, even though the sufficient condition $\nu_N(D_k) \leq \nu(D_k)$ would not (it is evident, however, that this or some other condition is needed to ensure countable additivity, since $\nu$ is not a measure). This example will appear briefly again in Section \ref{null_modification_section}, which contains a method for modifying sets by removing a null set so that Property 9 can be applied. 

Since $\calF$ contains $\emptyset$ and is closed under complements and finite disjoint unions, it is an object known as an {\em additive class}.

\section{Uniform convergence of chains in $\calF$}
\label{uniform_convergence_section}

A real-valued function on any set $X$ generates a chain of subsets consisting of inverse images of rays in $\mathbb{R}$. Thus the properties of chains consisting of subsets of $X$ are relevant to the study of real-valued functions on $X$. In this section and Section~\ref{null_modification_section}, two analysis tools for studying chains in $\calF$ are developed. The first of these is a characterisation of a certain class of chains in $\calF$ in terms of uniform convergence of partial averages over the sets in the chain. 

The following notation is helpful to describe this characterisation. Consider a chain of sets $\calT \subset \calP(\N)$ (that is, a collection of sets that is totally ordered by set inclusion). Let $\calT_{\cup}$ and $\calT_{\cap}$ denote the closure of $\calT$ under unions and intersections, respectively. That is, $\calT_{\cup} := \{ \bigcup \calC : \calC \subseteq \calT \}$ and $\calT_{\cap} := \{ \bigcap \calC : \calC \subseteq \calT \}$. Also set $\calT_* := \calT_{\cup} \cup \calT_{\cap}$. Some basic properties of $\calT_*, \calT_{\cup}$, and $\calT_{\cap}$ are the following.

\begin{proposition} \label{basicfacts}
Suppose $\calT \subset \calP(\N)$ is a chain. Then
\begin{enumerate}
    \item $\calT_{\cup}$, $\calT_{\cap}$, and $\calT_*$ are chains.
    
    \item $\calT_*$ is closed under unions and intersections.
\end{enumerate}
\end{proposition}

\noindent {\bf Proof: } 
To prove 1, first consider $\calT_{\cup}$. If $A \in \calT$ and $B = \bigcup \calC$ with $\calC \subseteq \calT$, then either $(i)$ $A \subseteq C$ for some $C \in \calC$, in which case $A \subseteq B$, or $(ii)$ $A \supseteq C$ for all $C \in \calC$, in which case $A \supseteq B$. Alternatively, suppose $A = \bigcup \calC_1$ and $B = \bigcup \calC_2$ with $\calC_1, \calC_2 \subseteq \calT$. By the previous argument, for every $C \in \calC_1$ either $C \subseteq \bigcup \calC_2$ or $C \supseteq \bigcup \calC_2$; if $C \subseteq \bigcup \calC_2$ for every $C \in \calC_1$ then $\bigcup \calC_1 \subseteq \bigcup \calC_2$, otherwise $C \supseteq \bigcup \calC_2$ for some $C \in \calC_1$, in which case $\bigcup \calC_1 \supseteq \bigcup \calC_2$. Thus $\calT_{\cup}$ is a chain. A complementary argument shows $\calT_{\cap}$ is a chain. To show $\calT_*$ is a chain, one must identify an ordering between $A = \bigcap \calC_1 \in \calT_{\cap}$ and $B=\bigcup \calC_2\in \calT_{\cup}$ with $\calC_1, \calC_2 \subseteq \calT$. If $C \supseteq D$ for every $C \in \calC_1, D \in \calC_2$, then $B \subseteq A$, otherwise $C \subseteq D$ for some $C \in \calC_1, D \in \calC_2$, in which case $A \subseteq B$. 

As for 2, suppose first that $\calC \subseteq \calT_{\cap}$ and consider $B = \bigcup \calC$. It is possible $B \in \calT_\cap$, but in this case there is nothing to prove since $\calT_\cap \subseteq \calT_*$, so assume $B \notin \calT_\cap$. Any $A \in \calC$ can be expressed as an intersection of sets in $\calT$, and if each of these sets contained $B$ then one would have $B \subseteq A$ and hence $B = A$, contradicting $B \notin \calT_\cap$. Hence there exists $A' \in \calT$ that contains $A$ but not $B$. Then $A \subseteq A' \subseteq B$, since $(\calT_{\cap})_{\cup}$ is a chain by the first part of this lemma. It follows that $B = \bigcup \{ A' \in \calT : A' \subseteq B \}$, and hence $B \in \calT_\cup$. This shows $(\calT_{\cap})_{\cup} \subseteq \calT_{\cap} \cup \calT_{\cup} = \calT_*$. Now suppose $\calC \subseteq \calT_*$. Then $\bigcup \calC = B_1 \cup B_2$ where $B_1 \in (\calT_{\cap})_{\cup} \subseteq \calT_*$ and $B_2 \in (\calT_{\cup})_{\cup} = \calT_{\cup} \subseteq \calT_*$. But $\calT_*$ is a chain, hence $B_1 \cup B_2$ is either $B_1$ or $B_2$. Either way $\bigcup \calC \in \calT_*$, hence $\calT_*$ is closed under unions. A complementary argument shows $\calT_*$ is closed under intersections.
\qed 

The following theorem identifies three alternative characterisations of a class of well behaved chains in $\calF$. The first characterisation implies countable additivity of the restriction of $\nu$ to the chain: it thus identifies chains in $\calF$ on which $\nu$ behaves like a measure. The other characterisations identify other useful properties of such chains, in particular, uniform convergence of partial averages of elements of the chain.

\begin{theorem} \label{uniformconvergence}
Let $\calT \subset \calF$ be a chain of sets. Then the following statements are logically equivalent.
\begin{enumerate}
\item $\calT_* \subset \calF$ and for any $\calC \subseteq \calT$, $\nu(\bigcup \calC) = \sup_{C \in \calC} \nu(C)$ and $\nu(\bigcap \calC) = \inf_{C \in \calC} \nu(C)$.
\item There exists a chain $\calU \subset \calF$ such that $\calT \subseteq \calU$ and $\nu(\calU) := \{ \nu(A) : A \in \calU \}$ is dense in $[0,1]$.
\item For every $\epsilon > 0$ there exists $N_{\epsilon} \in \N$ such that $\lvert \nu_N(A) - \nu(A) \rvert < \epsilon$ for all $A \in \calT$ and all $N > N_{\epsilon}$. 
\end{enumerate}
Moreover, if any of the three statements holds then 
\begin{enumerate}
\item $\calU_* \subset \calF$ and for any $\calC \subseteq \calU_*$, $\nu(\bigcup \calC) = \sup_{C \in \calC} \nu(C)$ and $\nu(\bigcap \calC) = \inf_{C \in \calC} \nu(C)$, 
\item $\nu(\calU_*) = [0,1]$, 
\item For every $\epsilon > 0$ there exists $N_{\epsilon} \in \N$ such that $\lvert \nu_N(A) - \nu(A) \rvert < \epsilon$ for all $A \in \calU_*$ and all $N > N_{\epsilon}$. 
\end{enumerate}
\end{theorem}

\noindent {\bf Proof: } 
($1 \implies 2$) Set $\calT_0 = \calT$. If there is no open interval $(b,c) \subset [0,1]$ of width at least $2^{-1}$ such that $(b,c) \cap \nu(\calT_0) = \emptyset$, then set $\calT_1 := \calT_0$. If there is exactly one such open interval, define sets $B := \bigcup \{ A \in \calT_0 : \nu(A) < b \}$ and  $C := \bigcap \{ A \in \calT_0 : \nu(A) > c \}$. Statement~1 implies $B, C \in \calT_* \subset \calF$, $\nu(B) \leq b < c \leq \nu(C)$ and $(\nu(B),\nu(C)) \cap \nu(\calT_0) = \emptyset$. Form a set $A$ (called a {\em midpoint set}) containing $B$ and every second element of the sequence generated by listing the elements of $C \setminus B$ in increasing order. Then $B \subset A \subset C$ and it is straightforward to show $A \in \calF$ with $\nu(A) = (\nu(B) + \nu(C))/2$. Set $\calT_1$ to be $\calT_0$ plus the midpoint set thus formed. If there are two disjoint open intervals of width at least $2^{-1}$, both of which have empty intersection with $\nu(\calT_0)$, then find the midpoint sets for both intervals and add them to $\calT_0$ to form $\calT_1$. Note there cannot be more than two such intervals. Then $\calT_1 \subset \calF$ is a chain that satisfies $\calT_1^* \subset \calF$ and for any $\calC \subseteq \calT_1$, $\nu(\bigcup \calC) = \sup_{C \in \calC} \nu(C)$ and $\nu(\bigcap \calC) = \inf_{C \in \calC} \nu(C)$. Moreover, $\nu(\calT_1) \subseteq [0,1]$ does not exclude any open intervals in $[0,1]$ of width at least $2^{-1}$. Proceeding inductively, one can generate a non-decreasing sequence of chains $\calT_1 \subseteq \calT_2 \subseteq \ldots$ such that $\nu(\calT_k)$ does not exclude any open intervals in $[0,1]$ of width at least $2^{-k}$. (Note $\nu(\calT_{k-1})$ cannot exclude more than $2^k$ disjoint open intervals in $[0,1]$ of width at least $2^{-k}$, so at most $2^k$ mid-point sets are added to $\calT_{k-1}$ to form $\calT_k$.) Thus the chain $\calU := \bigcup_{k=1}^{\infty} \calT_k$ contains $\calT$ and $\nu(\calU)$ is dense in $[0,1]$.

($2 \implies 3$) Suppose without loss of generality that $\emptyset, \N \in \calU$ (if not, simply add them). Fix $\epsilon > 0$. Then there exists finite $F \subseteq \nu(\calU)$ such that for every $x \in [0,1]$ there are $b, c \in F$ with $b \leq x \leq c$ and $c - b < \epsilon/2$, since $\nu(\calU)$ us dense in $[0,1]$. For each $b \in F$, there exists $A_b \in \calU$ such that $\nu(A_b) = b$. 

Since $F$ is finite, there exists $N_{\epsilon} \in \N$ such that $\lvert \nu_N(A_b) - \nu(A_b) \rvert < \epsilon/2$ for all $N > N_{\epsilon}$ and all $b \in F$. Now for any $A \in \calU$, there exist $b, c \in F$ with $A_b \subseteq A \subseteq A_c$, $b \leq \nu(A) \leq c$ and $c - b < \epsilon/2$. Thus,
\[
\lvert \nu_N(A_b) - \nu(A) \rvert \leq \lvert \nu_N(A_b) - b \rvert + \lvert \nu(A) - b \rvert \leq \lvert \nu_N(A_b) - \nu(A_b) \rvert + \lvert c - b \rvert < \epsilon
\]
and similarly $\lvert \nu_N(A_c) - \nu(A) \rvert < \epsilon$. 

Since $A_b \subseteq A \subseteq A_c$,  
\[
\nu(A) - \epsilon < \nu_N(A_b) \leq \nu_N(A) \leq \nu_N(A_c) < \nu(A) + \epsilon,
\] 
which implies $\lvert \nu_N(A) - \nu(A) \rvert < \epsilon$. Hence Condition~3 holds for all $A \in \calU$, and thus for all $A \in \calT$.

($3 \implies 1$) Consider $A \in \calT_{\cup}$ and fix $\epsilon > 0$. Define $\nu_{\cup}(A) := \sup \{ \nu(C) : C \in \calT, C \subseteq A \}$. Then there exists $C_1 \in \calT$ such that $C_1 \subseteq A$ and $\lvert \nu(C) - \nu_{\cup}(A) \rvert < \epsilon/2$ for all $C \in \calT$ such that $C_1 \subseteq C \subseteq A$. By assumption, there exists $N_{\epsilon} \in \N$ such that $\lvert \nu_N(C) - \nu(C) \rvert < \epsilon/2$ for all $C \in \calT$ and for all $N > N_{\epsilon}$. For any $N > N_{\epsilon}$, there exists $C_2 \in \calT$ such that $C_1 \subseteq C_2 \subseteq A$ and $I_A(n) = I_{C_2}(n)$ for $n = 1, \ldots, N$, so that $\lvert \nu_N(A) - \nu_N(C_2) \rvert = 0$. Thus
\[
\lvert \nu_N(A) - \nu_{\cup}(A) \rvert \leq \lvert \nu_N(A) - \nu_N(C_2) \rvert + \lvert \nu_N(C_2) - \nu(C_2) \rvert + \lvert \nu(C_2) - \nu_{\cup}(A) \rvert < \epsilon.
\]
Hence $\nu_N(A) \rightarrow \nu_{\cup}(A)$, implying $A \in \calF$ with $\nu(A) = \nu_{\cup}(A)$. Similarly, for all $A \in \calT_{\cap}$, $A \in \calF$ with $\nu(A) = \nu_{\cap}(A) :=  \inf \{ \nu(C) : C \in \calT, A \subseteq C \}$. Hence $\calT_{*} = \calT_{\cup} \cup \calT_{\cap} \subset \calF$. 

Now consider $\calC \subseteq \calT$ and define $A := \bigcup \calC$. Then $A \in \calT_{\cup} \subset \calF$ and $\nu(A) = \nu_{\cup}(A)$ as shown in the preceding paragraph. It is straightforward to check $\nu_{\cup}(A) = \sup_{C \in \calC} \nu(C)$. Similarly, $\nu(\bigcap \calC) = \inf_{C \in \calC} \nu(C)$. 

If any of the three statements hold for $\calT$, then Statement~2 also holds with $\calT$ replaced by $\calU$, since trivially $\calU \subseteq \calU$. Thus $\calU_* \subset \calF$, by Statement~1. But then Statement~2 holds with both $\calT$ and $\calU$ replaced by $\calU_*$. Hence Statements~1 and~3 hold with $\calT$ replaced by $\calU_*$. Finally, $\nu(\calU_*) = [0,1]$, since for any $x \in [0,1]$, the set $A_x := \bigcup \{ A \in \calU : \nu(A) \leq x \} \in \calU_*$ with $\nu(A_x) = x$. 
\qed

Theorem~\ref{uniformconvergence} is used in the proof of Theorem~\ref{nullmodification_fn} (more specifically, it is used in the proof of Lemma~\ref{Kp_Cesaro_integral}, on which the proof of Theorem~\ref{nullmodification_fn} depends).

The following two corollaries respectively provide a simplification of Theorem~\ref{uniformconvergence} for chains that are also sequences (Corollary~\ref{uniformconvergencesequence}), and a fourth characterisation of the class of chains described in Theorem~\ref{uniformconvergence} in terms of maximal chains (Corollary~\ref{maximal_chains}).

\begin{corollary} \label{uniformconvergencesequence}
Consider pairwise disjoint sets $\{ A_k \}_{k=1}^{\infty}$ in $\calF$. Let $B_k = \cup_{i=1}^k A_i$ for each $k$ and let $B = \cup_{i=1}^{\infty} A_i$. The following conditions are logically equivalent:
\begin{enumerate}
\item $B \in \calF$ and $\nu(B) = \sum_{k=1}^{\infty} \nu(A_k)$.
\item For any $\epsilon > 0$ there exists a positive integer $N_{\epsilon} = N_{\epsilon}(B_1, B_2, \ldots)$ such that 
\[
\lvert \nu_N(B_k) - \nu(B_k) \rvert < \epsilon
\]
for all $N \geq N_{\epsilon}$ and for all $k$. 
\end{enumerate}
Moreover, if either statement holds then $\lvert \nu_N(B) - \nu(B) \rvert < \epsilon$ for all $N \geq N_{\epsilon}$.
\end{corollary}

\noindent {\bf Proof: } 
Note $\calT = \{ B_k \}_{k=1}^{\infty}$ is a chain in $\calF$. Note also $\calT_{\cap} = \calT$ and $\calT_{*} = \calT_{\cup} = \calT \cup \{ B \}$. 

($1 \implies2$) Statement 1 gives $\calT_{*} = \calT \cup \{ B \} \subset \calF$. For any $\calC \subseteq \calT$, $\bigcap \calC$ is the smallest element of $\calC$, hence $\nu(\bigcap \calC) = \inf_{C \in \calC} \nu(C)$. If there is a largest element of $\calC$, then $\bigcup \calC$ is that largest element, otherwise $\bigcup \calC = B$. In the case of a largest element, $\nu(\bigcup \calC) = \sup_{C \in \calC} \nu(C)$. In the case $\bigcup \calC = B$, $\nu(\bigcup \calC) = \nu(B) =  \sum_{k=1}^{\infty} \nu(A_k) = \sup_{C \in \calC} \nu(C)$ again. Hence the uniform convergence condition holds on all $\calT_{*}$ by Theorem~\ref{uniformconvergence}.

($2 \implies 1$) Statement 2 is the uniform convergence condition of Theorem~\ref{uniformconvergence} as it applies to $\calT$. Hence $B \in \calT_{*} \subset \calF$ and $\nu(B) = \sup_{k=1}^{\infty} \nu(B_k) = \sum_{k=1}^{\infty} \nu(A_k)$.
\qed

\begin{corollary} \label{maximal_chains}
Let $\calT \subset \calF$ be a chain of sets. Then $\calT$ satisfies the equivalent conditions of Theorem~\ref{uniformconvergence} if and only if there exists a {\em maximal} chain $\calU_{**} \subset \calF$ (maximal in the sense that it is not a proper subset of any other chain in $\calP(\N)$) such that $\calT \subseteq \calU_{**}$ and $\nu(\calU_{**}) = [0,1]$.
\end{corollary}

\noindent {\bf Proof: } 
$(\implies)$
First define the chain $\calU_*$ as described in Theorem \ref{uniformconvergence}, and assume without loss of generality that $\emptyset, \N \in \calU_*$. For every $k \in \N$, let $B_k := \bigcup \{ A \in \calU_* : k \notin A \}$ and $C_k := \bigcap \{ A \in \calU_* : k \in A \}$. Let $D_k = C_k \setminus B_k$ (it is straightforward to verify $B_k \subset C_k$). Note $C_k$ is the smallest set in $\calU_*$ containing $k$, and $B_k$ is the largest set in $\calU_*$ not containing $k$. These sets have the following properties, which are left to the reader to verify.

\begin{enumerate}
    \item If $k' \in D_k$, then $D_k = D_{k'}$.
    \item If $B_k \subseteq A \subseteq C_k$ for some $A \in \calU_*$, $k \in \N$, then either $A = B_k$ or $A=C_k$.
    \item If $A \in \calU_*$, then for any $k \in \N$ either $A \subseteq B_k$ or $C_k \subseteq A$.
\end{enumerate}


For each $k \in \N$, let $D_k = \{x_{k1}, x_{k2}, \ldots \}$, ordered by increasing magnitude; this set may be finite or infinite. Now let $\calU_{**}$ be $\calU_*$ together with all sets of the form $B_k \cup \{x_{k1}, \ldots, x_{kN}\}$, for any $k \in \N$ and $N \in \N$ (if $D_k$ is finite, restrict $N$ accordingly). Then $\calU_{**}$ is a maximal chain, shown as follows. Let $E, F \in \calU_{**}$. If $E, F \in \calU_{*}$, they are comparable since $\calU_*$ is a chain. If $E \in \calU_*, F \in \calU_{**} \setminus \calU_{*}$, then choose $k \in \N$ so that $F$ is of the form $B_k \cup \{x_{k1}, \ldots, x_{kN}\}$. By $3$ above either $E \subseteq B_k$, in which case $E \subseteq F$, or $C_k \subseteq E$, in which case $F \subseteq E$, so in either case $E$ and $F$ are comparable. If $E,F \in \calU_{**} \setminus \calU_{*}$, then they must be of the form $E = B_k \cup \{x_{k1}, \ldots, x_{kN}\}$ and $F=B_{k'} \cup \{x_{k'1}, \ldots, x_{k'N'}\}$; if $k' \in D_k$ then one must contain the other by 1 above, whereas if $k' \notin D_k$ the result follows by noting that in this case either $C_k \subseteq B_{k'}$ or $C_{k'} \subseteq B_{k}$ (by 2 above). Thus $\calU_{**}$ is a chain. 

Suppose there exists a chain $\calV \subseteq \calP(\N)$ with $\calU_{**} \subseteq \calV$, and let $E \in \calV$. Suppose $C_k \subseteq E$ for all $k \in E$. Then $\cup_{k \in E} C_k \subseteq E$, implying $E = \cup_{k \in E} C_k \in \calU_*$. Alternatively, suppose there exists $k \in E$ such that $E \subseteq C_k$. Then $B_k \subseteq E \subseteq C_k$, since $k \notin B_k$. Either $E \in \{ B_k, C_k \} \subset \calU_*$, or there is a largest $N$ such that $x_{kN} \in E$, in which case $E = B_k \cup \{x_{k1}, \ldots, x_{kN}\} \in \calU_{**}$. Hence $\calU_{**}$ is maximal.


Finally, $\calU_{**} \subseteq \calF$, since $\calU_{*} \subseteq \calF$, and every set in $\calU_{**}$ differs from a set in $\calU_{*}$ by at most a finite (and therefore null) set.

$(\impliedby)$ This is immediate from Statement~2 of Theorem \ref{uniformconvergence}.
\qed

\section{Boolean algebras, quotients and the monotone class theorem}
\label{Boolean_review}

The set $\calF$ can in a certain sense be factored by the null sets $\calN$ to produce a simple structure known as a monotone class, on which the induced function $\nu$ is countably additive. This useful result is  Corollary~\ref{countableadditivity} below. The proof involves a technique for manipulating chains in $\calP(\N)$ that is here called {\em null modification}, described in Section~\ref{null_modification_section}. Both  sections involve Boolean quotients, and while the theory of Boolean algebras and their quotients will be familiar to many readers, it may nevertheless be helpful to briefly review key definitions and results. That is the purpose of this section. There are no new results in this section, but it does contain a slight generalisation of the monotone class theorem for Boolean algebras (Theorem~\ref{monotone_for_Boolean}), based on the proof for fields of sets given in Paul Halmos' classic text on Measure Theory.

A {\em Boolean algebra} \cite{givant2009} is an abstraction of a field of sets consisting of a non-empty set $\calA$ equipped with two binary operators called {\em join} $\vee$ and {\em meet} $\wedge$, a unary {\em complement} operator $^{\prime}$ and containing special elements called the {\em zero} $0$ and {\em unit} (or {\em one}) $1$, satisfying the following axioms:
\begin{eqnarray*}
p \wedge 1 = p, & p \vee 0 = p, \\
p \wedge p^{\prime} = 0, & p \vee p^{\prime} = 1, \\
p \wedge q = q \wedge p, & p \vee q = q \vee p, \\
p \wedge (q \vee r) = (p \wedge q) \vee (p \wedge r), & p \vee (q \wedge r) = (p \vee q) \wedge (p \vee r).
\end{eqnarray*} 
These four pairs of axioms are known as the {\em identity laws}, {\em complement laws}, {\em commutative laws} and {\em distributive laws} respectively, and entail a number of other well known identities including associative laws and De Morgan's laws. Other common Boolean operators and relations can be composed from the meet, join and complement, for example $p - q := p \wedge q^{\prime}$ and $p + q := (p \wedge q^{\prime}) \vee (p^{\prime} \wedge q)$.  Another example is the partial order defined by $p \leq q \iff p \vee q = q$.  

The simplest example of a Boolean algebra is the set $\{0, 1\}$, with basic Boolean operations defined by 
\begin{eqnarray*}
0 \wedge 0 = 0, & 1 \wedge 1 = 1, & 0 \wedge 1 = 1 \wedge 0 = 0,\\
0 \vee 0 = 0, & 1 \vee 1 = 1, & 0 \vee 1 = 1 \vee 0 = 1,\\
0^{\prime} = 1, & 1^{\prime} = 0.
\end{eqnarray*}

Any field of sets $\calA \subseteq \calP(X)$ on an arbitrary set $X$ is a Boolean algebra with pairwise intersection $\cap$ as the {\em meet} operator, pairwise union $\cup$ as the {\em join} operator, set complement $^c$ as the Boolean complement operator, the empty set $\emptyset$ as the zero and $X$ as the unit. Note also $p - q$ is the set difference $p \setminus q$, $p + q$ is the symmetric difference $p \triangle q$ and the partial order $p \leq q$ is the subset relation $p \subseteq q$.

A {\em Boolean homomorphism} is a mapping $f : \calA \rightarrow \calA^{\prime}$ between Boolean algebras $\calA$ and $\calA^{\prime}$ that respects the basic set operations. Specifically, a homomorphism satisfies
\begin{align*}
f(A \vee B) &= f(A) \vee f(B), \\
f(A \wedge B) &= f(A) \wedge f(B), \\
f(0) &= 0, \mbox{ and } \\
f(1) &= 1, 
\end{align*}
for all $A, B \in \calA$. It follows that $f(A^{\prime}) = f(A)^{\prime}$, and indeed all finite combinations of basic Boolean operations are respected, including the partial order, that is $p \leq q \implies f(p) \leq f(q)$. A {\em Boolean isomorphism} is a homomorphism with an inverse homomorphism.  

A {\em Boolean ideal} $\calM$~\cite{givant2009} is a non-empty subset of a Boolean algebra $\calA$ satisfying the following axioms: 
\begin{eqnarray*}
p,q \in \calM & \implies & p \vee q \in \calM, \\
p \in \calM, q \in \calA & \implies & p \wedge q \in \calM. 
\end{eqnarray*} 
For example, for any charge space $(X,\calA,\mu)$, the set $\mu^{-1}(0) := \{ A \in \calA : \mu(A) = 0 \}$, called the {\em kernel} of $\mu$, is an ideal of $\calA$, and the set $\{ A \in \calP(X) : \mu^*(A) = 0 \}$, where $\mu^*(A) := \inf \{ \mu(B) : B \in \calA, A \subseteq B \}$, forms an ideal of $\calP(X)$. The set $\calN$ defined in Definition~\ref{null_sets} is a Boolean ideal of $\calP(\N)$, since by Proposition~\ref{calFproperties}(6), $A \cup B \in \calN$ for all $A,B \in \calN$, and $A \cap C \in \calN$ for all $A \in \calN$ and $C \in \calP(\N)$. 

A Boolean ideal $\calM$ induces an equivalence relation $\sim$ on the containing Boolean algebra $\calA$ such that
\[
p \sim q \iff p+q \in \calM.
\]
The collection of equivalence classes $\calA / \calM := \{ [ p ] : p \in \calA \}$, where $[ p ]$ denotes the equivalence class of $p$ under the equivalence relation induced by $\calM$, is called the {\em quotient} of $\calA$ by $\calM$. When the Boolean algebra in question is ambiguous it is convenient to write $[ p ]_{\calA/\calM}$ to identify both the underlying algebra $\calA$ and the ideal $\calM$.  

A key example in this paper is the Boolean quotient
\[
\calP(\N) / \calN := \{ [A]  : A \in \calP(\N) \},
\] 
where $[A] $ denotes the equivalence class of $A$ under the equivalence relation $A \sim B \iff A \triangle B \in \calN$. 

A quotient is itself a Boolean algebra when equipped with the Boolean operators $[ p ] \wedge [ q ] := [ p \wedge q ]$, $[ p ] \vee [ q ] := [ p \vee q ]$, $[ p ]^{\prime} := [ p^{\prime} ]$, and with $[ 0 ]$ and $[ 1 ]$ as the zero and unit respectively. The map $p \mapsto [p]$ is a Boolean homomorphism. This map respects the partial order, and in fact $[ p ] \leq [ q ]$ if and only if there exists $p^* \in [ p ]$ such that $p^* \leq q$, or equivalently there exists $q^* \in [ q ]$ such that $p \leq q^*$.

For any $\calA^{\prime} \subseteq \calA$ (not necessarily a sub-algebra), define $[\calA^{\prime}] := \{ [A] : A \in \calA^{\prime} \}$. If $(X,\calA,\mu)$ is a charge space and $\calM$ is the kernel of $\mu$, the induced function $\mu : \calA / \calM \rightarrow \mathbb{R}$ given by $\mu[ A ] := \mu(A)$ for all $A \in \calA$ is finitely additive. (Here and throughout the paper, parentheses delimiting a function argument are omitted when the argument is contained in square brackets.)

A new version of the monotone class theorem, which generalises the version in \cite{keisler1977monotone}, is presented below. The new result makes use of the following definitions, some of which are non-standard. 

A Boolean algebra $\calA$ is said to be {\it countably complete} if every countable subset $\{ p_k \}_{k=1}^{\infty}$ in $\calA$ has a least upper bound in $\calA$. A subalgebra $\calB \subseteq \calA$ will here be called {\it countably complete} if every countable subset $\{ p_k \}_{k=1}^{\infty}$ in $\calB$ has an upper bound in $\calB$ that is less than any other upper bound of this subset in $\calA$. This upper bound is called the {\it supremum} of the subset and denoted $\vee_{k=1}^{\infty} p_k$. In that case, it is straightforward to show (by taking complements) that every countable subset also has a lower bound in $\calB$ that is greater than any other lower bound of the subset in $\calA$, called the {\it infimum} of the subset and denoted $\wedge_{k=1}^{\infty} p_k$. By definition these two elements are unique. 

An important subtlety is that a proper subalgebra $\calB \subset \calA$ will not here be called countably complete if it is only true that every countable subset of $\calB$ has an upper bound in $\calB$ that is less than any other upper bound of that subset in $\calB$: it must be less than any other upper bound of that subset in $\calA$. The reason for this requirement is that, without it, the supremum of a countable subset $\calC$ of a subalgebra $\calB_1$ could differ from the supremum of $\calC$ when viewed as a subset of a distinct subalgebra $\calB_2$. This can occur even if $\calB_1$ and $\calB_2$ are both countably complete algebras when the containing algebra $\calA$ is ignored. Thus the requirement is needed to ensure the supremum of $\calC$ is uniquely defined across all countably complete subalgebras of $\calA$. 

A subset $\calM$ of $\calA$ will be called a {\it monotone class} if: 
\begin{enumerate}
\item for any non-decreasing sequence $p_1 \leq p_2 \leq \ldots$ in $\calM$, there is an upper bound in $\calM$ that is less than any other upper bound of the sequence in $\calA$, and
\item for any non-increasing sequence $p_1 \geq p_2 \geq \ldots$ in $\calM$, there is a lower bound in $\calM$ that is greater than any other lower bound of the sequence in $\calA$.
\end{enumerate}
Similarly to countably complete subalgebras, this least upper bound will be called the supremum of the sequence, denoted $\vee_{k=1}^{\infty} p_k$, and this greatest lower bound will be called the infimum of the sequence, denoted $\wedge_{k=1}^{\infty} p_k$. The same word of caution is necessary here as for countably complete subalgebras: $\vee_{k=1}^{\infty} p_k$ must be less than any other upper bound of the sequence in $\calA$, not just in $\calM$, and similarly $\wedge_{k=1}^{\infty} p_k$ must be greater than any other lower bound in $\calA$, not just in $\calM$.

In fact, the version of the monotone class theorem presented in \cite{keisler1977monotone} also requires countably complete subalgebras and monotone classes to be understood in this sense, though this is not explicitly stated. Note the version of the monotone class theorem presented in that paper differs from the one below in requiring the containing algebra $\calA$ to be countably complete.




As stated above, the version of the monotone class theorem below is adapted from \cite{keisler1977monotone}. There it is claimed that the result is proved in \cite{halmos2013measure}; however this may be an example of mathematical folklore, as the result in that reference applies only to fields of sets, which are less general than Boolean algebras. A full proof is therefore provided.

\begin{theorem} \label{monotone_for_Boolean}
Let $\calA$ be a Boolean algebra and let $\calM \subseteq \calA$ be a monotone class. Let $\calA_0 \subseteq \calM$ be a subalgebra of $\calA$, and define $\sigma(\calA_0) \subseteq \calM$ to be the smallest monotone class in $\calA$ that contains $\calA_0$. Then $\sigma(\calA_0)$ is also the smallest countably complete subalgebra of $\calA$ that contains $\calA_0$.
\end{theorem}

Proof of the monotone class theorem depends on the following lemma, which is analogous to \cite[Thm. A, p. 27]{halmos2013measure}.

\begin{lemma} \label{BsubA_is_monotone}
Suppose $\calA$ is a Boolean algebra, and $\calM$ is a subalgebra that is also a monotone class. Then $\calM$ is a countably complete subalgebra.
\end{lemma}

\noindent {\bf Proof: }
Let $\calC=\{p_1,p_2, \ldots \}$. Then, since $\calM$ is a Boolean algebra, the elements $p_1, p_1 \vee p_2, p_1 \vee p_2 \vee p_3, \ldots$ are also in $\calM$. These elements form an increasing sequence, hence this sequence has an upper bound in $\calM$ that is less than any other upper bound in $\calA$. It may be checked this upper bound is also the supremum of $\calC$, implying $\calM$ is countably complete. 
\qed 

The following proof is derived from the proof of \cite[Thm. B, p. 27]{halmos2013measure}.

\noindent {\bf Proof of Theorem~\ref{monotone_for_Boolean}: }
It will be sufficient to show that $\sigma := \sigma(\calA_0)$ is a Boolean subalgebra, for then it will be countably complete by Lemma~\ref{BsubA_is_monotone}, and in fact it will be the smallest countably complete subalgebra containing $\calA_0$ because any smaller countably complete subalgebra containing $\calA_0$ would also be a smaller monotone class containing $\calA_0$. 

For $q \in \calA$, let
$$ 
K(q) = \{p \in \calA: p-q, q-p, q \vee p \in \sigma \}.
$$
These sets possess a convenient symmetry: $p \in K(q)$ if, and only if, $q \in K(p)$. Suppose $p_1 \leq p_2 \leq \ldots$ is a non-decreasing sequence of elements in $K(q)$. Then $p_k \vee q \in \sigma$ for all $k$, and since $\sigma$ is a monotone class it follows that 
$$
\vee_{k=1}^{\infty} (p_k \vee q) = \Big(\vee_{k=1}^{\infty} p_k\Big) \vee q \in \sigma.
$$

Similar arguments show that $\vee_{k=1}^{\infty} p_k - q, q-\vee_{k=1}^{\infty} p_k \in \sigma$ as well, and it follows that $\vee_{k=1}^{\infty} p_k \in K(q)$. A parallel argument shows that $\wedge_{k=1}^{\infty} p_k \in K(q)$, and thus $K(q)$ is a monotone class. If $q \in \calA_0$, then $\calA_0 \subseteq K(q)$, and thus, since $K(q)$ is a monotone class, $\sigma \subseteq K(q)$. However, the symmetry mentioned above now implies that if $p \in \sigma$ then $q \in K(p)$ for any $q \in \calA_0$, and then, since $K(p)$ is a monotone class, that $\sigma \subseteq K(p)$. This implies in particular that, for any $p,q \in \sigma$, the elements $p-q, q-p, q \vee p$ are all in $\sigma$, and hence $q \wedge p \in \sigma$ as well, since $q \wedge p = (q \vee p) - (p-q) - (q-p)$. It follows that $\sigma$ is a Boolean subalgebra.
\qed

\section{Null modification}
\label{null_modification_section}

This section develops another analytic tool for studying chains in $\calP(\N)$: a construction that is here called a {\em null modification}. A null modification takes a set in $\calP(\N)$ and constructs a new set of a form described in Property~9 of Proposition~\ref{calFproperties}, using Algorithm~1.

\begin{proposition} \label{nullmodification}
For any $A \in \calP(\N)$, Algorithm~1 decomposes $A$ into disjoint sets $A^{\prime} \in \calP(\N)$ and $F \in \calN$ such that
\begin{enumerate}
\item $A = A^{\prime} \cup F$,
\item $\nu^+(A^{\prime}) = \nu^+(A)$, and
\item $\nu_N(A^{\prime}) \leq \nu^+(A^{\prime})$ for all $N \in \calN$.
\end{enumerate}
Moreover, if $A \in \calF$, then $A^{\prime} \in \calF$.

\bigskip

\noindent {\bf Algorithm 1} Given $A \in \calP(\N)$, construct $A^{\prime} \in \calP(\N)$ and $F \in \calN$

\noindent Set $A^{\prime} = F = \emptyset$.\\
    {\bf for} $N=1, 2, \ldots$\\
      \indent {\bf if} $N \in A$ then\\
    	\indent \indent Add $N$ to $A^{\prime}$\\
      	\indent \indent {\bf if} $\nu_N(A^{\prime}) > \nu^+(A)$ then\\
	       \indent \indent \indent Remove $N$ from $A^{\prime}$ and add it to $F$.\\
	    \indent \indent {\bf endif}\\
        \indent {\bf endif}\\
    {\bf endfor}
\end{proposition}

\noindent {\bf Proof: } 
Algorithm~1 trivially ensures $A^{\prime}$ and $F$ are disjoint, $A = A^{\prime} \cup F$ and $\nu_N(A^{\prime}) \leq \nu^+(A)$ for all $N$. Moreover, if $N \in F$, then 
\[
\nu_N(A^{\prime}) > \nu^+(A) - \frac{1}{N}.
\]
Next show $F \in \calN$ as follows. This is trivial if $F$ is a finite set, so suppose it is infinite. Fix $\epsilon > 0$ and choose $N_{\epsilon}$ so that the following conditions are met:
\begin{enumerate}
\item $N_{\epsilon} \in F$,
\item $N_{\epsilon} > 2/ \epsilon$, and
\item $\nu_N(A) < \nu^+(A) + \epsilon /2$ for all $N \geq N_{\epsilon}$.
\end{enumerate}
Now for any $N \geq N_{\epsilon}$, if $N \in F$ then
\begin{eqnarray*}
\nu_N(F) & = & \nu_N(A) - \nu_N(A^{\prime}) \\
& < & \nu^+(A) + \frac{\epsilon}{2} - \left( \nu^+(A) - \frac{1}{N} \right) \\
& \leq & \frac{\epsilon}{2} + \frac{1}{N_{\epsilon}} \\
& < & \epsilon.
\end{eqnarray*}
If $N \notin F$ then
\[
\nu_N(F) < \nu_{N^{\prime}}(F) < \epsilon
\]
where $N^{\prime}$ is the largest integer less than $N$ for which $N^{\prime} \in F$, noting that $N^{\prime} \geq N_{\epsilon}$. Hence $\lim_{N \rightarrow \infty} \nu_N(F) = 0$, implying $F \in \calN$. Proposition~\ref{calFproperties}(6) gives $\nu^+(A^{\prime}) = \nu^+(A)$. If $A \in \calF$, Proposition~\ref{calFproperties}(6) gives $A^{\prime} \in \calF$ with $\nu(A^{\prime}) = \nu(A)$.  
\qed

Recalling the example given at the end of Section \ref{F_and_nu}, it can be checked that the set $\calO$ defined there, the set of all odd numbers at least 3, can be obtained by applying this algorithm to the set of all odd numbers. The null set removed by the algorithm is simply $\{1\}$.

Null modification can be used to transform chains in $\calF$ to acquire a useful topological property, defined in terms of the following pseudo-metric. Let $d_{\nu}(B,C) := \nu^+(B \triangle C)$ for all $B, C \in \calP(\N)$. Trivially, $d_{\nu}(B,B)=0$ and $d_{\nu}(B,C)=d_{\nu}(C,B)$. The triangle inequality $d_{\nu}(A,B) + d_{\nu}(B,C) \geq d_{\nu}(A,C)$ follows from the fact that
\[
\begin{split}
    A \triangle C & = (A \setminus C) \cup (C \setminus A) \subseteq (A \setminus B \cup B \setminus C) \cup (C \setminus B \cup B \setminus A) \\
    & = (A \setminus B \cup B \setminus A) \cup (C \setminus B \cup B \setminus C) = (A \triangle B) \cup (B \triangle C),
\end{split}
\]
so that $\nu_N(A \triangle C) \leq \nu_N(A \triangle B) + \nu_N(B \triangle C)$. 


This pseudo-metric is related to the continuity of the set functions $\nu^+$, $\nu^-$ and $\nu$ on chains in $\calP(\N)$ or $\calF$ in the following sense.

\begin{proposition} \label{sup_and_inf}
Consider a chain $\calS \subseteq \calP(\N)$. Then
\begin{enumerate}
\item if $\inf \{ d_{\nu}(\bigcup \calS, A) : A \in \calS \} = 0$, then $\nu^+(\bigcup \calS) = \sup \{ \nu^+(A) : A \in \calS \}$, and
\item if $\inf \{ d_{\nu}(A, \bigcap \calS) : A \in \calS \} = 0$, then $\nu^-(\bigcap \calS) = \inf \{ \nu^-(A) : A \in \calS \}$.
\end{enumerate}
Moreover, if $\calS \subseteq \calF$ and $\bigcup \calS \in \calF$, the converse of the first result holds, and if $\calS \subseteq \calF$ and $\bigcap \calS \in \calF$, the converse of the second result holds.
\end{proposition}

\noindent {\bf Proof: } 
For 1, first note $\sup \{ \nu^+(A) : A \in \calS \} \leq \nu^+(\bigcup \calS)$, since $\nu^+(A) \leq \nu^+(\bigcup \calS)$ for all $A \in \calS$. To show the reverse inequality, fix $\epsilon > 0$ and choose $A \in \calS$ such that $d_{\nu}(\bigcup \calS, A) < \epsilon$. By Proposition~\ref{calFproperties}(5), $d_{\nu}(\bigcup \calS, A) \geq \nu^+(\bigcup \calS) - \nu^+(A)$, hence
\[
\nu^+(\bigcup \calS) < \nu^+(A) + \epsilon \leq \sup \{ \nu^+(A) : A \in \calS \} + \epsilon.
\]
Let $\epsilon \rightarrow 0$ to obtain the first result. Then 2 follows by taking complements.

If $\calS \subseteq \calF$ and $\bigcup \calS \in \calF$, then $\bigcup \calS \setminus A \in \calF$ with $\nu(\bigcup \calS \setminus A) = \nu(\bigcup \calS) - \nu(A)$ for all $A \in \calS$, by Proposition~\ref{calFproperties}(5). Taking infima gives 
\[
\inf \{ d_{\nu}(\bigcup \calS, A) : A \in \calS \} = \nu(\bigcup \calS) - \sup \{ \nu(A) : A \in \calS \},
\]
implying the converse of 1. The converse of 2 follows by taking complements.
\qed 

The following theorem is the main result in this section, and it establishes that any chain $\calT$ in $\calF$ can be modified to produce a chain satisfying the conditions of Proposition~\ref{sup_and_inf} on all subchains $\calS \subseteq \calT$. This technique will later be used to show that any integrable function $f \in \L_p(\calA)$ can be modified to construct a sequence with a Ces\`aro limit $h$ that is equal to $f$ almost everywhere. To be specific, Theorem~\ref{nullmodification_thm} is invoked in the proof of Lemma~\ref{nullmodification_posfn}, which in turn is key to the proof of Theorem~\ref{nullmodification_fn} below. 

Theorem~\ref{nullmodification_thm} refers to equivalence classes in $\calP(\N)$: these are elements of the Boolean quotient $\calP(\N) / \calN$. The subscripted asterisk appearing in the theorem reprises the notation introduced in the paragraph before Proposition~\ref{basicfacts}.

\begin{theorem} \label{nullmodification_thm}
Given a chain $\calT \subset \calF$, there exists a map $\phi : \calT_* \cap \calF \rightarrow \calF$ such that:
\begin{enumerate}
\item for all $A \in \calT_* \cap \calF$, $\phi(A) \triangle A \in \calN$, 
\item for all $A, B \in \calT_* \cap \calF$, 
\begin{align*}
[A] = [B] &\iff \phi(A) = \phi(B) \mbox{ and } \\
[A] < [B] &\iff \phi(A) \subset \phi(B),
\end{align*}
\item for all $\calS \subseteq \calT_* \cap \calF$, 
\begin{enumerate}
\item $\bigcup \phi(\calS) \in \calF$ with $\nu(\bigcup \phi(\calS)) = \sup \{ \nu(A) : A \in \phi(\calS) \}$, and
\item $\bigcap \phi(\calS) \in \calF$ with $\nu(\bigcap \phi(\calS)) = \inf \{ \nu(A) : A \in \phi(\calS) \}$,
\end{enumerate}
\item for all $A \in \phi(\calT_* \cap \calF)_*$ and $N \in \N$, $\nu_N(A) \leq \nu(A)$, 
\item if there exists a second map $\rho : \calT_* \cap \calF \rightarrow \calF$ such that $\rho(A) \triangle A \in \calN$ for all $A \in \calT_* \cap \calF$, then for all $\calS \subseteq \calT_* \cap \calF$, 
\begin{enumerate}
\item if $\bigcup \rho(\calS) \in \calF$ with $\nu(\bigcup \rho(\calS)) = \sup \{ \nu(A) : A \in \rho(\calS) \}$, then 
\[
\left( \bigcup \phi(\calS) \right) \triangle \left( \bigcup \rho(\calS) \right) \in \calN, \mbox{ and }
\]
\item if $\bigcap \rho(\calS) \in \calF$ with $\nu(\bigcap \rho(\calS)) = \inf \{ \nu(A) : A \in \rho(\calS) \}$, then 
\[
\left( \bigcap \phi(\calS) \right) \triangle \left( \bigcap \rho(\calS) \right) \in \calN.
\]
\end{enumerate}
\end{enumerate}
\end{theorem}


The construction of the map $\phi$ involves first modifying the chain so that the conditions of Proposition~\ref{sup_and_inf}(1) hold for all subsets of the chain. To describe this construction, it will be convenient to introduce the following notation. Define
\begin{eqnarray*}
\calC(A,B] & := & \{ C \in \calC : A \subset C \subseteq B \} 
\end{eqnarray*}
to represent {\em sub-intervals} of a chain $\calC \subset \calP(N)$. Here $A, B \in \calP(\N)$ but are not necessarily elements of $\calC$. 

It will also be convenient to define the {\em left end-points} of a chain $\calC \subset \calF$ to be sets of the form $\bigcap \calS$, where $\calS \subset \calC$, such that:
\begin{enumerate}
\item $\bigcap \calS \in \calF$, and
\item $\exists \epsilon > 0$ such that for any $A \in \calC$,  $d_{\nu}(A,\bigcap \calS) < \epsilon \implies A \supseteq \bigcap \calS$.
\end{enumerate}
In other words, there is a ``gap'' of width at least $\epsilon$ to the left of $\nu(\bigcap \calS)$ in $\nu(\calC)$. A chain $\calC$ can have at most countably many left endpoints because there can be at most countably many disjoint sub-intervals in the interval $[0,1]$, corresponding to these gaps.

The subscripted $\cup$ in Property~4 of the following lemma represents closure under unions, reprising the notation introduced in the paragraph before Proposition~\ref{basicfacts}.

\begin{lemma} \label{nullmodification_oneside}
Consider a countable chain $\calC \subset \calF$ such that $\nu(A) \neq \nu(B)$ for distinct $A, B \in \calC$. There exists a map $\psi : \calC \rightarrow \calF$ such that:
\begin{enumerate}
\item for all $A \in \calC$, $\psi(A) \subseteq A$ with $A \setminus \psi(A) \in \calN$, 
\item for all $A, B \in \calC$, $A \subseteq B \implies \psi(A) \subseteq \psi(B)$, 
\item for all $\calS \subseteq \calC$, 
\begin{enumerate}
\item $\bigcup \psi(\calS) \in \calF$ with $\nu(\bigcup \psi(\calS)) = \sup \{ \nu(A) : A \in \psi(\calS) \}$, 
\item if $\bigcap \calS \in \calF$ with $\nu(\bigcap \calS) = \inf \{ \nu(A) : A \in \calS \}$, then $\bigcap \psi(\calS) \in \calF$ with $\nu(\bigcap \psi(\calS)) = \inf \{ \nu(A) : A \in \psi(\calS) \}$, 
\end{enumerate}
\item for all $A \in \psi(\calC)_{\cup}$ and $N \in \N$, $\nu_N(A) \leq \nu(A)$, and
\item if there exists a second map $\rho : \calC \rightarrow \calF$ such that $\rho(A) \triangle A \in \calN$ for all $A \in \calC$, then for any $\calS \subseteq \calC$ with $\bigcup \rho(\calS) \in \calF$ and $\nu(\bigcup \rho(\calS)) = \sup \{ \nu(A) : A \in \rho(\calS) \}$,
\[
\left( \bigcup \psi(\calS) \right) \triangle \left( \bigcup \rho(\calS) \right) \in \calN.
\] 
\end{enumerate}
\end{lemma}

\noindent {\bf Proof: }
Without loss of generality, suppose $\calC$ contains its left end-points. Note the chain remains countable if its left end-points are added, and also retains the property $\nu(A) \neq \nu(B)$ for distinct $A, B \in \calC$. To see the latter claim, note that if a left end-point has the same charge as some $A \in \calC$, then in fact $A$ is already a left end-point of $\calC$, and no new set need be added. The added end-points and their images under $\psi$ can be discarded at the end of the construction, and the stated properties will be retained.

Arbitrarily order the elements of $\calC$ as a sequence $\{ C_k \}_{k=1}^{\infty}$. (This assumes $\calC$ is infinite, but the proof that follows also applies if $\calC$ is finite, with minimal modification.) For each $k \in \N$, define $B_k$ to be the largest set in the sequence $C_1, \ldots, C_{k-1}$ that is a proper subset of $C_k$, if such a set exists (that is, $B_k := \bigcup \{ C_j : j < k, C_j \subset C_k  \}$). Otherwise, set $B_k := \emptyset$. 

Set $\psi_0(C) = C$ for all $C \in \calC$, and assume inductively that $\psi_{k-1}(\calC) \subset \calF$. This is trivially true for $k=1$, since $\calC \subset \calF$. Sequentially define $\psi_k$ for each $k \in \N$ as follows:
\begin{enumerate}
\item Apply Algorithm~1 to decompose $A_k := \psi_{k-1}(C_k) \setminus \psi_{k-1}(B_k)$ into disjoint sets $A_k ^{\prime} \in \calF$ and $F_k \in \calN$. 
\item For each $C \in \calC$, define
\[
\psi_k(C) := \left\{
\begin{array}{ll}
\psi_{k-1}(C) \setminus F_k & \mbox{ if } C \in \calC(B_k, C_k] \\
\psi_{k-1}(C) & \mbox{ otherwise.}
\end{array}
\right.
\]
\end{enumerate}
Note $\psi_k(C)$ differs from $C$ by the removal of at most $k$ null sets for each $k \in \N$ and $C \in \calC$. Hence $\psi_k(\calC) \subset \calF$.

For each $k \in \N$, define $\psi(C_k) := \psi_k(C_k)$. Then $\psi(C_k) \in \calF$ and Property~1 holds.

Also note that for all $A, B \in \calC$ and $k \in \N$, $A \subseteq B \implies \psi_k(A) \subseteq \psi_k(B)$. Moreover, $\psi_j(C_k) = \psi_k(C_k)$ for all $j > k$. Hence Property~2 holds because if $C_j \subseteq C_k$, then 
\[
\psi(C_j) = \psi_{\max\{j,k\}} (C_j) \subseteq \psi_{\max\{j,k\}} (C_k) = \psi(C_k).
\]

Next note that for all $k \in \N$,
\[
\psi(C_k) = \bigcup \{ A_j^{\prime} : j \leq k, C_j \subseteq C_k \},
\] 
where the components of the union are disjoint. This is shown by induction. It is trivially true for $k = 1$, since $\psi(C_1) = A_1^{\prime}$. Given it is true for all $j < k$, then since $B_k \in \{ \emptyset, C_1, \ldots, C_{k-1} \}$, 
\[
\psi_{k-1}(B_k) = \bigcup \{ A_j^{\prime} : j \leq k-1, C_j \subseteq B_k \},
\]
where the union is disjoint. Moreover, 
\begin{eqnarray*}
\psi(C_k) &=& \psi_k(C_k) = \psi_{k-1}(C_k) \setminus F_k = \psi_{k-1}(B_k) \cup A_k^{\prime} \\ 
&=& \bigcup \{ A_j^{\prime} : j \leq k, C_j \subseteq C_k \}
\end{eqnarray*}
with $A_k^{\prime}$ disjoint from $\psi_{k-1}(B_k)$.

For all $k, N \in \N$, $\nu_N(A_k^{\prime}) \leq \nu(A_k^{\prime})$ by Lemma~\ref{nullmodification}. Hence 
\[
\nu_N(\psi(C_k)) = \sum_{j \leq k: C_j \subseteq C_k} \nu_N(A_j^{\prime}) \leq \sum_{j \leq k: C_j \subseteq C_k} \nu(A_j^{\prime}) = \nu(\psi(C_k))
\]
by Proposition~\ref{calFproperties}(7). Property~3a thus follows by Proposition~\ref{calFproperties}(9).

To show 3b, first suppose $\bigcap \calS$ is a left end-point of $\calC$. Then $\bigcap \calS \in \calC$. Also note $\psi(\bigcap \calS) \subseteq \bigcap \psi(\calS)$, since $\psi(\bigcap \calS) \subseteq \psi(A)$ for all $A \in \calS$. Hence 
\begin{eqnarray*}
\nu(\bigcap \calS) &=& \nu(\psi(\bigcap \calS)) \\
&\leq& \nu^-(\bigcap\psi(\calS)) \\
&\leq& \nu^+(\bigcap\psi(\calS)) \\
&\leq& \inf \{ \nu(\psi(A)) : A \in \calS \} \\
&=& \inf \{ \nu(A) : A \in \calS \},
\end{eqnarray*}
implying $\bigcap\psi(\calS) \in \calF$ with $\nu(\bigcap\psi(\calS)) = \inf \{ \nu(A) : A \in \psi(\calS) \}$. On the other hand, if $\bigcap \calS$ is not a left end-point of $\calC$, then for any $\epsilon > 0$ there is $A \in \calC$ with $A \subset \bigcap \calS$ and $\nu(A) > \nu(\bigcap\calS) - \epsilon$. Hence 
\begin{eqnarray*}
\nu(\bigcap \calS) - \epsilon &<& \nu(A) \\
&=& \nu(\psi(A)) \\
&\leq& \nu^-(\bigcap\psi(\calS)) \\
&\leq& \nu^+(\bigcap\psi(\calS)) \\
&\leq& \inf \{ \nu(\psi(A)) : A \in \calS \} \\
&=& \inf \{ \nu(A) : A \in \calS \},
\end{eqnarray*}
again implying $\bigcap\psi(\calS) \in \calF$ with $\nu(\bigcap\psi(\calS)) = \inf \{ \nu(A) : A \in \psi(\calS) \}$. 

Consider any $\calS \subseteq \calC$, and define
\[
A := \bigcup \{ \psi(C) : C \in \calS \} = \cup_{i=1}^{\infty} \psi(C_{k_i})
\] 
where $k_1$ is the smallest positive integer $k$ such that $\psi(C_k) \subseteq A$, $k_2$ is the smallest positive integer $k$ such that $\psi(C_1) \subset \psi(C_k) \subseteq A$ and so on. This sequence is infinite for $A \not\in \calC$. It follows that $A = \cup_{j=1}^{\infty} A_{k_j}^{\prime}$,
where the union is disjoint, since for each $i \in \N$, $\psi(C_{k_i}) = \cup_{j=1}^{i} A_{k_j}^{\prime}$ is a disjoint union. Hence Property~4 holds because 
\[
\nu_N(A) = \sum_{j=1}^{\infty} \nu_N(A_{k_j}^{\prime}) \leq \sum_{j=1}^{\infty} \nu(A_{k_j}^{\prime}) \leq \nu(A)
\]
by Proposition~\ref{calFproperties}(8). 

For 5, define a third map $\tau(A) := \psi(A) \cap \rho(A)$ for all $A \in \calC$. Then since $\psi(A) \triangle A \in \calN$ and $\rho(A) \triangle A \in \calN$, one must have $\psi(A) \triangle \rho(A) \in \calN$ and then also $\tau(A) \triangle \psi(A) \in \calN$, using Proposition~\ref{calFproperties}(6). Moreover, $\rho(A) \in \calF$ and $\tau(A) \in \calF$ for $A \in \calC$, again by Proposition~\ref{calFproperties}(6). Now, for any $\calS \subseteq \calC$,
\begin{eqnarray*}
\nu^+ \left( \bigcup \psi(\calS) \setminus \bigcup \tau(\calS) \right) & \leq & \inf \left\{ \nu \left( \bigcup \psi(\calS) \setminus \tau(A) \right) : A \in \calS \right\} \\
& = & \nu(\bigcup \psi(\calS)) - \sup \{ \nu(\tau(A)) : A \in \calS \} \\
& = & \nu(\bigcup \psi(\calS)) - \sup \{ \nu(\psi(A)) : A \in \calS \} \\
& = & 0,
\end{eqnarray*}
using 3a. Noting $\bigcup \tau(\calS) \subseteq \bigcup \psi(\calS)$, this implies 
\[
\left(\bigcup \psi(\calS)\right) \triangle \left(\bigcup \tau(\calS)\right) \in \calN.
\]
If $\bigcup \rho(\calS) \in \calF$ with $\nu(\bigcup \rho(\calS)) = \sup \{ \nu(A) : A \in \rho(\calS) \}$, then one may apply a similar argument to $\rho$ instead of $\psi$, giving $\left(\bigcup \rho(\calS)\right) \triangle \left(\bigcup \tau(\calS)\right) \in \calN$, and hence 
\[
\left( \bigcup \psi(\calS) \right) \triangle \left( \bigcup \rho(\calS) \right) \in \calN.
\] 
\qed

One can now apply Lemma~\ref{nullmodification_oneside} twice - to a given chain and to the corresponding chain of complements - to obtain a chain that satisfies both the sufficient conditions of Proposition~\ref{sup_and_inf}(1) and~(2): this strategy will be used to prove Theorem~\ref{nullmodification_thm} below. The requirement that the chain be countable can also be removed, by an argument involving the following real analysis lemma.

\begin{lemma} \label{countable_subset} 
Totally ordered sets have the following properties.
\begin{enumerate}
\item Any $\calT \subseteq \mathbb{R}$ contains a countable subset $\calC$ such that for $a \in \calT$ and $\epsilon > 0$, there exist $b, c \in \calC$ with $b \leq a \leq c$ and $c - b < \epsilon$. 
\item Consider a totally ordered set $\calT$ and a strictly increasing function $\mu : \calT \rightarrow \mathbb{R}$. Then $\calT$ contains a countable subset $\calC$ such that for $A \in \calT$ and $\epsilon > 0$, there exist $B, C \in \calC$ with $B \leq A \leq C$ and $\mu(C) - \mu(B) < \epsilon$.
\end{enumerate}
\end{lemma}

\noindent {\bf Proof: }
For 1, let $\calC_0$ be a countable, dense subset of $\calT$. (Such a subset exists since any subset of the reals is separable. Standard proofs of this invoke the axiom of countable choice.) Define 
\begin{align*}
\calC_1 &:= \{ c \in \calC : \exists \epsilon > 0 \mbox{ with } (c - \epsilon/2,c) \cap \calC_0 = \emptyset \}, \mbox{ and } \\ 
\calC_2 &:= \{ c \in \calC : \exists \epsilon > 0 \mbox{ with } (c, c + \epsilon/2) \cap \calC_0 = \emptyset \}.
\end{align*}
Then $\calC_1$ and $\calC_2$ are both countable, since there cannot be an uncountable number of pairwise disjoint intervals of non-zero width contained in $\mathbb{R}$ (each must contain a distinct rational). Thus $\calC := \calC_0 \cup \calC_1 \cup \calC_2$ is a countable set with the claimed property.

Claim~2 then follows by applying Claim~1 to $\mu(\calT)$, and noting $\mu$ is one-to-one and order preserving, as a consequence of being strictly increasing.
\qed 

Theorem~\ref{nullmodification_thm} can now be proved as follows. 

\noindent {\bf Proof of Theorem~\ref{nullmodification_thm}: }
Without loss of generality, suppose $\calT = \calT_* \cap \calF$. (No generality is lost because for any chain $\calT \subset \calF$, the chain $\calS := \calT_* \cap \calF$ has the property $\calS = \calS_* \cap \calF$. Moreover, if the lemma holds for $\calS$, it holds for $\calT$.) 

Let $[\calC]$ be the countable subchain of $[\calT] := \{ [A] : A \in \calT \}$ obtained by applying Lemma~\ref{countable_subset}(2) to the totally ordered set $[\calT]$ and the strictly increasing function $\nu: [\calT] \rightarrow \mathbb{R}$ given by $\nu[A] := \nu(A)$. Then let $\calC \subseteq \calT$ be obtained by selecting exactly one element of $\calT$ from each of the equivalence classes in $[\calC]$. (This implicitly invokes the axiom of countable choice, in general.) 

Construct a map $\psi : \calC^c \rightarrow \calF$ as in the proof of Lemma~\ref{nullmodification_oneside}, and then a map $\psi^{\prime} : \psi(\calC^c)^c \rightarrow \calF$ also as in the proof of Lemma~\ref{nullmodification_oneside}. Define $\phi^{\prime}(A) := \psi^{\prime}(\psi(A^c)^c)$ for each $A \in \calC$. 

The map $\phi^{\prime}$ inherits Properties~2, 3a, and 4 of Lemma~\ref{nullmodification_oneside} from $\psi^{\prime}$ and $\psi$. Property~1 of that lemma must be weakened to $\phi^{\prime}(A) \triangle A$ for all $A \in \calC$, because $\phi^{\prime}$ effectively adds a null set to $A$ and then removes a null set. However, Property~3b of that lemma can be strengthened to the unconditional claim $\bigcap \phi^{\prime}(\calS) \in \calF$ with $\nu(\bigcap \phi^{\prime}(\calS)) = \inf \{ \nu(A) : A \in \phi^{\prime}(\calS) \}$ for all $\calS \subseteq \calC$. This follows because $\bigcup \psi(\calS^c) \in \calF$ with $\nu(\bigcup \psi(\calS^c)) = \sup \{ \nu(A): A \in \psi(\calS^c) \}$, by Lemma~\ref{nullmodification_oneside}(3a) as it applies to $\psi$. Hence the condition of Lemma~\ref{nullmodification_oneside}(3b), as it applies to $\psi^{\prime}$, is satisfied for any $\calS \subseteq \calC$, that is $\bigcap \psi(\calS^c)^c \in \calF$ with
\[
\nu(\bigcap \psi(\calS^c)^c) = \inf \{ \nu(A): A \in \psi(\calS^c)^c \}.
\]

Define
\[
\phi(A) := \bigcup \{ \phi^{\prime}(B) : B \in \calC, [B] \leq [A] \}
\]
for each $A \in \calT$. Then $\phi$ is an extension of $\phi^{\prime}$, because for any $A, B \in \calC$ with $[B] \leq [A]$, one must have $B \subseteq A$, since $\calC$ contains at most one element from each equivalence class. Hence $\phi(A) = \phi^{\prime}(A)$.

Note that for any $A, B \in \calT$,
\begin{align*}
[A] = [B] &\implies \phi(A) = \phi(B) \mbox{ and } \\
[A] < [B] &\implies \phi(A) \subseteq \phi(B).
\end{align*}
This is immediate from the definition of $\phi$, and entails that $\phi$ is order preserving. 

To show~1, consider $A \in \calT$ and $\epsilon > 0$. There exist $B, C \in \calC$ with $[B] \leq [A] \leq [C]$ and $\nu(C) - \nu(B) < \epsilon / 2$. If $[A] = [B]$ then $A \triangle B \in \calN$ and $\phi(A) = \phi^{\prime}(B)$. Moreover, $\phi^{\prime}(B) \triangle B \in \calN$. Hence $\phi(A) \triangle A \in \calN$. Similarly, if $[A] = [C]$ then $\phi(A) \triangle A \in \calN$. Suppose $[B] < [A] < [C]$, implying $B \subset A \subset C$. Then $\phi^{\prime}(B) \subseteq \phi(A) \subseteq \phi^{\prime}(C)$, since $\phi$ is order preserving, and $\nu(\phi^{\prime}(C)) - \nu(\phi^{\prime}(B)) < \epsilon / 2$, using Proposition~\ref{calFproperties}(6). It follows that
\[
\phi(A) \triangle A \subseteq (\phi^{\prime}(B) \triangle B) \cup (\phi^{\prime}(C) \setminus \phi^{\prime}(B)) \cup (C \setminus B),
\]
and
\[
\nu^+(\phi(A) \triangle A) \leq \nu(\phi(B) \triangle B) + \nu(\phi^{\prime}(C)) - \nu(\phi^{\prime}(B)) + \nu(C) - \nu(B) < \epsilon.
\]
Letting $\epsilon \rightarrow 0$ gives $\phi(A) \triangle A \in \calN$. This also implies $\phi(\calT) \subset \calF$, using Proposition~\ref{calFproperties}(6).

For 2, note
\[
\phi(A) = \phi(B) \implies [\phi(A)] = [\phi(B)] \implies [A] = [B],
\]
since $A \in [\phi(A)]$ and $B \in [\phi(B)]$. This in turn gives 
\[
[A] < [B] \implies \phi(A) \neq \phi(B) \implies \phi(A) \subset \phi(B)
\]
since it is already established that $[A] < [B] \implies \phi(A) \subseteq \phi(B)$. Moreover,
\begin{eqnarray*}
\phi(A) \subset \phi(B) &\implies& \phi(A) \neq \phi(B) \mbox{ and } \phi(B) \not\subset \phi(A) \\
&\implies& [A] \neq [B] \mbox{ and } [B] \not< [A] \\
&\implies& [A] < [B].
\end{eqnarray*}

For 3a, consider $\calS \subseteq \calT$. If $\bigcup \phi(\calS) = \phi(A)$ for any $A \in \calS$, then  the result holds trivially, so assume $\bigcup \phi(\calS) \supset \phi(A)$ for all $A \in \calS$. But then for any $A \in \calS$, there is $B \in \calS$ with $\phi(A) \subset \phi(B) \subset \bigcup \phi(\calS)$ and $C \in \calC$ with $\phi(B) \subseteq \phi^{\prime}(C) \subset \bigcup \phi(\calS)$. Define $\calB := \{ C \in \calC : \phi^{\prime}(C) \subset \bigcup \phi(\calS) \}$, then $\bigcup \phi(\calS) = \bigcup \phi^{\prime}(\calB) \in \calF$ and 
\begin{eqnarray*}
\nu \bigg( \bigcup \phi(\calS) \bigg) &=& \nu \bigg( \bigcup \phi^{\prime}(\calB) \bigg) \\
&=& \sup \{ \nu(\phi^{\prime}(C)) : C \in \calB \} \\
&=& \sup \{ \nu(\phi(A)) : A \in \calS \}.
\end{eqnarray*}
Property~3b follows by a similar argument.

For 4, consider $A \in \phi(\calT)_*$, $N \in \N$ and $\epsilon > 0$. Consider the case $A := \bigcup \phi(\calS)$ for some $\calS \subseteq \calT$. If $\nu_N(A) = 0$ then trivially $\nu_N(A) \leq \nu(A)$, so assume $A$ contains some element of $\{ 1, \ldots, N \}$. Then there is $B \in \calS$ such that $\nu_N(\phi(B)) = \nu_N(A)$, since each element of $A \cap \{ 1, \ldots, N \}$ must be contained in some member of $\phi(\calS)$. Moreover, there exists $C \in \calC$ such that $[B] \leq [C]$ and $\nu(C) - \nu(B) \leq \epsilon$. Putting this all together,
\[
\nu_N(A) = \nu_N(\phi(B)) \leq \nu_N(\phi(C)) \leq \nu(\phi(C)) \leq \nu(\phi(B)) + \epsilon \leq \nu(A) + \epsilon.
\]
Letting $\epsilon \rightarrow 0$ gives $\nu_N(A) \leq \nu(A)$.

On the other hand, if $A := \bigcap \phi(\calS)$ for some $\calS \subseteq \calT$, then there is $B \in \calS$ with $\nu(\phi(B)) - \nu(A) < \epsilon$, by Property 3b. Hence 
\begin{eqnarray*}
\nu_N(A) - \nu(A) &=& \big( \nu_N(A) - \nu_N(\phi(B)) \big) + \big( \nu_N(\phi(B)) - \nu(\phi(B)) \big) \\
&& + \big( \nu(\phi(B)) - \nu(A) \big) \\
&\leq& 0 + 0 + \epsilon, 
\end{eqnarray*}
where the first summand is less than or equal 0 because $A \subseteq \phi(B)$, and the second summand is less than or equal 0 by the first case. Letting $\epsilon \rightarrow 0$ gives $\nu_N(A) \leq \nu^+(A)$.

Property~5a follows by essentially the same argument as Lemma~\ref{nullmodification_oneside}(5), merely replacing $\psi$ with $\phi$. Property~5b also follows by a similar argument, but using Property~3b instead of 3a.
\qed

In the process of proving Theorem~\ref{nullmodification_thm}, the following corollary emerges as an additional benefit. In essence, it establishes that any countable, pairwise disjoint collection of sets in $\calF$ can be transformed into a similar collection on which $\nu$ is countably additive. This corollary will later be used in the proofs of Theorems~\ref{calF_quotient_complete} and~\ref{nuisameasure}. 

\begin{corollary} \label{nullmodification_sequence}
Consider pairwise disjoint sets $\{ A_k \}_{k=1}^{\infty}$ in $\calF$. For each $k$, there exists $A_k^{\prime} \in \calF$ such that
\begin{enumerate}
\item $A_k^{\prime} \subseteq A_k$ with $F_k := A_k \setminus A_k^{\prime} \in \calN$,  
\item $A := \cup_{k=j}^{\infty} A_k^{\prime} \in \calF$ with $\nu(A) = \sum_{k=1}^{\infty} \nu(A_k^{\prime})$,
\item $\nu_N(A_k^{\prime}) \leq \nu(A_k^{\prime})$ for all $k,N \in \calN$, and $\nu_N(A) \leq \nu(A)$ for all $N \in \calN$,
\item if there exists $\{ A_k^{\prime\prime} \}_{k=1}^{\infty} \subset \calF$ such that 
\begin{enumerate}
\item $A_k \triangle A_k^{\prime\prime} \in \calN$ for each $k \in \N$, and
\item $\cup_{k=j}^{\infty} A_k^{\prime\prime} \in \calF$ with $\nu(\cup_{k=j}^{\infty} A_k^{\prime\prime}) = \sum_{k=1}^{\infty} \nu(A_k^{\prime\prime})$,
\end{enumerate}
then $\left( \cup_{k=1}^{\infty} A_k^{\prime} \right) \triangle  \left( \cup_{k=1}^{\infty} A_k^{\prime\prime} \right) \in \calN$.
\end{enumerate}
\end{corollary}

\noindent {\bf Proof: } 
First consider the case $\nu(A_k) > 0$ for all $k \in \N$. Then one may apply Lemma~\ref{nullmodification_oneside} to the chain $\calC := \{ B_j \}_{j=1}^{\infty}$, where $B_j := \cup_{k=1}^j A_k$ for $j \in \N$, $A_1^{\prime} = \psi(B_1)$ and $A_{k+1}^{\prime} := \psi(B_{k+1}) \setminus \psi(B_k)$ for $k \in \N$. 

Property~1 follows because the construction in the proof of Lemma~\ref{nullmodification_oneside} uses Algorithm~1 to remove the null set $F_k$ from $A_k$ to produce $A_k^{\prime}$. Property~2 follows from Lemma~\ref{nullmodification_oneside}(3a) with $\calS = \calC$. Property~3 is immediate from Lemma~\ref{nullmodification_oneside}(4). Property~4 follows from Lemma~\ref{nullmodification_oneside}(5), with $\rho(B_j) := \cup_{k=1}^j A_k^{\prime\prime}$ for $j \in \N$. 

Now consider the case $\nu(A_k) = 0$ for some $k \in  \N$. For each such $k$, set $A_k^{\prime} := \emptyset$ and apply Lemma~\ref{nullmodification_oneside} as above to the remaining elements of $\{ A_k \}_{k=1}^{\infty}$ (ie. those with non-zero charge). Then Properties~1 and the first part of~3 hold for $A_k^{\prime} = \emptyset$, as for the sets with non-zero charge. Property~2 and the latter part of~3 are properties of $A := \cup_{k=j}^{\infty} A_k^{\prime}$, to which sets with $A_k^{\prime} = \emptyset$ make no contribution. Property~4 is also a property of $A$, and will therefore hold provided
\[
\nu \bigg( \bigcup \{ A_k^{\prime\prime} : \nu(A_k) = 0 \} \bigg) = 0.
\]
This must be the case, otherwise condition~4b could not hold.
\qed

\section{The space $[\calF]$}
\label{quotient_space}

This section reveals a useful property of the set $\calF$, specifically that the quotient map $\xi: \calP(\N) \rightarrow \calP(\N) / \calN$ maps $\calF$ to a monotone class of $\calP(\N) / \calN$. 

Recall the notation
\[
[\calF] := \{ [A]  \in \calP(\N) / \calN : A \in \calF \}.
\]
The set function $\nu: \calF \rightarrow [0,1]$ induces a corresponding function on $[\calF]$ defined by
\[
\nu[A] := \nu(A),
\]
This function is well defined, since if $[A]  = [B]$ and $A \in \calF$, then $A \triangle B \in \calN$ and $B \in \calF$ with $\nu(A) = \nu(B)$ by Properties~5 and~6 of Proposition~\ref{calFproperties}. Similarly the set functions $\nu^+: \calP(\N) \rightarrow [0,1]$ and $\nu^- : \calP(\N) \rightarrow [0,1]$ induce functions on $\calP(\N) / \calN$ defined by
\[
\nu^+[A] := \nu^+(A) \mbox{ and } \nu^-[A] := \nu^-(A).
\]
These set functions are well defined as a consequence of the following result.

\begin{proposition} \label{dance}
If $A \triangle B \in \calN$ then $\nu^+(A) = \nu^+(B)$ and $\nu^-(A) = \nu^-(B)$.
\end{proposition}

\noindent {\bf Proof: } 
For any $N \in \N$,
\[
\nu_N(A \cup B) = \nu_N(A) + \nu_N(B \setminus A) = \nu_N(B) + \nu_N(A \setminus B).
\]
Consequently,
\[
\lvert \nu_N(A) - \nu_N(B) \rvert = \lvert \nu_N(A \setminus B) - \nu_N(B \setminus A) \rvert \leq \nu_N(A \setminus B) + \nu_N(B \setminus A) = \nu_N(A \triangle B).
\]
But then
\[
\nu_N(A) \leq \nu_N(B) + \lvert \nu_N(A) - \nu_N(B) \rvert \leq \nu_N(B) + \nu_N(A \triangle B).
\]
Taking the $\limsup$ as $N \rightarrow \infty$ gives
\[
\nu^+(A) \leq \nu^+(B) + \nu(A \triangle B) = \nu^+(B).
\]
Similarly $\nu^+(B) \leq \nu^+(A)$ and thus $\nu^+(A) = \nu^+(B)$. To show the corresponding result for $\nu^-$, note that $A \triangle B \in \calN$ implies $A^c \triangle B^c \in \calN$, and so
\[
\nu^-(A) = 1 - \nu^+(A^c) = 1 - \nu^+(B^c) = \nu^-(B).
\]
\qed 

Note $\nu^+[A] = \nu^-[A]$ if and only if $[A]  \in [\calF]$, and if either statement holds then $\nu[A] = \nu^+[A]$. These properties are inherited from the corresponding set functions on $\calP(\N)$. 

As Proposition~\ref{dance} and the preceding discussion suggest, many of the properties of $\calF$ carry over naturally to $[\calF]$, in some cases with simplified or stronger statements. These are summarized in Proposition~\ref{quotient_properties} below, the proof of which is straightforward and omitted. Note in particular that the final claim of Statement~5 below differs from the corresponding claim in Proposition~\ref{calFproperties}, because $\nu[A] = 0 \iff [A] = [\emptyset]$ for $A \in \calF$.

Property~7 below refers to disjoint elements in $\calP(\N) / \calN$. This is conventional terminology for elements of an abstract Boolean algebra that have a meet of 0, but nevertheless the following definition may clarify the meaning of `disjoint' in the present context.

\begin{definition}
Equivalence classes $[A], [B] \in \calP(\N) / \calN$ are said to be {\em disjoint} if $[A] \wedge [B] = [\emptyset]$, or equivalently $A \cap B \in \calN$.
\end{definition}
Note the two definitions are equivalent because $[A] \wedge [B] = [A \cap B]$, and the latter is equal to $[\emptyset]$ if and only if $A \cap B \in \calN$.

\begin{proposition} \label{quotient_properties}
The collection $[\calF]$ and functions $\nu^+$, $\nu^-$ and $\nu$ have the following properties.
\begin{enumerate}
\item $[\emptyset], [\N] \in [\calF]$, with $\nu[\emptyset] = 0$ and $\nu[\N ] = 1$.
\item For all $[A] \in [\calF]$, $[A]^c \in [\calF]$ with $\nu([A]^c) = 1 - \nu[A]$.
\item For all $[A], [B] \in [\calF]$, $[A] \vee [B] \in [\calF]$ if and only if $[A] \wedge [B] \in [\calF]$, and if either is true then $\nu([A] \vee [B]) = \nu[A] + \nu[B] - \nu([A] \wedge [B])$.
\item If $[A], [B] \in \calP(\N) / \calN$, then $\nu^+([A] \vee [B]) \leq \nu^+[A] + \nu^+[B]$. If $[A]$, $[B]$, and $[A] \vee [B]$ are all in $\calF$, then $\nu([A] \vee [B]) \leq \nu[A] + \nu[B]$.
\item For $[A], [B] \in \calP(\N) / \calN$ such that $[A] \leq [B]$, 
\begin{enumerate}
\item $\nu^+[A] \leq \nu^+[B]$,  
\item $\nu^+([B] - [A]) \geq \nu^+[B] - \nu^+[A]$, 
\item $\nu^-[A] \leq \nu^-[B]$, and 
\item $\nu^-([B] - [A]) \leq \nu^-[B] - \nu^-[A]$. 
\end{enumerate}
If in addition $[A], [B] \in [\calF]$, then 
\begin{enumerate}
\item $\nu[A] \leq \nu[B]$, 
\item $[B] - [A] \in [\calF]$,
\item $\nu([B] - [A]) = \nu[B] - \nu[A]$, and
\item $\nu[A] = \nu[B] \iff [A] = [B]$.
\end{enumerate}
\item For $[A] \in \calP(\N) / \calN$, $\nu^+[A] = 0 \iff [A] = [\emptyset]$. 
\item For pairwise disjoint elements $[A_1], \ldots, [A_K] \in [\calF]$, $\vee_{k=1}^K [A_k] \in [\calF]$ with 
\[
\nu(\vee_{k=1}^K [A_k]) = \sum_{k=1}^K \nu[A_k].
\] 
\end{enumerate}
\end{proposition}

The properties related to countable additivity - Properties 8 to 10 of Proposition~\ref{calFproperties} - are not carried over to this setting because Theorem~\ref{calF_quotient_complete} below establishes something stronger: that $\nu$ is countably additive on pairwise disjoint sequences in $[\calF]$.

\begin{theorem} \label{calF_quotient_complete}
For pairwise disjoint equivalence classes $\{ [ A_k ]  \}_{k=1}^{\infty} \subset [\calF]$, the representative sets $\{ A_k \}_{k=1}^{\infty} \subset \calF$ can be chosen so that
\begin{enumerate}
\item $\{ A_k \}_{k=1}^{\infty}$ are pairwise disjoint, and
\item $A := \cup_{k=1}^{\infty} A_k \in \calF$ with $\nu(A) = \sum_{k=1}^{\infty} \nu(A_k)$.
\end{enumerate}
Moreover, $[ A ]$ is the least upper bound of $\{ [A_k] \}_{k=1}^{\infty}$ in $\calP(\N) / \calN$. 
\end{theorem}

\noindent {\bf Proof: } 
Consider pairwise disjoint $\{ [ A_k ]  \}_{k=1}^{\infty} \subset [\calF]$. Then $A_j \cap A_k \in \calN$ for $j, k \in \N$ with $j \neq k$, hence $A_k \cap (\cup_{j=1}^{k-1} A_j) \in \calN$ and $[ A_k ]  = [ A_k \setminus (\cup_{j=1}^{k-1} A_j) ]$. Thus one may instead choose $\{ A_k \setminus (\cup_{j=1}^{k-1} A_j) \}_{k=1}^{\infty}$ to be the representative sets from each equivalence class to ensure Claim~1 holds. One may then use Corollary~\ref{nullmodification_sequence} to replace each representative set with a subset in the same equivalence class to ensure Claim~2 holds. 

To see that $[ A ]$ is the least upper bound  of $\{ [A_k]  \}_{k=1}^{\infty}$ in $\calP(\N) / \calN$, first note $[ A ]$ is an upper bound for this set because $A_k \subseteq A$ for all $k$, and the quotient map respects the partial order. Suppose there is some other $B \in \calP(\N)$ such that $[ A_k ]  \leq [ B ]$ for all $k$. Define $A_k^{\prime} := B \cap A_k$ for each $k$, and $A^{\prime} := \cup_{k=1}^{\infty} A_k^{\prime}$. Then $[ A^{\prime} ]  = [ B \cap A ]  = [ B ]  \wedge [ A ]$, giving $[ A^{\prime} ] \leq [B]$. Moreover, since both $[A]$ and $[B]$ are upper bounds for $[A_k]$, so is $[ A^{\prime} ]$, which implies $A_k \setminus A^{\prime} \in \calN$ for all $k \in \N$. Now, for any $K \in \N$,
\begin{eqnarray*}
\nu^+(A \setminus A^{\prime}) &\leq& \nu^+(A \setminus \cup_{k=1}^K A_k) + \nu^+( \cup_{k=1}^K A_k \setminus A^{\prime}) \\
&=& \nu(\cup_{k=K+1}^{\infty} A_k) \\
&=& \lim_{K \rightarrow \infty} \sum_{k = K+1}^{\infty} \nu(A_k).
\end{eqnarray*}
Letting $K \rightarrow \infty$ gives $\nu^+(A \setminus A^{\prime}) = 0$ and hence $[A] \leq [A^{\prime}] \leq [B]$.
\qed 

This theorem has the following useful corollary.

\begin{corollary} \label{countableadditivity}
$[\calF]$ is a monotone class in $\calP(\N) / \calN$ on which $\nu$ is countably additive.
\end{corollary}

\noindent {\bf Proof: } 
If $B_1 \leq B_2 \leq \ldots$ is a non-decreasing sequence in $[\calF]$, Theorem~\ref{calF_quotient_complete} implies the sequence has a supremum $B \in [\calF]$. Similarly, if $C_1 \geq C_2 \geq \ldots$ is a non-increasing sequence in $[\calF]$, then $C_1^{\prime} \leq C_2^{\prime} \leq \ldots$ is a non-decreasing sequence in $[\calF]$ with supremum $C^{\prime} \in [\calF]$, and then $C$ is the infimum for $\{ C_k \}_{k=1}^{\infty}$. Thus $[\calF]$ is a monotone class. Countable additivity of $\nu$ on $[\calF]$ follows from Theorem~\ref{calF_quotient_complete}(2).
\qed 

As a consequence of this corollary and the monotone class theorem for Boolean algebras (Theorem~\ref{monotone_for_Boolean}), $[\calF]$ contains the countably complete Boolean algebras generated by any Boolean algebra contained in $[\calF]$. This fact is used in the proofs of Theorem~\ref{isomorphism}, Theorem~\ref{chain_generated_field} and Corollary~\ref{Lp_Polish} to construct complete and separable $L_p$ spaces.

\section{Ces\`aro limits as integrals}
\label{cesaro_integrals}



Up to this point, the paper has focused on binary sequences with Ces\`aro limits, but this and the remaining sections study general real-valued sequences with Ces\`aro limits. First, the definitions of $\nu_N$, $\nu^+$, $\nu^-$ and $\nu$ can be generalised as follows.

\begin{definition}
For any function $f : \N \rightarrow \mathbb{R}$, define
\[
\nu_N(f) := \frac{1}{N} \sum_{n=1}^N f(n).
\]
Also define
\[
\nu^+(f) := \limsup_{N \rightarrow \infty} \nu_N(f),
\]
and
\[
\nu^-(f) := \liminf_{N \rightarrow \infty} \nu_N(f).
\]
If the limit as $N \rightarrow \infty$ exists, or equivalently if $\nu^+(f) = \nu^-(f)$, define
\[
\nu(f) := \lim_{N \rightarrow \infty} \nu_N(f) = \nu^+(f) = \nu^-(f).
\]
\end{definition}

Note $\nu_N(I_A) = \nu_N(A)$, $\nu^+(I_A) = \nu^+(A)$ and $\nu^-(I_A) = \nu^-(A)$ for all $A \in \calP(\N)$. Similarly, $\nu(I_A) = \nu(A)$ for all $A \in \calF$.

Some key concepts from the theory of charges will be required throughout what follows: {\em Peano-Jordan completion}, {\em outer charges}, {\em null functions}, {\em equivalence almost everywhere (a.e.)} and {\em dominance almost everywhere}. The following definitions are reproduced from~\cite{keith2022}. Similar definitions are standard in the literature on charges (see for example~\cite{bhaskararao1983} and~\cite{basile2000}). 

\begin{definition} \label{PJ_completion}
Let $(X,\calA,\mu)$ be a charge space. The {\em Peano-Jordan completion} of $(X,\calA,\mu)$ is the charge space $(X,\overline{\calA},\overline{\mu})$ where
\[
\overline{\calA} := \{ A \subseteq X : \forall \epsilon > 0, \exists B, C \in \calA \mbox{ such that } B \subseteq A \subseteq C \mbox{ and } \mu(C \setminus B) < \epsilon \}
\]
and
\[
\overline{\mu}(A) := \sup \{\mu(B): B \subseteq A, B \in \calA\} = \inf \{\mu(C): A \subseteq C, C \in \calA\}.
\]
A charge space is said to be {\em Peano-Jordan} complete if it is equal to its Peano-Jordan completion. 
\end{definition}

Peano-Jordan completion respects $\calF$ and $\nu$ in the following sense.

\begin{proposition} \label{PJlemma}
Consider any field of sets $\calA \subset \calF$. Then 
\begin{enumerate}
\item $\overline{\calA} \subset \calF$,
\item For all $A \in \overline{\calA}$, 
\begin{eqnarray*}
\nu(A) &=& \sup \{\nu(B): B \subseteq A, B \in \calA\} \\
&=& \inf \{\nu(C): A \subseteq C, C \in \calA\}, \mbox{ and }
\end{eqnarray*}
\item $(\N,\overline{\calA},\nu)$ is the Peano-Jordan completion of $(\N,\calA,\nu)$.
\end{enumerate}
\end{proposition}

\noindent {\bf Proof: } 
Given $A \in  \overline{\calA}$ and $\epsilon>0$, choose $B, C \in \calA \mbox{ such that } B \subseteq A \subseteq C \mbox{ and } \nu(C \setminus B) < \epsilon$. Then $\nu(B)=\nu^-(B) \leq \nu^-(A) \leq \nu^+(A) \leq \nu^+(C) = \nu(C)$. Since $\nu^+(A) - \nu^-(A) \leq \nu(C) - \nu(B) < \epsilon$, letting $\epsilon \rightarrow 0$ gives $\nu^+(A) = \nu^-(A)$, and thus $A \in \calF$. 

For 2, given $A \in  \overline{\calA}$ and $\epsilon>0$, choose $B, C \in \calA$ as before. Then $\nu(C)-\epsilon < \nu(B) \leq \nu(A) \leq \nu(C) < \nu(B)+\epsilon$. Taking $\epsilon$ to $0$ yields the result. 

For 3, it is enough to note Statement 2 implies $\overline{\nu} = \nu$ on $\overline{\calA}$.
\qed

\begin{definition}
A charge $\mu$ on a field of subsets $\calA$ of a set $X$ can be extended to an {\em outer charge} on $\calP(X)$ as follows
\[
\mu^*(A) := \inf\{ \mu(B): B \in \calA, A \subseteq B \}
\]
for all $A \in \calP(X)$. 
\end{definition}

\begin{definition} \label{null_function}
Let $(X,\calA,\mu)$ be a charge space. A {\em null function} is a function $f : X \rightarrow \mathbb{R}$ such that
\[
\mu^*(\{ x \in X : \lvert f(x) \rvert > \epsilon \}) = 0
\]
for all $\epsilon > 0$.
\end{definition}

\begin{definition} \label{equal_ae}
Let $(X,\calA,\mu)$ be a charge space. Two functions $f,g: X \rightarrow \mathbb{R}$ are said to be {\em equal almost everywhere} (abbreviated as $f = g$ a.e.) if $f - g$ is a null function. The function $f$ is said to be {\it dominated almost everywhere} by $g$ (abbreviated as $f \leq g$ a.e.)  if $f \leq g + h$, where $h$ is a null function. 
\end{definition}

The $K_p(\calA)$ function spaces that are the focus of this paper can now be defined.

\begin{definition} \label{Kp}
Let $\calA$ be a field of subsets of $\N$ and let $p \in [1,\infty)$. Define $K_p(\calA)$ to be the space consisting of all functions $h : \N \rightarrow \mathbb{R}$ such that the following conditions hold.
\begin{enumerate}
\item There exists a countable set $C \subset \mathbb{R}$ such that
\[
\{ h^{-1}(y,\infty) : y \in \mathbb{R} \setminus C \} \subseteq \overline{\calA}. 
\]
\item For all $\epsilon > 0$, there exists $y \in (0, \infty)$ such that 
\[ 
\nu^+( g I_{g^{-1}(y^p,\infty)} ) < \epsilon,
\]
where $g := \lvert h \rvert^p$. 
\end{enumerate}
Let $\calK_p(\calA) := K_p(\calA) / \sim$, that is, the collection of equivalence classes of $K_p(\calA)$ under the equivalence relation $h_1 \sim h_2 \iff h_1 = h_2$ a.e., for $h_1, h_2 \in K_p(\calA)$. 
\end{definition}

The two conditions above resemble the characterisation of an integrable function in Theorem~3.9 of~\cite{keith2022}. The importance of the second condition can be seen in the following example. Let $f(k) = 0$ for all $k \in \N$, and define
\[
g(k):=\begin{cases}
m & \mbox{ if $k=m^2$ is a perfect square}  \\
0 & \mbox{ if $k$ is not a perfect square}
\end{cases}.
\]
Then $f$ and $g$ differ only on the set of squares, which is a null set. Thus, $f=g$ a.e. However, clearly $\nu(f) = 0$, while $\nu(g) = \frac{1}{2}$, as can be shown in the following way. Given any $k \in \N$, let $m^2$ be the largest perfect square that is no bigger than $k$. Since $g(k) = 0$ for $m^2< k < (m+1)^2$, it follows that $\nu_{(m+1)^2-1}(g) \leq \nu_k(g) \leq \nu_{m^2}(g)$. But then
$$
\nu_{m^2}(g) = \frac{1}{m^2} \sum_{j=1}^m j = \frac{m(m+1)}{2m^2},
$$
and 
$$
\nu_{(m+1)^2-1}(g) = \frac{1}{(m+1)^2-1} \sum_{j=1}^m j = \frac{m(m+1)}{2((m+1)^2-1)},
$$
so $\lim_{m \to \infty} \nu_{m^2}(g) = \lim_{m \to \infty} \nu_{(m+1)^2-1}(g) = \frac{1}{2}$, and it follows that $\nu(g) = \frac{1}{2}$. Thus the analogy between Ces\`aro limits and integrals fails for some sequences, in the sense that two sequences that are equal almost everywhere may have different  Ces\`aro limits. But the second condition in Definition~\ref{Kp} excludes $g$ (and similar anomalies) from $K_p(\calA)$ for any field of sets $\calA \subset \calF$, since $\lim_{N \to \infty} \nu_N(g I_{g^{-1}(y,\infty)}) = \frac{1}{2}$ for all $y \in (0,\infty)$.

Many of the properties of $L_p$ spaces also hold for $K_p$ spaces. To describe some of these, further definitions used in the theory of charges are required (again see~\cite{keith2022,bhaskararao1983} for more details).

\begin{definition}
A {\em simple function} on a charge space $(X,\calA,\mu)$ is a function $f: X \rightarrow \mathbb{R}$ of the form $\sum_{k=1}^K c_k I_{A_k}$, for $K \in \N$, real numbers $c_1, \ldots, c_K$ and a partition of $X$ into subsets $\{ A_1, \ldots, A_K \} \subseteq \calA$.
\end{definition}

\begin{definition}
Let $f$ and $\{ f_n \}_{n=1}^{\infty}$ be real-valued functions on a charge space $(X,\calA,\mu)$. The sequence $\{ f_n \}_{n=1}^{\infty}$ is said to {\em converge hazily} to $f$, abbreviated as $f_n \xrightarrow{h,\calA} f$ or $f_n \xrightarrow{h} f$ if $\calA$ is clear from the context, if for every $\epsilon > 0$,
\[
\mu^*(\{ x: \lvert f_n(x) - f(x) \rvert > \epsilon \}) \rightarrow 0
\]
as $n \rightarrow \infty$. If $\{ f_n \}_{n=1}^{\infty}$ is a sequence of simple functions with $f_n \xrightarrow{h} f$, then $f$ is said to be {\em $T_1$-measurable}. The set of all $T_1$-measurable functions on $(X,\calA,\mu)$ is denoted $L_0(X,\calA,\mu)$.
\end{definition}

\begin{definition} \label{integrable}
Any simple function on a charge space $(X,\calA,\mu)$ is {\em integrable} with integral
\[
\int \sum_{k=1}^K c_k I_{A_k} := \sum_{k=1}^K c_k \mu(A_k).
\]
A general function $f: X \rightarrow \mathbb{R}$ is said to be {\em integrable} if there is a sequence of simple functions $\{ f_n \}_{n=1}^{\infty}$ such that:
\begin{enumerate}
\item $f_n \xrightarrow{h} f$, and
\item $\int \lvert f_n - f_m \rvert d\mu \rightarrow 0$ as $n,m \rightarrow \infty$.
\end{enumerate}
The integral is given by 
\[
\int f d\mu := \lim_{n \rightarrow \infty} \int f_n d\mu.
\]
For $p \in [1,\infty)$, the function space $L_p(X,\calA,\mu)$ is the set of all $T_1$-measurable functions $f: X \rightarrow \mathbb{R}$ such that $\lvert f \rvert^p$ is integrable.
\end{definition}

In what follows, $X = \N$ and $\mu$ is the restriction of $\nu$ to $\calA$, hence $L_0(\N,\calA,\nu)$ will be abbreviated as $L_0(\calA)$ and $L_p(\N,\calA,\nu)$ will be abbreviated as $L_p(\calA)$.

\begin{proposition} \label{Kp_basic_props}
Consider a field $\calA \subset \calF$. Let $r,p \in [1,\infty)$. Let $h \in K_p(\calA)$ and suppose $f$ is a real-valued function on $\N$.
\begin{enumerate}
\item $K_p(\calA)$ contains the simple functions with respect to $\overline{\calA}$.
\item $K_p(\overline{\calA}) = K_p(\calA)$.
\item if $\calA^{\prime} \subseteq \calA$ is a field, then $K_p(\calA^{\prime}) \subseteq K_p(\calA)$.
\item $h$ is $T_1$-measurable.
\item $K_p(\calA) \subseteq L_p(\calA)$.
\item If $r \leq p$ then $K_r(\calA) \supseteq K_p(\calA)$.
\item If $\lvert f \rvert \leq h$, then $f \in K_p(\calA)$ if and only if it is $T_1$-measurable. 
\item $f \in K_p(\calA)$ if and only if $f I_A \in K_p(\calA)$ for all $A \in \overline{\calA}$.
\item $\lvert h \rvert \in K_p(\calA)$ and $\lvert h \rvert^p \in K_1(\calA)$.
\item $f \in K_p(\calA) \iff (f^+)^p, (f^-)^p \in K_1(\calA) \iff f^+, f^- \in K_p(\calA)$.
\item $K_p(\calA)$ is a real vector space.
\end{enumerate}
\end{proposition}

The standard notation $f^+ := \max \{ f, 0 \}$ and $f^- := \max\{ -f, 0 \}$ is used in Property~10, and throughout the rest of the paper.

\noindent {\bf Proof: } 
Properties~1 to~3 are immediate from Definition~\ref{Kp}.

For 4, note $h$ immediately satisfies Statement~2a of Theorem~3.4 in~\cite{keith2022}, since it is the same as the first condition in Definition~\ref{Kp}. Moreover, $h$ is smooth, since
\[
\nu^*(\lvert h \rvert^{-1}(y,\infty)) = \nu(\lvert h \rvert^{-1}(y,\infty)) = \nu^+(g^{-1}(y^p,\infty)) \leq \nu^+(g I_{g^{-1}(y^p,\infty))})
\]
for any $y \in [1,\infty) \setminus C$, where $g := \lvert h \rvert^p$ and $C$ is the countable set referenced in Definition~\ref{Kp}. (The first equality follows by Proposition~\ref{PJlemma}.) Thus $h$ is $T_1$-measurable, by Theorem~3.4 in~\cite{keith2022}. 

For 5, use the fact $h$ is $T_1$-measurable by Property~4. Let $\{ h^+_n \}_{n=1}^{\infty}$ and $\{ h^-_n \}_{n=1}^{\infty}$ be the sequences obtained by applying Statement~3 of Theorem~3.4 in~\cite{keith2022} to $h$ and let $Y := \{ y_{n,j} : n \in \N, j \in \{1, \ldots, K_n \} \}$ where $y_{n,j}$ is as described in that statement. Let $y_n := y_{n,K_n}$. 

First note $h^+_n \xrightarrow{h,\calA} h^+$, hence $(h^+_n)^p \xrightarrow{h,\calA} (h^+)^p$. Moreover, given $\epsilon > 0$, one may choose $m$ sufficiently large that 
\begin{enumerate}
\item $\lvert (h_m^+)^p(x) - (h^+)^p(x) \rvert < \epsilon/2$ for $x \in (h^+)^{-1}[0,y_m^p]$, and
\item $\nu^+(g I_{g^{-1}(y_m^p,\infty)}) < \epsilon/2$, 
\end{enumerate}
where $g := \lvert h \rvert^p$. (To see 1, recall $y_{m,1} \leq 2^{-m}$ and $y_{m,j+1} - y_{m,j}  \leq 2^{1-m}$ for $j \in \{ 1, \ldots, K_m - 1 \}$. Hence $y_{m,1}^p \leq 2^{-mp}$ for any $p \in [1,\infty)$. 
Also, $y^p$ is convex, hence by the mean value theorem $y_{m,j+1}^p - y_{m,j}^p  \leq 2^{1-m} p m^{p-1}$, since $y_m \leq m$.)

Hence for $n > m$,
\begin{eqnarray*}
\int \lvert (h^+_m)^p - (h^+_n)^p \rvert d\nu_{\overline{\calA}} & \leq & \int \lvert (h^+_m)^p - (h^+)^p \rvert I_{g^{-1}[0,y_m^p]} d\nu_{\overline{\calA}} \\
& &+ \int (h_n^+)^p I_{g^{-1}(y_m^p,\infty)} d\nu_{\overline{\calA}} \\
&\leq& \frac{\epsilon}{2} \nu(\lvert h \rvert^{-1}([0,y_m])) + \nu((h^+_n)^p I_{g^{-1}(y_m^p,\infty)}) \\
& \leq & \frac{\epsilon}{2} + \nu^+(g I_{g^{-1}(y_m^p,\infty)}) \\
& < & \epsilon,
\end{eqnarray*}
where the second inequality uses the fact that $(h^+_n)^p$ is a simple function with respect to $\overline{\calA}$, and the definition of the integral of a simple function. 
Hence $\{ (h^+_n)^p \}_{n=1}^{\infty}$ is a determining sequence for $(h^+)^p$ with respect to $\overline{\calA}$, and $h^+ \in L_p(\overline{\calA})$. By similar reasoning, $h^- \in L_p(\overline{\calA})$, and thus $h \in L_p(\overline{\calA})$. Finally, $h \in L_p(\calA)$ by Proposition~1.8 of~\cite{basile2000}.

For 6, note the first condition in Definition~\ref{Kp} is the same for both $K_r(\calA)$ and $K_p(\calA)$. Also,
\[
g_r I_{g_r^{-1}(y^r,\infty)} \leq g_p I_{g_p^{-1}(y^p,\infty)}
\]
for $y \in (1,\infty)$, where $g_p := \lvert h \rvert^p$ and $g_r := \lvert h \rvert^r$. Hence
\[
\nu^+(g_r I_{g_r^{-1}(y^r,\infty)}) \leq \nu^+(g_p I_{g_p^{-1}(y^p,\infty)}),
\]
and the second condition in the definition of $K_r(\calA)$ holds for $h$. 

For 7, the forward implication is immediate from Property~4. For the reverse implication, first consider the case $p = 1$, and observe that since both $f$ and $h$ are $T_1$-measurable ($h$ by Property~4), there is a countable set $C \subset \mathbb{R}$ such that for $y \in \mathbb{R} \setminus C$, both $f^{-1}(y,\infty)$ and $h^{-1}(y,\infty)$ are in $\overline{\calA}$, by Theorem~3.4 in~\cite{keith2022}. Fix $\epsilon > 0$. Then there exists $y \in (0,\infty)$ such that $\nu^+(h I_{h^{-1}(y,\infty)}) < \epsilon$. Note
\[
\lvert f \rvert I_{\lvert f \rvert^{-1}(y,\infty)} \leq h I_{h^{-1}(y,\infty)}
\]
and hence $\nu^+(\lvert f \rvert I_{\lvert f \rvert^{-1}(y,\infty)}) < \epsilon$. Thus $f \in K_1(\calA)$. The case $p > 1$ will be proved after Properties~9 and~10.

Property~8 can be shown by a similar argument to Corollary~2.28 of~\cite{keith2022}, using Property~7.

For 9, recall $h$ is $T_1$-measurable by Property~4, hence $\lvert h \rvert$ is $T_1$-measurable, and by Theorem~3.4 of~\cite{keith2022} there exists a countable set $C \subset (0,\infty)$ such that for $y \in (0, \infty) \setminus C$, $\lvert h \rvert^{-1}(y,\infty) \in \overline{\calA}$. The second condition of Definition~\ref{Kp} is satisfied for $\lvert h \rvert$ since it is satisfied for $h$, hence $\lvert h \rvert \in K_p(\calA)$. But then $\lvert h \rvert^p$ satisfies the conditions for membership of $K_1(\calA)$, with the required countable set being $C^p := \{ y^p : y \in C \}$.

For 10, the first claim implies the second because if $f \in K_p(\calA)$ then $f$ is $T_1$-measurable by Property~4, hence $(f^+)^p$ and $(f^-)^p$ are $T_1$-measurable. Moreover, $(f^+)^p$ and $(f^-)^p$ are dominated by $\lvert f \rvert^p$, hence $(f^+)^p, (f^-)^p \in K_1(\calA)$ by Property~7. The second claim implies the third by the definition of $K_p(\calA)$. The third claim implies the first because $h$ is $T_1$-measurable, since $h := h^+ - h^-$, and hence the first condition of Definition~\ref{Kp} holds by Theorem~3.4 of~\cite{keith2022}. The second condition of that theorem holds because
\[
g I _{g^{-1}(y^p,\infty)} = g_+ I _{g_+^{-1}(y^p,\infty)} + g_- I _{g_-^{-1}(y^p,\infty)}
\]
where $g := \lvert h \rvert^p$, $g_+ := (h^+)^p$ and $g_- := (h^-)^p$.

Returning to 7, consider the case $p > 1$ and suppose $f$ is $T_1$-measurable. Then $(f^+)^p$ and $(f^-)^p$ are $T_1$-measurable, $\lvert f \rvert^p \leq h^p$, and $h^p \in K_1(\calA)$ by Property~9. Thus Property~7 with $p=1$ (proved above) gives $(f^+)^p, (f^-)^p \in K_1(\calA)$. But then Property~10 gives $f \in K_p(\calA)$.

For 11, consider $h_1, h_2 \in K_p(\calA)$, and let $C_1, C_2$ be the respective countable sets asserted by Definition~\ref{Kp}. Then $h_1, h_2$ are $T_1$-measurable by Property~4, hence $h_3 := h_1 + h_2$ is $T_1$-measurable. By Theorem~3.4 of~\cite{keith2022}, there exists countable $C_3 \subset \mathbb{R}$ such that $h_3^{-1}(y,\infty) \in \overline{\calA}$ for $y \in \mathbb{R} \setminus C_3$.

Fix $\epsilon > 0$, then there exists $y \in (0,\infty)$ such that 
\begin{align*}
\nu^+( g_1 I_{g_1^{-1}(y^p/2,\infty)} ) &< \frac{\epsilon}{3 \cdot 2^p}, \mbox{ and } \\
\nu^+( g_2 I_{g_2^{-1}(y^p/2,\infty)} ) &< \frac{\epsilon}{3 \cdot 2^p} 
\end{align*}
where $g_1 := \lvert h_1 \rvert^p$ and $g_2 := \lvert h_2 \rvert^p$. Set $g_3 := \lvert h_3 \rvert^p$, then
\[
g_3^{-1}(y^p,\infty) \subseteq g_1^{-1}(y^p/2,\infty) \cup g_2^{-1}(y^p/2,\infty).
\]
For $p \in [1,\infty)$, the map $x \mapsto \lvert x \rvert^p$ is convex, and hence
\[
\left\lvert \frac{ h_1 + h_2 }{2} \right\rvert^p \leq \frac{ \lvert h_1 \rvert^p + \lvert h_2 \rvert^p }{2}.
\]
That is, $g_3 \leq 2^{p-1} (g_1 + g_2)$. Putting this together gives
\begin{eqnarray*}
g_3 I_{g_3^{-1}(y^p,\infty)} &\leq& 2^{p-1} (g_1 + g_2)(I_{g_1^{-1}(y^p/2,\infty)} + I_{g_2^{-1}(y^p/2,\infty)}) \\
&=& 2^{p-1} (g_1 I_{g_1^{-1}(y^p/2,\infty)} + g_1 I_{g_2^{-1}(y^p/2,\infty)} + \\
& & g_2 I_{g_1^{-1}(y^p/2,\infty)} + g_2 I_{g_2^{-1}(y^p/2,\infty)}) \\
&\leq& 3 \cdot 2^{p-1} (g_1 I_{g_1^{-1}(y^p/2,\infty)} + g_2 I_{g_2^{-1}(y^p/2,\infty)}) 
\end{eqnarray*}
where the last line follows because
\[
g_1 I_{g_2^{-1}(y^p/2,\infty)} \leq g_1 I_{g_1^{-1}(y^p/2,\infty)} + g_2 I_{g_2^{-1}(y^p/2,\infty)}
\]
since either $g_1(x) \leq g_2(x)$ or $I_{g_2^{-1}(y^p/2,\infty)}(x) \leq I_{g_1^{-1}(y^p/2,\infty)}(x)$, and similarly
\[
g_2 I_{g_1^{-1}(y^p/2,\infty)} \leq g_1 I_{g_1^{-1}(y^p/2,\infty)} + g_2 I_{g_2^{-1}(y^p/2,\infty)}
\]
since either $g_2(x) \leq g_1(x)$ or $I_{g_1^{-1}(y^p/2,\infty)}(x) \leq I_{g_2^{-1}(y^p/2,\infty)}(x)$. But then 
\[\nu^+(g_3 I_{g_3^{-1}(y^p,\infty)}) < \epsilon.
\]
Thus $K_p(\calA)$ is closed under addition. The other axioms of a real vector space are straightforward to prove.
\qed  

Proposition~\ref{Kp_basic_props}(7) requires strict dominance of $\lvert f \rvert$ by $g$: the condition $\lvert f \rvert \leq g$ a.e. is not sufficient. To see this, let $h$ be a null function that is not in $K_1(\calA)$ (an example of such a null function is given following Definition~\ref{Kp}). Let $f := g + h$. Then $f \leq g$ a.e., but one cannot have both $f$ and $g$ in $K_1(\calA)$, by Proposition~\ref{Kp_basic_props}(11).

The next theorem extends the concept of null modification to integrable functions, adding to the null modifications obtained for sets (Proposition~\ref{nullmodification}), chains (Theorem~\ref{nullmodification_thm}), and sequences (Corollary~\ref{nullmodification_sequence}) in Section~\ref{null_modification_section}. The value of this null modification for functions is that it maps any $f \in L_p(\calA)$ to a sequence $h$ with a Ces\`aro limit, equal to $f$ almost everywhere. Thus every equivalence class in $\calL_p(\calA)$ contains at least one sequence with a Ces\`aro limit. Moreover, that Ces\`aro limit must be equal to the integral of $f$ (since functions that are equal almost everywhere have the same integral).

Here and throughout the rest of the paper, the notation $\alpha(\calC)$ represents the field of sets generated by a collection of sets $\calC$, that is, the smallest field of sets containing $\calC$. 

\begin{theorem} \label{nullmodification_fn}
Consider a field of sets $\calA \subset \calF$ and $f \in L_p(\calA)$ for some $p \in [1,\infty)$. There exists $h \in K_p(\alpha(\calA \cup \calN))$ such that:
\begin{enumerate}
\item $h = f$ a.e. with respect to $\alpha(\calA \cup \calN)$,
\item $\nu_N(\lvert h \rvert^p)$ converges as $N \rightarrow \infty$, with $\int \lvert h \rvert^p d\nu = \nu (\lvert h \rvert^p)$, and
\item $\nu_N(h)$ converges as $N \rightarrow \infty$, with $\int h d\nu = \nu(h)$.
\end{enumerate}
\end{theorem}

The proof requires several lemmas. To begin with, the following lemma identifies special conditions under which an element of a $K_p$ space can be shown to have a Ces\`aro limit equal to its integral.

\begin{lemma} \label{Kp_Cesaro_integral}
Consider a field $\calA \subset \calF$, $p \in [1,\infty)$ and a function $h \in L_p(\calA)$ such that
\begin{enumerate}
\item $\calR := \{ h^{-1}(y,\infty) : y \in \mathbb{R} \} \subset \calF$, 
\item for any $\calC \subseteq \calR$, 
\begin{enumerate}
\item $\bigcup \calC \in \calF$ with $\nu(\bigcup \calC) = \sup \{ \nu(A) : A \in \calC \}$,
\item $\bigcap \calC \in \calF$ with $\nu(\bigcap \calC) = \inf \{ \nu(A) : A \in \calC \}$, and
\end{enumerate}
\item $\nu_N(A) \leq \nu(A)$ for all $A \in \calR$ and $N \in \N$.
\end{enumerate}
Let $\{ h^+_n \}_{n=1}^{\infty}$ and $\{ h^-_n \}_{n=1}^{\infty}$ be the sequences of functions obtained by applying Statement~3 of Theorem~3.4 in~\cite{keith2022} to $h$. Then the following statements hold:
\begin{enumerate}
\item $\nu_N((h^+_n)^p)$ and $\nu_N((h^-_n)^p)$ converge uniformly over $n \in \N$ as $N \rightarrow \infty$, 
\item $\nu_N((h^+)^p)$ converges as $N \rightarrow \infty$, with $\int (h^+)^p d\nu = \nu((h^+)^p)$,
\item $\nu_N((h^-)^p)$ converges as $N \rightarrow \infty$, with $\int (h^-)^p d\nu = \nu((h^-)^p)$,
\item $\nu_N(\lvert h \rvert^p)$ converges as $N \rightarrow \infty$, with $\int \lvert h \rvert^p d\nu = \nu (\lvert h \rvert^p)$,
\item $\nu_N(h)$ converges as $N \rightarrow \infty$, with $\int h d\nu = \nu(h)$, and
\item $h \in K_p(\calA)$.
\end{enumerate}
\end{lemma}

\noindent {\bf Proof: }
Let $Y := \{ y_{n,j} : n \in \N, j \in \{1, \ldots, K_n \} \}$ where $y_{n,j}$ is as described in Theorem~3.4 of~\cite{keith2022}. Let $y_n := y_{n,K_n}$. Then
\[
h^+_n := \sum_{j=1}^{K_n - 1} y_{n,j} I_{C^+_{n,j} \setminus C^+_{n,j+1}} + y_n I_{C^+_n}
\mbox{ and }
h^-_n := \sum_{j=1}^{K_n - 1} y_{n,j} I_{C^-_{n,j} \setminus C^-_{n,j+1}} + y_n I_{C^-_n}
\]
where $C^+_{n,j} := (h^+)^{-1}(y_{n,j},\infty)$, $C^+_n := (h^+)^{-1}(y_n,\infty)$, $C^-_{n,j} := (h^-)^{-1}(y_{n,j},\infty)$ and $C^-_n := (h^-)^{-1}(y_n,\infty)$ for $n \in \N$ and $j \in \{ 1, \ldots, K_n \}$. 

Fix $\epsilon > 0$. By Property~2b of Theorem~3.9 in~\cite{keith2022}, there exists $k \in \N$ such that 
\[ 
\int (h^+)^p I_{C^+_k} d\nu < \frac{\epsilon}{4}.
\]
Moreover, by Theorem~\ref{uniformconvergence} and Condition~2 above, there exists $N_{\epsilon} \in \N$ such that 
\[
\lvert (\nu_N - \nu) \big( C^+_{n,j} \big) \rvert = \lvert \nu_N(C^+_{n,j}) - \nu(C^+_{n,j}) \rvert < \frac{\epsilon}{2y_k^p}
\]
for all $n \in \N$, $j \in \{ 1, \ldots, K_n \}$ and $N > N_{\epsilon}$. Choose $N > N_{\epsilon}$ and $n \in \N$. Let $J$ be the largest integer in the set $\{ 1, \ldots, K_n \}$ such that $y_{n,J} \leq y_k$. Note $J = K_n$ for $n \leq k$ and $y_{n,J} = y_k$ for $n \geq k$. Then
\begin{eqnarray*}
\left\lvert (\nu_N - \nu) \big( (h^+_n)^p I_{(C^+_k)^c} \big) \right\rvert & = & \bigg\lvert \sum_{j=1}^{J - 1} y_{n,j}^p (\nu_N-\nu) \big(C^+_{n,j} \setminus C^+_{n,j+1} \big) \\
& & + y_n^p (\nu_N - \nu) \big( C^+_n \big)  \bigg\rvert \\
& = & \bigg\lvert y_{n,1}^p (\nu_N - \nu) \big( C^+_{n,1} \big)  \\
& & + \sum_{j=2}^J (y_{n,j}^p - y_{n,j-1}^p) (\nu_N - \nu) \big( C^+_{n,j} \big) \bigg\rvert \\
& < & \frac{\epsilon}{2y_k^p} \bigg( y_{n,1}^p +  \sum_{j=2}^J \big( y_{n,j}^p - y_{n,j-1}^p \big) \bigg) \\
& \leq & \frac{\epsilon}{2}.
\end{eqnarray*}
For $n \leq k$, $(h^+_n)^p I_{C^+_k} = y_n^p I_{C^+_k}$, hence 
\[
\big\lvert (\nu_N - \nu) \big( (h^+_n)^p I_{C^+_k} \big) \big\rvert = y_n^p \lvert (\nu_N - \nu)(C^+_k) \rvert < y_n^p \frac{\epsilon}{2 y_k^p} \leq \frac{\epsilon}{2}.
\] 
For $n > k$,
\begin{eqnarray*}
\left\lvert (\nu_N - \nu) \big( (h^+_n)^p I_{C^+_k} \big) \right\rvert & = & \left\lvert \int (h^+_n)^p I_{C^+_k} d\nu_N - \int (h^+_n)^p I_{C^+_k} d\nu  \right\rvert \\
& \leq & \int (h^+)^p I_{C^+_k} d\nu_N + \int (h^+)^p I_{C^+_k} d\nu \\
& \leq & 2 \int (h^+)^p I_{C^+_k} d\nu \\
& < & \frac{\epsilon}{2}
\end{eqnarray*}
where Line~1 uses the fact that $(h^+_n)^p I_{C^+_k}$ is a simple function and Line~3 uses Corollary~3.10 of~\cite{keith2022} with $E = Y$ (recalling that $\nu_N$ is a measure, and hence a charge). Putting these inequalities together gives $\lvert \nu_N((h^+_n)^p) - \nu((h^+_n)^p) \rvert < \epsilon$, as required. Uniform convergence of $\nu_N((h^-_n)^p)$ follows by similar reasoning.

To show Statement~2, fix $\epsilon > 0$. By Statement~1, there exists $N_{\epsilon}$ such that 
\[
\lvert \nu_N((h^+_n)^p) - \nu((h^+_n)^p) \rvert < \frac{\epsilon}{3}
\]
for $N > N_{\epsilon}$ and all $n \in \N$. For any $N > N_{\epsilon}$, one may choose $k$ sufficiently large that the following conditions are met.
\begin{enumerate}
\item $\left\lvert \nu_N((h^+)^p) - \nu_N((h^+_k)^p) \right\rvert < \epsilon / 3$, and
\item $\left\lvert \int (h^+_k)^p d\nu - \int (h^+)^p d\nu \right\rvert < \epsilon/3$,
\end{enumerate}
where the former condition holds because $(h^+_k)^p$ can approximate $(h^+)^p$ arbitrarily closely on the finite set $\{ 1, \ldots, N \}$ and the latter condition holds by dominated convergence (Theorem~4.6.14 of~\cite{bhaskararao1983}) with $(h^+)^p$ as the dominating function. Then
\begin{eqnarray*}
\bigg\lvert \nu_N((h^+)^p) - \int (h^+)^p d\nu \bigg\rvert & \leq & \lvert \nu_N((h^+)^p) - \nu_N((h^+_k)^p) \rvert \\
& & + \lvert \nu_N((h^+_k)^p) - \nu((h^+_k)^p) \rvert \\
& &  + \left\lvert \nu((h^+_k)^p) - \int (h^+_k)^p d\nu \right\rvert \\
& & + \left\lvert \int (h^+_k)^p d\nu - \int (h^+)^p d\nu \right\rvert \\
&<& \epsilon
\end{eqnarray*}
The third summand is zero by the definition of the integral of a simple function. The other three terms are each smaller than $\epsilon/3$ by construction. Hence $\nu_N((h^+)^p)$ converges to a finite limit $\nu((h^+)^p) = \int (h^+)^p d\nu$ as $N \rightarrow \infty$. 

Statement~3 follows by similar reasoning. Statement~4 follows by the additivity of integrals and limits, noting that $\lvert h \rvert^p = (h^+)^p + (h^-)^p$. Statement~5 follows since $h \in L_1(\calA)$ by Corollary~4.6.5 of~\cite{bhaskararao1983} and $h = h^+ - h^-$.

For 6, note the first condition of Definition~\ref{Kp} holds for $h$ by Theorem~3.9 of~\cite{keith2022}. To show the second condition, fix $\epsilon > 0$, then by Theorem~3.9 of~\cite{keith2022} there exists $z \in (0,\infty)$ such that $\int \lvert h \rvert^p I_{(\lvert h \rvert^p)^{-1}(z^p, \infty)} d\nu < \epsilon$. By Theorem~3.4 of~\cite{keith2022}, there exists $y > z$ such that $h^{-1}(y^p, \infty) \in \overline{\calA}$ and $h^{-1}(-y^p, \infty) \in \overline{\calA}$. Hence
\begin{align*}
h_1 &:= h I_{h^{-1}(y^p, \infty)} \in L_p(\calA) \mbox{ and } \\
h_2 &:= h I_{h^{-1}(-y^p, \infty)} \in L_p(\calA),
\end{align*}
by Corollary~2.28 of~\cite{keith2022}. Define
\begin{align*}
\calR_1 &:= \{ h_1^{-1}(y,\infty) : y \in \mathbb{R} \}, \mbox{ and } \\
\calR_2 &:= \{ h_2^{-1}(y,\infty) : y \in \mathbb{R} \}
\end{align*}
and note $\calR_1 \subseteq \calR$ and $\calR_2 \subseteq \calR$. Thus Statement~4 of this lemma may be applied to $h_1$ and $h_2$, giving that $\nu_N(\lvert h_1 \rvert^p)$ and $\nu_N(\lvert h_2 \rvert^p)$ converge to $\int \lvert h_1 \rvert^p d\nu$ and $\int \lvert h_2 \rvert^p d\nu$, respectively. It follows that:
\begin{eqnarray*}
\nu^+ \big( \lvert h \rvert^p I_{(\lvert h \rvert^p)^{-1}(z^p,\infty)} \big) &\leq& \nu(\lvert h \rvert^p) - \big( \nu(\lvert h_2 \rvert^p) - \nu(\lvert h_1 \rvert^p) \big) \\
&=& \int \lvert h \rvert^p d\nu - \bigg( \int \lvert h_2 \rvert^p d\nu - \int \lvert h_1 \rvert^p d\nu \bigg) \\
&=& \int \lvert h \rvert^p I_{(\lvert h \rvert^p)^{-1}(y^p, \infty)} d\nu \\
&<& \epsilon,
\end{eqnarray*}
as required.
\qed 

The question remains whether non-trivial functions satisfying the conditions of Lemma~\ref{Kp_Cesaro_integral} exist. The next lemma uses null modification for chains (Theorem~\ref{nullmodification_thm}) to construct functions with the desired properties. 

\begin{lemma} \label{nullmodification_posfn}
Consider a field of sets $\calA \subset \calF$ and a non-negative function $f \in L_0(\calA)$. Then there exists a non-negative function $h \in L_0(\alpha(\calA \cup \calN))$ such that:
\begin{enumerate}
\item $f = h$ a.e. (with respect to $\alpha(\calA \cup \calN)$),
\item $\calR := \{ h^{-1}(y,\infty) : y \in [0,\infty) \} \subset \calF$, 
\item for any $\calC \subseteq \calR$, 
\begin{enumerate}
\item $\bigcup \calC \in \calF$ with $\nu(\bigcup \calC) = \sup \{ \nu(A) : A \in \calC \}$,
\item $\bigcap \calC \in \calF$ with $\nu(\bigcap \calC) = \inf \{ \nu(A) : A \in \calC \}$, and
\end{enumerate}
\item $\nu_N(A) \leq \nu(A)$ for all $A \in \calR$ and $N \in \N$.
\end{enumerate}
Moreover, if $f \in L_p(\calA)$ for $p \in [1,\infty)$, then $h \in K_p(\alpha(\calA \cup \calN))$.
\end{lemma}

\noindent {\bf Proof: }
Without loss of generality, one may assume $\overline{\calA} = \calA$ and $\calN \subset \calA$, since if the lemma is proved with $\calA$ replaced by $\overline{\alpha(\calA \cup \calN)}$, then it is true also for $\calA$. To see this, note the condition $f \in L_p(\calA)$ implies $f \in L_0(\overline{\alpha(\calA \cup \calN)})$ for any $p \in \{0\} \cup [1,\infty)$, and the conclusions $h \in L_0(\overline{\alpha(\calA \cup \calN)})$ and $h \in K_p(\overline{\alpha(\calA \cup \calN)})$ respectively imply $h \in L_0(\alpha(\calA \cup \calN))$ (by Proposition~1.8 of~\cite{basile2000}) and $h \in K_p(\alpha(\calA \cup \calN))$ (by Property~2 of Proposition~\ref{Kp_basic_props}).

As in the proof of Lemma~5.9 of~\cite{keith2022}, define
\begin{align*}
A_y &:= f^{-1}(y, \infty) \mbox{ for each } y \in (0,\infty) \setminus C, \mbox{ and }\\
\calT &:= \{ A_y : y \in (0,\infty) \setminus C \},
\end{align*}
where $C$ is the countable set referenced in Statement~2a of Theorem~3.4 in~\cite{keith2022}, hence $\calT \subseteq \calA$. The overall strategy in what follows is to use null modification to modify $\calT$ to create a new chain $\calR$ for which Properties~2 to~4 hold, and then construct a new function $h$ from $\calR$. 

Let $\phi : \calT \rightarrow \calF$ be the order preserving map asserted in Theorem~\ref{nullmodification_thm}, as it applies to $\calT$. Define
\begin{align*}
B_y &:= \phi(A_y) \mbox{ for each } y \in (0,\infty) \setminus C, \\
\calS &:= \{ B_y : y \in (0,\infty) \setminus C \}, \mbox{ and } \\
\calS_y &:= \{ B \in \calS : B \subset B_y \} \mbox{ for each } y \in (0,\infty).
\end{align*}
Then, for each $y \in (0,\infty) \setminus C$, $B_y \triangle A_y \in \calN$ and $\nu_N(B_y) \leq \nu(B_y)$ for all $N \in \N$. Moreover, $\calS \subseteq \calA$, since $\calN \subset \calA$, and Properties~3a and~3b hold for any $\calC \subseteq \calS$. 

Next define:
\begin{align*}
C_y &:= \bigcup \calS_y \setminus \bigcap \calS \mbox{ for each } y \in (0,\infty),  \\
\calR &:= \{ C_y : y \in (0,\infty) \}, \mbox{ and } \\
\calR_y &:= \{ C \in \calR : C \subset C_y \} \mbox{ for each } y \in (0,\infty).
\end{align*}
Note the following properties:
\begin{enumerate}
\item $C_y \in \calA$ for each $y \in (0,\infty) \setminus C$,
\item $C_y = \bigcup \calR_y$ for each $y \in (0,\infty)$,
\item $\bigcap \calR = \emptyset$, 
\item $C_y \triangle A_y \in \calN$ for each $y \in (0,\infty) \setminus C$,
\item Properties~3a and 3b of the lemma hold for any $\calC \subseteq \calR$, and
\item $\nu_N(C_y) \leq \nu(C_y)$ for each $y \in (0,\infty)$ and $N \in \N$.
\end{enumerate}

Claims~1 to~4 (of the preceding list) follow as in the proof of Lemma~5.9 in~\cite{keith2022}, and also $\bigcap \calS \in \calA$ with $\nu(\bigcap \calS) = 0$, as in that proof.

For 5, note that for any $Y \subseteq (0,\infty) \setminus C$,
\[
\bigcup \{ C_y : y \in Y \} = \bigcup \left\{ \bigcup \calS_y : y \in Y \right\} \setminus \bigcap \calS, 
\]
which is in $\calF$, since Properties~3a and~3b (of this lemma) hold for any $\calC \subseteq \calS$. Moreover, 
\begin{eqnarray*}
\nu \left( \bigcup \{ C_y : y \in Y \} \right) & = & \nu \left( \bigcup \{ \bigcup \calS_y : y \in Y \} \right) - \nu \left( \bigcap \calS \right)\\
& = & \sup \{ \sup \left\{ \nu(B) : B \in \calS_y \} : y \in Y \right\} \\
& = & \sup \left\{ \nu(C_y) : y \in Y \right\}.
\end{eqnarray*}
Similarly $\bigcap \{ C_y : y \in Y \} \in \calF$ with $\nu( \bigcap \{ C_y : y \in Y \}) = \inf \{ \nu(C_y) : y \in Y \}$.

For 6, observe that $C_y \subseteq \bigcup \calS_y \subseteq B_y$ and $\nu(C_y) = \nu(\bigcup \calS_y) = \nu(B_y)$, hence
\[
\nu_N(C_y) \leq  \nu_N(B_y) \leq  \nu(B_y) = \nu(C_y)
\]
for each $y \in (0,\infty) \setminus C$ and $N \in \N$. Moreover,
\begin{eqnarray*}
\nu_N(C_y) &=& \sup \{ \nu_N(C_z) : z \in (y,\infty) \setminus C \} \\ 
&\leq& \sup \{ \nu(C_z) : z \in (y,\infty) \setminus C \} \\
&=& \nu(C_y)
\end{eqnarray*}
for each $y \in C$ and $N \in \N$. (The first equality holds because for any $N$ there is $z > y$ such that $I_{C_z}(n) = I_{C_y}(n)$ for $n \in \{ 1, \ldots, N \}$. The last equality holds by~5.)

By Lemma~5.8 of~\cite{keith2022}, there is a function $h : \N \rightarrow [0,\infty)$ such that $C_y = h^{-1}(y,\infty)$ for each $y \in (0,\infty)$. As in the proof of Lemma~5.9 in~\cite{keith2022}, $h = f$ a.e. (with respect to $\alpha(\calA \cup \calN)$), and $h \in L_0(\calA)$.

Finally, if $f \in L_p(\calA)$, then $f \in L_p(\alpha(\calA \cup \calN))$, hence $h \in L_p(\alpha(\calA \cup \calN))$, and then $h \in K_p(\alpha(\calA \cup \calN))$ by Lemma~\ref{Kp_Cesaro_integral}.
\qed 

Lemmas~\ref{Kp_Cesaro_integral} and~\ref{nullmodification_posfn} can be combined to give a proof of Theorem~\ref{nullmodification_fn}, as follows.

\noindent {\bf Proof of Theorem~\ref{nullmodification_fn}: }
Applying Lemma~\ref{nullmodification_posfn} to $f^+$ and $f^-$ yields non-negative functions $h^+, h^- \in K_p(\alpha(\calA \cup \calN))$ that satisfy the conditions of Lemma~\ref{Kp_Cesaro_integral}. Then $h := h^+ - h^- \in K_p(\alpha(\calA \cup \calN))$ has the required properties. \qed

Armed with this null modification, one of the main theorems of this paper can now be proved. Informally, the theorem establishes that $K_p$ spaces are comprised of sequences for which Ces\`aro limits exist, that these limits determine a natural pseudonorm on $K_p(\calA)$ and a norm on $\calK_p(\calA)$, and that $K_p(\calA)$ can be isometrically embedded in $L_p(\calA)$. 

To state this theorem, the following notation is required. Let $[h]_{K_p(\calA)}$ denote the equivalence class of $h$ in $K_p(\calA)$ (under equality almost everywhere), to distinguish it from the equivalence class $[h]_{L_p(\calA)}$ of $h$ in $L_p(\calA)$. In a similar manner, the pseudonorm of $h$ in $L_p(\calA)$ will sometimes be denoted $\| h \|_{L_p(\calA)}$ and the norm of $[h]$ in $\calL_p(\calA)$ will sometimes be denoted $\| [h] \|_{L_p(\calA)}$, where needed to distinguish them from the pseudonorm and norm defined in the following theorem. 

\begin{theorem} \label{Kp_pseudonorm}
Let $\calA \subset \calF$ be a field of sets and $p \in [1,\infty)$. Then the following statements hold. 
\begin{enumerate}
\item $K_p(\calA)$ is a vector subspace of $L_p(\calA)$. 
\item For all $h \in K_p(\calA)$, $\nu_N(\lvert h \rvert^p)$ converges to $\int \lvert h \rvert^p d\nu$ and $\nu_N(h)$ converges to $\int h d\nu$ as $N \rightarrow \infty$.
\item The function $\| \cdot \|_{K_p(\calA)} : K_p(\calA) \rightarrow \mathbb{R}$ given by:
\[
\| h \|_{K_p(\calA)} := \left( \nu(\lvert h \rvert^p) \right)^{1/p}
\]
is a pseudonorm on $K_p(\calA)$.
\item For all $h \in K_p(\calA)$, 
\[
\| h \|_{K_p(\calA)} = \left( \int \lvert h \rvert^p d\nu \right)^{1/p} = \| h \|_{L_p(\calA)}.
\]
That is, the embedding of $K_p(\calA)$ in $L_p(\calA)$ is an isometry.
\item For $h_1, h_2 \in K_p(\calA)$, $\| h_1 - h_2 \|_{K_p(\calA)} = 0$ if and only if $h_1$ and $h_2$ are equal a.e. 
\item The quotient space $\calK_p(\calA)$ is a vector space with norm
\[
\| [h] \|_{K_p(\calA)} := \| h \|_{K_p(\calA)}.
\]
\end{enumerate}
\end{theorem}


The proof requires the following lemma, which establishes that Ces\`aro limits exist for all elements of $K_p(\calA)$, not just those that satisfy the conditions of Lemma~\ref{Kp_Cesaro_integral}.

\begin{lemma} \label{K1_equal_ae}
Consider a field of sets $\calA \subset \calF$. Then
\begin{enumerate}
\item $\nu^+(A) \leq \nu_{\calA}^*(A)$ for all $A \in \calP(\N)$,
\item for any null function $h \in K_1(\calA)$, $\nu_N(h)$ and $\nu_N(\lvert h \rvert)$ converge as $N \rightarrow \infty$, with $\nu(h) = \nu(\lvert h \rvert) = 0$, 
\item for any $f, g \in K_1(\calA)$ with $f = g$ a.e.
\begin{enumerate}
\item $\nu^+(f) = \nu^+(g)$,
\item $\nu^-(f) = \nu^-(g)$, and
\item $\nu_N(f)$ converges if and only if $\nu_N(g)$ converges, in which case $\nu(f) = \nu(g)$.
\end{enumerate}
\end{enumerate}
\end{lemma}

\noindent {\bf Proof: }
For 1, note for any $B \in \calA$ with $A \subseteq B$, $\nu^+(A) \leq \nu(B)$. Hence
\[
\nu^+(A) \leq \inf \{ \nu(B) : B \in \calA, A \subseteq B \} = \nu_{\calA}^*(A). 
\]

For 2, it will be sufficient to prove the statement for $h \geq 0$, since then it applies to $h^+$ and $h^-$, and thus to $h^+ - h^-$. Fix $\epsilon > 0$. Since $h \in K_1(\calA)$, there is $y \in (0,\infty)$ such that $\nu^+(h I_{h^{-1}(y,\infty)}) < \epsilon/2$. Then
\begin{eqnarray*}
0 \leq \nu^-(h) \leq \nu^+(h) &\leq& \nu^+(h I_{h^{-1}(0,\epsilon/2)}) + \nu^+(h I_{h^{-1}(\epsilon/2,y)}) + \nu^+(h I_{h^{-1}(y,\infty)}) \\
&<& \frac{\epsilon}{2} + y \nu^+(h^{-1}(\epsilon/2,y)) + \frac{\epsilon}{2} \\
&=& \epsilon,
\end{eqnarray*}
where $\nu^+(h^{-1}(\epsilon/2,y)) \leq \nu^*(h^{-1}(\epsilon/2,\infty)) = 0$, since $h$ is a null function. Let $\epsilon \rightarrow 0$ to obtain the result.

For 3, note that if $f = g$ a.e., then $h:= f - g$ is a null function. Moreover, $h \in K_1(\calA)$ by Proposition~\ref{Kp_basic_props}(11). Thus $\nu^+(f) \leq \nu^+(g) + \nu(h) = \nu^+(g)$ and $\nu^-(f) \geq \nu^-(g) + \nu(h) = \nu^-(g)$. Similarly, $\nu^+(g) \leq \nu^+(f)$ and $\nu^-(g) \geq \nu^-(f)$, giving 3a and 3b, which together imply 3c.
\qed 

Theorem~\ref{Kp_pseudonorm} can now be proved as follows.

\noindent {\bf Proof of Theorem~\ref{Kp_pseudonorm}: }
For 1, Proposition~\ref{Kp_basic_props}(11) gives that $K_p(\calA)$ is a vector space, and Proposition~\ref{Kp_basic_props}(5) gives $K_p(\calA) \subseteq L_p(\calA)$. 

For 2, consider $h \in K_p(\calA) \subseteq L_p(\calA)$. Then by Theorem~\ref{nullmodification_fn}, there is $h^{\prime} \in K_p(\alpha(\calA \cup \calN))$ such that $h = h^{\prime}$ a.e. with respect to $\alpha(\calA \cup \calN)$, $\nu_N(\lvert h^{\prime} \rvert^p)$ converges to $\int \lvert h^{\prime} \rvert^p d\nu_{\alpha(\calA \cup \calN)}$ and $\nu_N(h^{\prime})$ converges to $\int h^{\prime} d\nu_{\alpha(\calA \cup \calN)}$ as $N \rightarrow \infty$. But $h \in K_p(\alpha(\calA \cup \calN))$, and thus Lemma~\ref{K1_equal_ae} gives the result.

Properties~3 and~4 follow from~2, since $\left( \nu( \lvert h \rvert^p ) \right)^{1/p}$ inherits the properties of a pseudonorm from $\left( \int \lvert h \rvert^p d\nu \right)^{1/p}$.

For 5, note $h_1, h_2 \in L_p(\calA)$ with $\| h_1 - h_2 \|_{L_p(\calA)} = \| h_1 - h_2 \|_{K_p(\calA)}$. The result then follows by Theorem~4.4.13(ix) of~\cite{bhaskararao1983} (see also Comment~1.5 of \cite{basile2000}). 

Claim~6 follows from~5, since then the equivalence relation $\sim$ used in Definition~\ref{Kp} (ie. almost everywhere equivalence) is the same equivalence relation induced by the pseudonorm $\| \cdot \|_p$. 
\qed

The next proposition lists additional properties of $K_p$ spaces that are similar to familiar properties of $L_p$ spaces. The proof of these is facilitated by Theorem~\ref{Kp_pseudonorm}.

\begin{proposition} \label{integrable_properties} 
Consider a field of subsets $\calA \subset \calF$, and let $f$ and $g$ be real-valued functions on $\N$. 
\begin{enumerate}
\item If $f, g \in K_1(\calA)$ and $c, d \in \mathbb{R}$, then $\nu(cf + dg) = c \nu(f) + d \nu(g)$.
\item If $f, g \in K_1(\calA)$ and $f \leq g$ a.e., then $\nu(f) \leq \nu(g)$.
\item If $f, g \in K_1(\calA)$ and $f = g$ a.e. then $\nu(f) = \nu(g)$.
\item ({\bf Dominated Convergence Theorem I}) Suppose $g \in K_p(\calA)$ for some $p \in [1,\infty)$, and let $\{ f_k \}_{k=1}^{\infty}$ be a sequence of $T_1$-measurable functions on $\N$ such that $\lvert f_k \rvert \leq g$ for each $k \in \N$. Then $f_k \xrightarrow{h} f$ if and only if $f \in K_p(\calA)$ and $\| f_k - f \|_{K_p(\calA)} \rightarrow 0$.
\item ({\bf Dominated Convergence Theorem II}) Suppose $f, g \in K_p(\calA)$ for some $p \in [1,\infty)$, and let $\{ f_k \}_{k=1}^{\infty} \subseteq K_p(\calA)$ satisfy $\lvert f_k \rvert \leq g$ a.e. for each $k \in \N$. Then $f_k \xrightarrow{h} f$ if and only if $\| f_k - f \|_{K_p(\calA)} \rightarrow 0$.

\end{enumerate}
\end{proposition}

\noindent {\bf Proof: }
For 1, recall $cf + dg \in K_1(\calA)$ by Proposition~\ref{Kp_basic_props}(11). The convergence of $\nu_N(f)$, $\nu_N(g)$ and $\nu_N(cf + dg)$ as $N \rightarrow \infty$ follows from Theorem~\ref{Kp_pseudonorm}(2), and then $\nu(cf + dg) = c \nu(f) + d \nu(g)$ by the linearity of limits. 

For 2, Theorem~\ref{Kp_pseudonorm}(2) gives $\nu(f) = \int f d\nu \leq \int g d\nu = \nu(g)$. Property~3 then follows straightforwardly from Property~2.

For 4, note Proposition~\ref{Kp_basic_props}(7) gives $\{ f_k \}_{k=1}^{\infty} \subseteq K_p(\calA)$ and Proposition~\ref{Kp_basic_props}(5) gives $g \in L_p(\calA)$ and $\{ f_k \}_{k=1}^{\infty} \subseteq L_p(\calA)$. For the forward implication, Corollary~4.4.9 of \cite{bhaskararao1983} gives that $f$ is $T_1$-measurable, and then Proposition~\ref{Kp_basic_props}(7)
gives $f \in K_p(\calA) \subseteq L_p(\calA)$. Hence 
\[
\| f_k - f \|_{K_p(\calA)} = \| f_k - f \|_{L_p(\calA)} \rightarrow 0.
\]
where the first equality follows by Theorem~\ref{Kp_pseudonorm}, and convergence to zero follows by the dominated convergence theorem for $L_p$ spaces (Theorem~4.6.14 of~\cite{bhaskararao1983}). For the reverse implication, Theorem~\ref{Kp_pseudonorm} gives $\| f_k - f \|_{L_p(\calA)} = \| f_k - f \|_{K_p(\calA)}$ for each $k \in \N$, hence dominated convergence gives $f_k \xrightarrow{h} f$. The proof of 5 is similar to that of Property~4.
\qed

Theorem~\ref{Kp_pseudonorm} provides an isometric embedding of $K_p(\calA)$ in $L_p(\calA)$, and by implication, of $\calK_p(\calA)$ in $\calL_p(\calA)$. However, it leaves open the question of whether this embedding is surjective, and thus an isomorphism. The next theorem characterises fields of sets for which this embedding is surjective.

\begin{theorem} \label{Kp_Lp_isomorphism}
Let $\calA \subset \calF$ be a field of sets and let $p \in [1,\infty)$. 
\begin{enumerate}
\item If $\calN \subset \overline{\calA}$, then $\calK_p(\calA) \cong \calL_p(\calA)$ with the isomorphism given by $[h]_{K_p(\calA)} \mapsto [h]_{L_p(\calA)}$ for each $h \in K_p(\calA)$.
\item The following statements are logically equivalent.
\begin{enumerate}
\item $\calL_p(\calA) \cong \calK_p(\alpha(\calA \cup \calN))$ with the isomorphism mapping $[f]_{L_p(\calA)}$ to $[h]_{K_p(\alpha(\calA \cup \calN))}$, where $f$ and $h$ are as described in Theorem~\ref{nullmodification_fn}.
\item $\overline{\alpha(\calA \cup \calN)} = \alpha(\overline{\calA} \cup \calN)$.
\end{enumerate}
\end{enumerate}
\end{theorem}

\noindent {\bf Proof: } 
For 1, Theorem~\ref{Kp_pseudonorm} gives that $[h]_{K_p(\overline{\calA})} \mapsto [h]_{L_p(\overline{\calA})}$ is an isometry, and Theorem~\ref{nullmodification_fn} implies this map is surjective. Hence $\calK_p(\overline{\calA}) \cong \calL_p(\overline{\calA})$. The result then follows by Proposition~1.8 of~\cite{basile2000} and Proposition~\ref{Kp_basic_props}(2) (of this paper).

For 2, first suppose 2a. Then Statement~1 implies $\calK_p(\alpha(\calA \cup \calN)) \cong \calL_p(\alpha(\calA \cup \calN))$, with the isomorphism given by $[h]_{K_p(\alpha(\calA \cup \calN))} \mapsto [h]_{L_p(\alpha(\calA \cup \calN))}$ for each $h \in K_p(\alpha(\calA \cup \calN))$. But Theorem~\ref{nullmodification_fn} gives $f = h$ a.e. with respect to $\alpha(\calA \cup \calN)$, hence $[h]_{L_p(\alpha(\calA \cup \calN))} = [f]_{L_p(\alpha(\calA \cup \calN))}$. Chaining the isomorphisms gives $\calL_p(\calA) \cong \calL_p(\alpha(\calA \cup \calN))$ with the isomorphism given by $[f]_{L_p(\calA)} \mapsto [f]_{L_p(\alpha(\calA \cup \calN))}$, and then Theorem~5.10(2) of~\cite{keith2022} gives 2b.

Conversely, suppose 2b. Then Theorem~5.10(2) of~\cite{keith2022} gives $\calL_p(\calA) \cong \calL_p(\alpha(\calA \cup \calN))$ with the isomorphism given by $[f]_{L_p(\calA)} \mapsto [f]_{L_p(\alpha(\calA \cup \calN))}$ for each $f \in L_p(\calA)$. Apply Theorem~\ref{nullmodification_fn} to $f$ to construct $h \in K_p(\alpha(\calA \cup \calN))$, then $[f]_{L_p(\alpha(\calA \cup \calN))} = [h]_{L_p(\alpha(\calA \cup \calN))}$. Statement~1 then gives $\calK_p(\alpha(\calA \cup \calN)) \cong \calL_p(\alpha(\calA \cup \calN))$, with the isomorphism mapping $[h]_{K_p(\alpha(\calA \cup \calN))} \mapsto [h]_{L_p(\alpha(\calA \cup \calN))} = [f]_{L_p(\alpha(\calA \cup \calN))}$. Inverting the latter isomorphism and chaining the isomorphisms gives 2a.
\qed

Statement~2 of Theorem~\ref{Kp_Lp_isomorphism} identifies conditions under which $\calL_p(\calA)$ is isomorphic to $\calK_p(\alpha(\calA \cup \calN))$, but this requires augmenting the null sets to $\calA$. It leaves unanswered whether $\calK_p(\alpha(\calA \cup \calN))$ is in turn isomorphic to $\calK_p(\calA)$. The final theorem of this section addresses this question.

\begin{theorem} \label{Kp_null_sets_add_no_functions}
Consider a field of sets $\calA \subset \calF$ and $p \in [1,\infty)$. Then $K_p(\calA)$ is a dense subspace of $K_p(\alpha(\calA \cup \calN))$. Moreover,
\begin{enumerate}
\item $K_p(\alpha(\calA \cup \calN)) = K_p(\calA)$ if and only if $\calN \subset \overline{\calA}$, and
\item $\calK_p(\calA) \cong \calK_p(\alpha(\calA \cup \calN))$ with isomorphism $[h]_{K_p(\calA)} \mapsto [h]_{K_p(\alpha(\calA \cup \calN))}$ if and only if both of the following conditions hold: 
\begin{enumerate}
\item $\overline{\alpha(\calA \cup \calN)} = \alpha(\overline{\calA} \cup \calN)$, and
\item $\calK_p(\calA) \cong \calL_p(\calA)$ with isomorphism $[h]_{K_p(\calA)} \mapsto [h]_{L_p(\calA)}$. 
\end{enumerate}
\end{enumerate}
\end{theorem}

\noindent {\bf Proof: } 
As argued in the proof of Theorem~5.10 in~\cite{keith2022}, for every simple function $s$ with respect to $\alpha(\calA \cup \calN)$, there is a simple function $s^{\prime}$ with respect to $\calA$, such that $\| s - s^{\prime} \|_{L_p(\calA)} = 0$. Recall $K_p(\alpha(\calA \cup \calN))$ contains the simple functions with respect to $\alpha(\calA \cup \calN)$ (and therefore also the simple functions with respect to $\calA$) by Proposition~\ref{Kp_basic_props}(1). Also note Theorem~\ref{Kp_pseudonorm}(5) gives $\| s - s^{\prime} \|_{K_p(\calA)} = 0$. Simple functions with respect to $\alpha(\calA \cup \calN)$ are dense in $L_p(\alpha(\calA \cup \calN))$ by Theorem~4.6.15 of~\cite{bhaskararao1983}, and hence also dense in the subspace $K_p(\alpha(\calA \cup \calN))$. But then so are simple functions with respect to $\calA$, implying $K_p(\calA)$ is dense in $K_p(\alpha(\calA \cup \calN))$.

To show~1, first note $\calN \subset \overline{\calA}$ if and only if $\overline{\alpha(\calA \cup \calN)} = \overline{\calA}$ (Lemma~5.7 of~\cite{keith2022}). So suppose $\overline{\alpha(\calA \cup \calN)} = \overline{\calA}$. Proposition~\ref{Kp_basic_props}(2) gives
\[
K_p(\alpha(\calA \cup \calN)) = K_p(\overline{\alpha(\calA \cup \calN)}) = K_p(\overline{\calA}) = K_p(\calA).
\]
Conversely, suppose $K_p(\alpha(\calA \cup \calN)) = K_p(\calA)$, and consider $A \in \overline{\alpha(\calA \cup \calN)}$. Then $I_A \in K_p(\alpha(\calA \cup \calN))$ by Proposition~\ref{Kp_basic_props}(1), hence $I_A \in K_p(\calA)$. By Definition~\ref{Kp}, there is $y \in (0,1)$ such that $A = I_A^{-1}(y,\infty) \in \overline{\calA}$.

For 2, first suppose $\calK_p(\alpha(\calA \cup \calN)) \cong \calK_p(\calA)$ with the isomorphism given by $[h]_{K_p(\calA)} \mapsto [h]_{K_p(\alpha(\calA \cup \calN))}$, and consider $A \in \overline{\alpha(\calA \cup \calN)}$. As in the proof of 1, $I_A \in K_p(\alpha(\calA \cup \calN))$, hence $[I_A] \in \calK_p(\alpha(\calA \cup \calN))$. By assumption, there exists $h \in K_p(\calA)$  such that $[h]_{K_p(\alpha(\calA \cup \calN))} = [I_A]_{K_p(\alpha(\calA \cup \calN))}$, which implies $h = I_A$ a.e. with respect to $\alpha(\calA \cup \calN)$. There exists some $y \in (0,1)$ such that $I_A^{-1}(y,\infty) \triangle h^{-1}(y,\infty) \in \calN$ (by Theorem~3.11 of~\cite{keith2022}) and also such that $h^{-1}(y,\infty) \in \overline{\calA}$ (by Theorem~3.4 of~\cite{keith2022}, since $h$ is $T_1$-measurable with respect to $\calA$). Hence $A = I_A^{-1}(y,\infty)$ differs from a set in $\overline{\calA}$ by a null set. This implies $A \in \alpha(\overline{\calA} \cup \calN)$, by Lemma~5.3 of~\cite{keith2022}. Thus $\overline{\alpha(\calA \cup \calN)} \subseteq \alpha(\overline{\calA} \cup \calN)$. Moreover, $\alpha(\overline{\calA} \cup \calN) \subseteq \overline{\alpha(\calA \cup \calN)}$, since $\overline{\calA} \subseteq \overline{\alpha(\calA \cup \calN)}$ and $\calN \subseteq \overline{\alpha(\calA \cup \calN)}$, giving $\overline{\alpha(\calA \cup \calN)} = \alpha(\overline{\calA} \cup \calN)$. 

By Theorem~5.10 of~\cite{keith2022}, $\calL_p(\calA) \cong \calL_p(\alpha(\calA \cup \calN))$ with isomorphism $[h]_{L_p(\calA)} \mapsto [h]_{L_p(\alpha(\calA \cup \calN))}$, and then by Theorem~\ref{Kp_Lp_isomorphism},
\[
\calK_p(\calA) \cong \calK_p(\alpha(\calA \cup \calN)) \cong \calL_p(\alpha(\calA \cup \calN)) \cong \calL_p(\calA).
\]
Moreover, the first two isomorphisms map 
\[
[h]_{K_p(\calA)} \mapsto [h]_{K_p(\alpha(\calA \cup \calN))} \mapsto [h]_{L_p(\alpha(\calA \cup \calN))}
\]
for $h \in K_p(\calA)$. Since $K_p(\calA) \subseteq L_p(\calA)$ (Proposition~\ref{Kp_basic_props}(5)), the third isomorphism maps $[h]_{L_p(\calA)} \mapsto [h]_{L_p(\alpha(\calA \cup \calN))}$, and hence its inverse maps $[h]_{L_p(\alpha(\calA \cup \calN))} \mapsto [h]_{L_p(\calA)}$. The combined isomorphism maps $[h]_{K_p(\calA)} \mapsto [h]_{L_p(\calA)}$.

For the converse, 2(a) implies $\calL_p(\calA) \cong \calL_p(\alpha(\calA \cup \calN))$ with isomorphism $[h]_{L_p(\calA)} \mapsto [h]_{L_p(\alpha(\calA \cup \calN))}$, by Theorem~5.10 of~\cite{keith2022}. Combining this with 2(b) gives $\calK_p(\calA) \cong \calL_p(\alpha(\calA \cup \calN))$, with the isomorphism given by $[h]_{K_p(\calA)} \mapsto [h]_{L_p(\alpha(\calA \cup \calN))}$. Theorem~\ref{Kp_Lp_isomorphism} gives $\calK_p(\alpha(\calA \cup \calN)) \cong \calL_p(\alpha(\calA \cup \calN))$ with isomorphism $[h]_{K_p(\alpha(\calA \cup \calN))} \mapsto [h]_{L_p(\alpha(\calA \cup \calN))}$ for $h \in K_p(\alpha(\calA \cup \calN))$. Since $K_p(\calA) \subseteq K_p(\alpha(\calA \cup \calN)) \subseteq L_p(\alpha(\calA \cup \calN))$ (by Proposition~\ref{Kp_basic_props}(5)), the inverse of the latter isomorphism maps $[h]_{L_p(\alpha(\calA \cup \calN))} \mapsto [h]_{K_p(\alpha(\calA \cup \calN))}$ for each $h \in K_p(\calA)$. Hence chaining these isomorphisms gives $\calK_p(\alpha(\calA \cup \calN)) \cong \calK_p(\calA)$ with isomorphism $[h]_{K_p(\calA)} \mapsto [h]_{K_p(\alpha(\calA \cup \calN))}$.
\qed

Theorems~\ref{Kp_Lp_isomorphism} and~\ref{Kp_null_sets_add_no_functions} imply that many of the properties of $\calK_p(\calA)$ spaces can be learned through study of $\calL_p(\calA)$ spaces. However, as several of the results in this section have demonstrated, care must be taken in adapting theorems for $L_p$ spaces to obtain theorems for $K_p$ spaces, when those theorems involve equality or dominance of functions almost everywhere (as many do). The $K_p$ spaces are in general proper subsets of $L_p$ spaces, and this may necessitate modifying the conditions or conclusions of key theorems when adapting them for $K_p$ spaces.

The remainder of this paper focuses on the completeness and separability of $L_p(\calA)$ spaces, and by implication the $K_p(\calA)$ spaces isomorphic to them. 

\section{$L_p$ spaces over $\sigma$-fields on which $\nu$ is countably additive}
\label{Lp_countably_additive}

This short section considers the situation in which $(\N,\calA,\nu)$ is a measure space, that is, $\calA$ is a $\sigma$-field and $\nu$ is countably additive. The $L_p$ spaces in this case have a particularly simple form. Theorem~\ref{nuisameasure} characterises those fields of sets for which this situation occurs, and Theorem~\ref{measure_space_functions} characterises functions that are measurable (in the standard measure theory sense) with respect to such a field.

\begin{theorem} \label{nuisameasure}
Suppose $\calA \subseteq \calP(\N)$ is a $\sigma$-field. Then the following statements are logically equivalent.
\begin{enumerate}
\item $\calA \subset \calF$ and $\nu$ is countably additive on $\calA$. 
\item $\calA$ is generated by a countable partition $\calS$ of $\N$ (that is, $\calA = \sigma(\calS)$) such that $\calS \subset \calF$ and $\sum_{A \in \calS} \nu(A) = 1$.
\end{enumerate}
\end{theorem}

\begin{theorem} \label{measure_space_functions}
Suppose $(\N,\calA,\nu)$ is a measure space, where $\calA \subset \calF$. Then measurable functions on $(\N,\calA,\nu)$ are of the form
\[
f = \sum_{j = 1}^J a_j I_{A_j} + g,
\]
where $J$ may be finite or (countably) infinite, $a_j \in \mathbb{R}$ and $A_j \in \calA$ with $\nu(A_j) > 0$ for each $j$, the sets $\{ A_j \}_{j=1}^{J}$ are pairwise disjoint, and $g$ is non-zero only on a null set disjoint from $\cup_{j=1}^J A_j$. 
\end{theorem}


The proofs of these two theorems take advantage of the fact $\N$ is countable. The following two lemmas identify consequences of this for all $\sigma$-fields in $\N$.

\begin{lemma} \label{arbitraryunions}
Any $\sigma$-field $\calA$ comprised of subsets of $\N$ is closed under arbitrary unions (that is, including uncountable unions).
\end{lemma}

\noindent {\bf Proof: }
Consider $\calC \subseteq \calA$ and let $C := \bigcup \calC$. For each $n \in C$, there exists some $C_n \in \calC$ such that $n \in C_n \subseteq C$. Hence $C = \cup_{n \in C} C_n$. Since $C$ is countable, this is a countable union and hence $C \in \calA$. 
\qed

\begin{lemma} \label{countablegenerators}
Any $\sigma$-field $\calA \subseteq \calP(\N)$ is generated by a countable partition of $\N$.
\end{lemma}

\noindent {\bf Proof: }
For each $k \in \N$, let $A_k$ be the intersection of all sets in $\calA$ that contain $k$. Then $A_k \in \calA$ by Lemma~\ref{arbitraryunions}. The collection $\{ A_k \}_{k=1}^{\infty}$ generates a $\sigma$-field $\{ \cup_{k \in \calI} A_k : \calI \subseteq \N \}$ that contains every element of $\calA$ and hence $\calA = \{ \cup_{k \in \calI} A_k : \calI \subseteq \N \}$. For distinct $j, k \in \N$, either $A_j = A_k$ or $A_j \cap A_k = \emptyset$, otherwise $A_j \setminus A_k \in \calA$ would be a proper subset of $A_j$ containing $j$. Thus $\{ A_k \}_{k=1}^{\infty}$ forms a partition of $\N$.
\qed

These properties of $\sigma$-fields over $\N$ entail that requiring $(\N,\calA,\nu)$ to be a measure space imposes severe restrictions on the structure of $\calA$. These are described in Theorem~\ref{nuisameasure}, which can now be proved.

\noindent {\bf Proof of Theorem~\ref{nuisameasure}: }
$(1 \implies 2)$ By Lemma~\ref{countablegenerators}, $\calA$ is generated by a countable partition $\calS$ of $\N$. Thus $\sum_{A \in \calS} \nu(A) = \nu(\N) = 1$ by the countable additivity of $\nu$ on $\calA$.

$(2 \implies 1)$ Suppose $\calS = \{ A_k \}_{k=1}^{\infty}$. (This assumes $\calS$ is infinite, but the proof that follows also works for finite $\calS$, with minor modifications.) Recall $\calA = \{ \cup_{k \in \calI} A_k : \calI \subseteq \N \}$. Thus it is sufficient to show $\cup_{k \in \calI} A_k \in \calF$ with $\nu(\cup_{k \in \calI} A_k) = \sum_{k \in \calI} \nu(A_k)$ for any $\calI \subseteq \N$.
 
Apply Corollary~\ref{nullmodification_sequence} to obtain sets $\{ A_k^{\prime} \}_{k=1}^{\infty} \subset \calF$ and $\{ F_k \}_{k=1}^{\infty} \subset \calN$ such that for each $k \in \N$:
\begin{enumerate}
\item $A_k = A_k^{\prime} \cup F_k$, 
\item $\nu(A_k^{\prime}) = \nu(A_k)$, 
\item $\nu_N(A_k^{\prime}) \leq \nu(A_k^{\prime})$ for all $N \in \N$, and
\item $\cup_{k=1}^{\infty} F_k \in \calF$ with
\[
\nu(\cup_{k=1}^{\infty} F_k) = \nu(\cup_{k=1}^{\infty} A_k \setminus \cup_{k=1}^{\infty} A_k^{\prime}) = \nu(\N) - \sum_{k=1}^{\infty} \nu(A_k) = 0.
\]
\end{enumerate}

By Proposition~\ref{calFproperties}(9), $\cup_{k \in \calI} A_k^{\prime} \in \calF$ with 
\[
\nu( \cup_{k \in \calI} A_k^{\prime}) = \sum_{k \in \calI} \nu(A_k^{\prime}) = \sum_{k \in \calI} \nu(A_k),
\]
for any $\calI \subseteq \N$. Hence 
\[
\cup_{k=1}^{\infty} A_k = (\cup_{k \in \calI} A_k^{\prime}) \cup (\cup_{k \in \calI} F_k) \in \calF
\] 
with
\[
\nu(\cup_{k \in \calI} A_k) = \nu( \cup_{k \in \calI} A_k^{\prime}) + \nu(\cup_{k \in \calI} F_k) =  \sum_{k \in \calI} \nu(A_k)
\]  
as required.
\qed 

As a consequence, $L_p(\calA)$ spaces are rather simple when $(\N,\calA,\nu)$ is a measure space, since all functions have the form described in Theorem~\ref{measure_space_functions}. This can be proved as follows.

\noindent {\bf Proof of Theorem~\ref{measure_space_functions}: }
Let $f \in L_0(\calA)$. Then the sets $\{f^{-1}(a)\}_{a \in \mathbb{R}}$ form a partition of $\N$, and only countably many of them can be non-empty by Theorem~\ref{nuisameasure}. Hence
\[
f = \sum_{j = 1}^J a_j I_{A_j} + \sum_{j = 1}^K b_j I_{B_j},
\]
where $J$ and $K$ may be finite or infinite, $\{a_j\}_{j=1}^J \cup \{b_j\}_{j=1}^K \subset \mathbb{R}$, and $\{A_j\}_{j=1}^J \cup \{B_j\}_{j=1}^K \subset \calA$ are disjoint sets with $\nu(A_j) > 0$ for each $j$ and $\nu(B_j) = 0$ for each $j$. But then $\cup_{j=1}^K B_j \in \calN$ by countable additivity.
\qed 

The $L_p$ spaces comprised of such functions provide first examples of complete, separable $L_p(\calA)$ spaces. For any $p \in \{0\} \cup [1,\infty)$, these spaces are complete by Corollary~3.7 of~\cite{keith2022}. They are also separable, since the simple functions with rational coefficients form a countable, dense subset of $L_p(\calA)$.

\section{Completeness of $L_p(\calA)$}
\label{complete_calF}

This section characterises complete $L_p(\calA)$ spaces. It contains two complementary theorems. The first of these embeds $\calL_p(\calA)$ in a conventional Lebesgue function space, and thus provides access to the familiar properties and theorems of Lebesgue integration for application to functions in $L_p(\calA)$. This is desirable, because although many of the properties and theorems of the Lebesgue integral have analogues in charge spaces (see~\cite{leader1953}, \cite{appling1974}, \cite{chen1976}, \cite{karandikar1982} and \cite{karandikar1988} for some examples), these analogous results typically have different or additional conditions that make them more difficult to use, or at least less familiar. 

The first theorem and its proof refer to the concept of a {\em representation} of the Boolean quotient $\calP(\N) / \calN$, that is, a Boolean isomorphism $\phi$ mapping $\calP(\N) / \calN$ to a field of subsets of some set $Y$. Stone's representation theorem (originally proved in~\cite{stone1937}, but~\cite{givant2009} and~\cite{sikorski1969boolean} contain helpful expositions) is a fundamental result asserting that every abstract Boolean algebra $\calA$ is isomorphic to a field of sets. More specifically, $\calA$ can be embedded in $\calP(S_{\calA})$, where $S_{\calA}$ is a set called the {\em Stone space} of $\calA$. Stone spaces can be constructed in various equivalent ways, but the details are not needed here. It will be sufficient throughout what follows to work with an unspecified representation of $\calP(\N) / \calN$.

The theorem and proof also refer to the set $[\calF]$ defined in Section~\ref{quotient_space}, and the set $\langle \calF \rangle_{\phi} := \{ \langle A \rangle_{\phi} : A \in \calF \}$, where $\langle A \rangle_{\phi} := \phi[A]$. 

\begin{theorem} \label{isomorphism}
Consider a field of sets $\calA \subset \calF$ and $p \in \{0\} \cup [1,\infty)$. Let $\phi : \calP(\N) / \calN \rightarrow \calP(Y)$ be a representation of the Boolean algebra $\calP(\N) / \calN$. Then 
\begin{enumerate}
\item $\phi(\sigma[\calA]) \subset \langle \calF \rangle_{\phi}$, 
\item $\big( Y, \overline{\phi(\sigma[\calA])}, \overline{\nu} \big)$ is a complete measure space, 
\item $\calL_p \big( Y, \overline{\phi(\sigma[\calA])}, \overline{\nu} \big)$ is a Lebesgue function space, 
\item $\calL_p(\calA)$ is isometrically isomorphic to a dense subspace of $\calL_p \big( Y,\phi(\sigma[\calA]),\nu \big)$, and 
\item $\calL_p(\calA)$ is complete if and only if $\calL_p(\calA) \cong \calL_p \big( Y,\phi(\sigma[\calA]),\nu \big)$.
\end{enumerate}
\end{theorem}


The proof uses the following lemma, which resembles a well known result sometimes set as an exercise for students of measure theory. Recall that the symbol `+' in this context represents exclusive disjunction, that is $p + q = (p \wedge q^{\prime}) \vee (p^{\prime} \wedge q)$, where $p$ and $q$ are elements of a Boolean algebra.

\begin{lemma} \label{sigma_approximation}
Suppose $\calA$ is a Boolean algebra, and $\calM \subseteq \calA$ is a monotone class on which a non-negative, bounded, countably additive function $\mu$ is defined. Let $\calA_0 \subseteq \calM$ be a sub-algebra. Then for any $B \in \sigma(\calA_0)$ and $\epsilon > 0$, there exists $A \in \calA_0$ such that $\mu(B + A) < \epsilon$.
\end{lemma}

\noindent {\bf Proof: }
Let $\calB$ be the collection of elements $B \in \calM$ such that for any $\epsilon > 0$, there exists $A \in \calA_0$ with $\mu(B + A) < \epsilon$. Let $\{ B_k \}_{k=1}^{\infty}$ be a non-decreasing sequence in $\calB$ and choose $\epsilon > 0$. Set $B := \vee_{k=1}^{\infty} B_k$. Then there is $B_k$ with $\mu(B + B_k) < \epsilon/2$ and $A \in \calA_0$ with $\mu(B_k + A) < \epsilon/2$. Hence $\mu(B + A) \leq \epsilon$, implying $\calB$ is a monotone class containing $\calA_0$, and hence $\sigma(\calA_0) \subseteq \calB$ by the monotone class theorem (Theorem~\ref{monotone_for_Boolean}).
\qed

\noindent {\bf Proof of Theorem~\ref{isomorphism}: }
First note $\sigma[\calA] \subset [\calF]$ with $\nu$ countably additive on $\sigma[\calA]$, by Corollary~\ref{countableadditivity} and the monotone class theorem for Boolean algebras (Theorem~\ref{monotone_for_Boolean}). Then $\phi(\sigma[\calA]) \subset \langle \calF \rangle_{\phi}$. Moreover $\big( Y,\overline{\phi(\sigma[\calA])},\overline{\nu} \big)$ is a complete measure space by Lemma~4.1 of~\cite{keith2022}. Corollaries~3.5 and~3.6 of~\cite{keith2022} give that $\calL_p\big( Y,\overline{\phi(\sigma[\calA])},\overline{\nu} \big)$ is a Lebesgue function space, which is therefore complete (Corollary~3.7 of~\cite{keith2022}). Proposition~1.8 of~\cite{basile2000} then gives $\calL_p\big( Y,\overline{\phi(\sigma[\calA])},\overline{\nu} \big) = \calL_p \big( Y,\phi(\sigma[\calA]),\nu \big)$.

Define a map $\rho: \calL_p(\calA) \rightarrow \calL_p \big( Y,\phi(\sigma[\calA]),\nu \big)$ in two stages: first for equivalence classes of simple functions, then for general equivalence classes of $\calL_p(\calA)$. To define $\rho$ for the equivalence class of a simple function $f = \sum_{k=1}^K c_{k} I_{A_{k}}$, suppose without loss of generality that $\nu(A_k) > 0$ and $c_k \neq 0$ for each $k$, and that the $\{ c_k \}$ are distinct and in increasing order. (One can discard terms for which $\nu(A_k) = 0$ or $c_k = 0$, merge sets with equal coefficients and reorder to obtain a simple function with these properties in the same equivalence class $[f]$.) Then define
\[
\rho [f] := \sum_{k=1}^K c_{k} [I_{\langle A_{k} \rangle }].
\]
This mapping is well defined since if $[f]$ contains another simple function $f ^{\prime} = \sum_{k=1}^{K^{\prime}} c_{k}^{\prime} I_{A_{k}^{\prime}}$ (again with $\nu(A_k^{\prime}) > 0$ and $c_k^{\prime} \neq 0$ for each $k$ and distinct $\{ c_k^{\prime} \}$ in increasing order) then
\[
\sum_{k=1}^{K^{\prime}} c_{k}^{\prime} [I_{\langle A_{k}^{\prime} \rangle }] = \sum_{k=1}^K c_{k} [I_{\langle A_{k} \rangle }].
\]
To see this, observe that since there are only finitely many values of $\{ c_k \}$ and $\{ c_k^{\prime} \}$, for sufficiently small $\epsilon > 0$ one must have 
\[
\nu^*(\{ n \in \N : f(n) \neq f^{\prime}(n) \}) = \nu^*(\{ n \in \N : \lvert f(n) - f^{\prime}(n) \rvert > \epsilon \}) = 0,
\]
since $f$ and $f^{\prime}$ are equal a.e. But this is only possible if $K = K^{\prime}$, $c_k = c_k^{\prime}$, and $A_k \triangle A_k^{\prime} \in \calN$ (hence $\langle A_k^{\prime} \rangle  = \langle A_k \rangle $) for $k \in \{1, \ldots, K \}$.

The mapping $\rho$ is an isometry for simple functions in the case $p \in [1,\infty)$, since 
\begin{eqnarray*}
\| \rho[f] \|_p^p &=& \int \sum_{k=1}^K \lvert c_{k} \rvert^p I_{\langle A_{k} \rangle } d\nu  = \sum_{k=1}^K \lvert c_{k} \rvert^p \nu\langle A_{k} \rangle \\
&=& \sum_{k=1}^K \lvert c_{k} \rvert^p \nu(A_{k}) = \| f \|_p^p = \| [f] \|_p^p.
\end{eqnarray*}
Likewise, $\rho$ is an isometry for simple functions in the case $p = 0$, where the metrics on $\calL_0(\calA)$ and $\calL_0 \big( Y,\phi(\sigma[\calA]),\nu \big)$ are of the form described immediately following Definition~2.10 of~\cite{keith2022}. To see this, note that for simple functions $f, f^{\prime} \in L_0(\calA)$, one can express $d(\rho[f],\rho[f^{\prime}])$ and $d(f,f^{\prime})$ as parallel functions of finite sets of real numbers of the form $\{\nu \langle A_k \rangle \}_{k=1}^K$  and $\{\nu(A_k) \}_{k=1}^K$ respectively, and then use $\nu \langle A_k \rangle = \nu(A_k)$ for each $k$.

To define $\rho$ for general $[f] \in \calL_p(\calA)$, note there is a sequence of simple functions $\{ f_n \}_{n=1}^{\infty}$ with respect to $\calA$ that converges hazily to $f$. Then
\[
d(\rho([f_n]), \rho([f_m])) = d(f_n, f_m) \rightarrow 0
\]
as $m, n \rightarrow 0$. Thus $\{ \rho([f_n]) \}_{n \geq 1}$ is a Cauchy sequence in the metric space $\calL_0 \big( Y,\phi(\sigma[\calA]),\nu \big)$, which is complete. Define $\rho([f])$ to be the limit of this Cauchy sequence. Then $\rho$ is well defined, since if there is some other sequence of simple functions $\{ f_n^{\prime} \}_{n=1}^{\infty}$ that converges to $f^{\prime} \in [f]$, it follows that
\begin{eqnarray*}
d(\rho([f_n]), \rho([f_n^{\prime}])) &=& d(f_n, f_n^{\prime}) \\
&\leq& d(f_n, f) + d(f, f^{\prime}) + d(f_n^{\prime}, f^{\prime}). 
\end{eqnarray*}
The middle summand is identically zero since $f \sim f^{\prime}$, and the other two terms converge to zero as $n \rightarrow \infty$. Hence $\rho([f_n^{\prime}])$ converges to the same limit as $\rho([f_n])$.

To show $\rho$ is an isometry for general $[f] \in \calL_p(\calA)$ in the case $p \in [1, \infty)$, note that for any $n$, 
\begin{eqnarray*}
\lvert \| \rho[f] \|_p - \| [f] \|_p \rvert &\leq& \lvert \| \rho[f] \|_p - \| \rho[f_n] \|_p \rvert + \lvert \| \rho[f_n] \|_p - \| [f_n] \|_p \rvert \\
& & + \lvert \| [f_n] \|_p - \| [f] \|_p \rvert \\
& \leq & \| \rho[f] - \rho[f_n] \|_p + \| f - f_n \|_p
\end{eqnarray*}
where the middle summand on the right hand side of the first inequality is identically zero since $\rho$ is an isometry for simple functions. But the final line goes to zero as $n \rightarrow \infty$. In particular, the first summand goes to zero by dominated convergence (Theorem~4.6.14 of~\cite{bhaskararao1983}), noting $\{ f_n \}_{n=1}^{\infty}$ can be chosen to be dominated by $f$, using Theorem~3.4 of~\cite{keith2022}. Hence $\| \rho[f] \|_p = \| [f] \|_p$. For $p = 0$ and general $[f], [g] \in \calL_p(\calA)$,
\begin{eqnarray*}
d(\rho[f], \rho[g]) &\leq& d(\rho[f], \rho[f_n]) + d(\rho[f_n], \rho[g_n]) + d(\rho[g_n], \rho[g]) \\
& = & d(\rho[f], \rho[f_n]) + d([f_n], [g_n]) + d(\rho[g_n], \rho[g]) \\
& \rightarrow & \lim_{n\rightarrow \infty} d([f_n], [g_n]) \\
& = & d([f],[g]).
\end{eqnarray*}
By a similar argument, $d([f],[g]) \leq d(\rho[f], \rho[g])$.

Lemma~\ref{sigma_approximation} implies any $[A] \in \sigma[\calA]$ can be approximated arbitrarily closely by some $[A_{\epsilon}] \in [\calA]$. Hence any simple function with respect to $\phi(\sigma[\calA])$ can be approximated arbitrarily closely (in the $p$ pseudo-norm or pseudo-metric $d$) by a simple function with respect to $\phi[\calA]$. The simple functions with respect to $\phi(\sigma[\calA])$ are dense in $L_p\big( Y,\phi(\sigma[\calA]),\nu \big)$ by Theorem~4.6.15 of~\cite{bhaskararao1983}, hence so are the simple functions with respect to $\phi[\calA]$. But then the equivalence classes of simple functions with respect to $\phi[\calA]$ are dense in $\calL_p\big( Y,\phi(\sigma[\calA]),\nu \big)$, that is, the images under $\rho$ of the equivalence classes of simple functions with respect to $\calA$ are dense in $\calL_p\big( Y,\phi(\sigma[\calA]),\nu \big)$. Consequently, $\rho(\calL_p(\calA))$ is dense in $\calL_p\big( Y,\phi(\sigma[\calA]),\nu \big)$.

Suppose $\calL_p(\calA)$ is complete. Then $\calL_p(\calA) \cong \calL_p \big( Y,\phi(\sigma[\calA]),\nu \big)$, since $\calL_p(\calA)$ is dense in $\calL_p \big( Y,\phi(\sigma[\calA]),\nu \big)$. Conversely, suppose $\calL_p(\calA) \cong \calL_p \big( Y,\phi(\sigma[\calA]),\nu \big)$. Then $\calL_p(\calA)$ is complete because $\calL_p \big( Y,\phi(\sigma[\calA]),\nu \big)$ is complete. 
\qed

The above theorem entails that $\calL_p$ spaces, and the $\calK_p$ spaces to which they are isomorphic, can ultimately be embedded in conventional Lebesgue function spaces on a complete measure space. Alternatively, instead of being ``mapped forward'' using a representation $\phi$, $\sigma[\calA]$ can be ``mapped back'' into $\calF$, from which a complete $\calL_p$ space on $\N$ can be constructed.

Recall the definition $\xi(A) := [A]$ for all $A \in \calP(\N)$. Note $\xi^{-1}[\calB] = \alpha(\calB \cup \calN)$ for any field of sets $\calB \subseteq \calP(\N)$, by Lemma~5.3 of \cite{keith2022}, since $\xi^{-1}[\calB] = \{ A \in \calP(\N) : A \triangle B \in \calN \mbox{ for some } B \in \calB \}$.

\begin{theorem} \label{completeness}
Consider a field of subsets $\calA \subset \calF$ and $p \in \{0\} \cup [1,\infty)$. Define $\lambda(\calA) := \xi^{-1}(\sigma[\calA])$. Then the following claims hold. 
\begin{enumerate}
\item $\lambda(\calA) \subset \calF$.
\item $\calL_p(\lambda(\calA))$ is complete.
\item $\calL_p(\calA)$ is isometrically isomorphic to a dense subspace of $\calL_p(\lambda(\calA))$, with the isomorphism mapping $[f]_{L_p(\calA)} \mapsto [f]_{L_p(\lambda(\calA))}$.
\item The following statements are logically equivalent:
\begin{enumerate}
\item $\calL_p(\calA)$ is complete,
\item $\calL_p(\lambda(\calA)) \cong \calL_p(\calA)$ with isomorphism $[f]_{L_p(\calA)} \mapsto [f]_{L_p(\lambda(\calA))}$, 
\item $\overline{\lambda(\calA)} = \alpha(\overline{\calA} \cup \calN)$, and
\item $\calK_p(\alpha(\calA \cup \calN))$ is complete and $\overline{\alpha(\calA \cup \calN)} = \alpha(\overline{\calA} \cup \calN)$.
\end{enumerate}
Moreover, if any of 4(a)-(d) holds, then 
\[
\calL_p(\calA) \cong \calL_p(\alpha(\calA \cup \calN)) \cong \calK_p(\alpha(\calA \cup \calN)).
\]
\end{enumerate}
\end{theorem}

\noindent {\bf Proof: } 
As in the proof of Theorem~\ref{isomorphism}, $\sigma[\calA] \subset [\calF]$, hence $\xi^{-1}(\sigma[\calA]) \subset \calF$.

Since $\lambda(\calA) / \calN = [\lambda(\calA)] = \sigma[\calA]$ and $\nu$ is countably additive on $\sigma[\calA]$, $\calL_p(\lambda(\calA))$ is complete by Theorem~4.2(2) of~\cite{keith2022} (noting $L_p(\lambda(\calA))$ is complete if and only if $L_1(\lambda(\calA))$ is complete by Theorem~3.4 of~\cite{basile2000}). Let $\phi : \calP(\N) / \calN \rightarrow \calP(Y)$ be a representation of the Boolean algebra $\calP(\N) / \calN$. By Theorem~\ref{isomorphism}, $\calL_p(\lambda(\calA)) \cong \calL_p(Y, \phi(\sigma[\calA]), \nu)$ and $\calL_p(\calA)$ is isometrically isomorphic to a dense subspace of $\calL_p(Y, \phi(\sigma[\calA]), \nu)$. Hence $\calL_p(\calA)$ is isometrically isomorphic to a dense subspace of $\calL_p(\lambda(\calA))$. Also note the former two isomorphisms map $[f]_{L_p(\calA)}$ and $[f]_{L_p(\lambda(\calA))}$ to the same element of $\calL_p(Y, \phi(\sigma[\calA]), \nu)$, for each $f \in L_p(\calA)$, hence the latter isomorphism maps $[f]_{L_p(\calA)} \mapsto [f]_{L_p(\lambda(\calA))}$.

$(4a \implies 4b)$ Suppose $\calL_p(\calA)$ is complete. Then $\calL_p(\calA) \cong \calL_p(\lambda(\calA))$ because $\calL_p(\calA)$ is isometrically isomorphic to a dense subspace of $\calL_p(\lambda(\calA))$.

$(4b \implies 4a)$ This is immediate from the fact $\calL_p(\lambda(\calA))$ is complete.

$(4b \implies 4c)$ This follows by an argument familiar from the proof of Theorem~\ref{Kp_null_sets_add_no_functions}(2) (and also Theorem~5.10(2) of~\cite{keith2022}). Consider $A \in \overline{\lambda(\calA)}$, which implies $I_A \in L_p(\lambda(\calA))$ by Proposition~1.8 of~\cite{basile2000}. By assumption, there exists $f \in L_p(\calA)$ such that $[f]_{L_p(\lambda(\calA))} = [I_A]_{L_p(\lambda(\calA))}$, which implies $f = I_A$ a.e. with respect to $\lambda(\calA)$. There exists some $y \in (0,1)$ such that $I_A^{-1}(y,\infty) \triangle f^{-1}(y,\infty) \in \calN$ (by Theorem~3.11 of~\cite{keith2022}) and also such that $f^{-1}(y,\infty) \in \overline{\calA}$ (by Theorem~3.4 of~\cite{keith2022}, since $f$ is $T_1$-measurable with respect to $\calA$). Hence $A = I_A^{-1}(y,\infty)$ differs from a set in $\overline{\calA}$ by a null set. This implies $A \in \alpha(\overline{\calA} \cup \calN)$, by Lemma~5.3 of~\cite{keith2022}. Thus $\overline{\lambda(\calA)} \subseteq \alpha(\overline{\calA} \cup \calN)$. Moreover, $\alpha(\overline{\calA} \cup \calN) \subseteq \overline{\lambda(\calA)}$, since $\calA \subseteq \lambda(\calA)$ and $\calN \subseteq \lambda(\calA)$, giving $\overline{\lambda(\calA)} = \alpha(\overline{\calA} \cup \calN)$. 

$(4c \implies 4b)$ If $\overline{\lambda(\calA)} = \alpha(\overline{\calA} \cup \calN)$ then $\alpha(\overline{\calA} \cup \calN)$ is Peano-Jordan complete, hence $\overline{\alpha(\calA \cup \calN)} = \alpha(\overline{\calA} \cup \calN)$ by Lemma~5.7(2) of~\cite{keith2022}. By Theorem~5.10(2) of~\cite{keith2022}, $\calL_p(\calA) \cong \calL_p(\alpha(\calA \cup \calN))$ with isomorphism $[f]_{L_p(\calA)} \mapsto [f]_{L_p(\alpha(\calA \cup \calN))}$, and by Proposition~1.8 of~\cite{basile2000}, 
\[
L_p(\alpha(\calA \cup \calN)) = L_p(\overline{\alpha(\calA \cup \calN)}) = L_p(\alpha(\overline{\calA} \cup \calN)) = L_p(\overline{\lambda(\calA)}) = L_p(\lambda(\calA)),
\]
giving~4b. Moreover, $\overline{\alpha(\calA \cup \calN)} = \alpha(\overline{\calA} \cup \calN)$ is also Condition~2b of Theorem~\ref{Kp_Lp_isomorphism}, hence $\calL_p(\calA) \cong \calK_p(\alpha(\calA \cup \calN))$

$(4a \implies 4d)$ The preceding parts already establish 4a implies $\calL_p(\calA) \cong \calK_p(\alpha(\calA \cup \calN))$ and $\overline{\alpha(\calA \cup \calN)} = \alpha(\overline{\calA} \cup \calN)$. But then $\calK_p(\alpha(\calA \cup \calN))$ must also be complete.

$(4d \implies 4a)$ By Theorem~\ref{Kp_Lp_isomorphism}(1), $\calK_p(\alpha(\calA \cup \calN)) \cong \calL_p(\alpha(\calA \cup \calN))$, and by Theorem~5.10(2) of~\cite{keith2022}, $\calL_p(\calA) \cong \calL_p(\alpha(\calA \cup \calN))$, since $\overline{\alpha(\calA \cup \calN)} = \alpha(\overline{\calA} \cup \calN)$. Hence $\calL_p(\calA)$ is complete.
\qed 

Condition~4c of Theorem~\ref{completeness} is a convenient and surprising characterisation of complete $\calL_p(\calA)$ spaces, which also implies $\calK_p(\alpha(\calA \cup \calN))$ is complete. Condition~4c does not appear to imply $\calK_p(\calA)$ is complete, but if it is, then Condition~4c implies
\[
\calK_p(\calA) \cong \calL_p(\calA) \cong \calL_p(\alpha(\calA \cup \calN)) \cong \calK_p(\alpha(\calA \cup \calN)),
\]
because $\calK_p(\calA)$ is dense in $\calK_p(\alpha(\calA \cup \calN))$, by Theorem~\ref{Kp_null_sets_add_no_functions}.

\section{Separability of $L_p(\calA)$}
\label{separability_section}

This section concerns sufficient conditions for $L_p(\calA)$ to be separable, that is, to have a countable dense subset. Even conventional Lebesgue function spaces are not separable in general; however, there is a well known sufficient condition, described in the following definition and proposition.

\begin{definition} \label{separable_measure_space}
A charge space $(X,\calA,\mu)$ is said to be {\em separable} if there is a countable subset $\calC \subseteq \calA$ such that for any $A \in \calA$ and $\epsilon > 0$, there is $A_{\epsilon} \in \calC$ with $\mu(A \triangle A_{\epsilon}) < \epsilon$.
\end{definition}

\begin{proposition} \label{separable_Lp}
If a charge space $(X,\calA,\mu)$ is separable, then $L_p(X,\calA,\mu)$ is separable for $p \in \{0\} \cup [1,\infty)$.
\end{proposition}

Proving the preceding lemma for measure spaces is sometimes set as an exercise for students of measure theory. The proof also applies with minimal modification to the $L_p$ spaces constructed on a charge space, and involves constructing a countable dense subset of $L_p(X,\calA,\mu)$, comprised of simple functions with rational coefficients and indicator functions only for sets in the countable set $\calC \subseteq \calA$ described in Definition~\ref{separable_measure_space}. Proving these simple functions are dense in $L_p(X,\calA,\mu)$ is straightforward with the aid of Theorem~4.6.15 of~\cite{bhaskararao1983}.

One natural way to generate a field in $\calF$ for which the corresponding $L_p$ and $K_p$ spaces are complete and separable, is to start with a countable field $\calA \subset \calF$, generate the countably complete Boolean algebra $\sigma[\calA]$, and then map it back to the inverse image $\xi^{-1}(\sigma[\calA])$. This strategy is implemented in the following theorem.

\begin{theorem} \label{chain_generated_field}
Consider a field of sets $\calA \subset \calF$. The following statements are logically equivalent.
\begin{enumerate}
\item There is a chain $\calT \subseteq \calA$ such that $[\calA] = \sigma[\alpha(\calT)]$.
\item There is a countable set $\calC \subseteq \calA$ such that $[\calA] = \sigma[\alpha(\calC)]$.
\end{enumerate}
Moreover, if either statement holds then $\calA$ is separable.
\end{theorem}


The proof of Theorem~\ref{chain_generated_field} requires the following lemma, which provides several alternative characterisations of countable fields. The lemma requires the following definitions. Consider a Boolean algebra $\calA$. A {\em partition of 1} is a collection $\Delta$ of pairwise disjoint elements of $\calA$ with supremum 1. A partition of 1 $\Delta_1$ is said to be a {\em refinement} of a partition of 1 $\Delta_2$, written $\Delta_1 \leq \Delta_2$, if for every $p_1 \in \Delta_1$ there is $p_2 \in \Delta_2$ such that $p_1 \leq p_2$.

Statement~4 of the following lemma uses the notation $\alpha(\calC)$ to represent the subalgebra of $\calA$ generated by $\calC \subseteq \calA$, generalising notation used earlier in this paper for fields of sets. 

\begin{lemma} \label{countablefield}
Consider a Boolean algebra $\calA$. The following statements are logically equivalent.
\begin{enumerate}
\item $\calA$ is countable.
\item $\calA$ has a countable generating set.
\item There is a non-decreasing sequence of finite subalgebras $\{ \calA_k \}_{k=1}^{\infty} \subseteq \calA$ such that $\calA = \cup_{k=1}^{\infty} \calA_k$.
\item There is a sequence $\{ \Delta_k \}_{k=1}^{\infty}$ of finite partitions of 1 such that $\Delta_{k+1} \leq \Delta_k$ for each $k \in \N$, and $\calA = \alpha(\cup_{k=1}^{\infty} \Delta_k)$.
\item $\calA$ is generated by a countable chain $\calT \subseteq \calA$.
\end{enumerate}
\end{lemma}

\noindent {\bf Proof: }
$(1 \implies 2)$ This is trivial, since $\calA$ can be its own generating set.

$(2 \implies 3)$ Let $\calC := \{ p_k \}_{k=1}^{\infty} \subseteq \calA$ be a countable generating set for $\calA$. (Note $\calC$ may be finite, since $p_j$ need not be distinct from $p_k$ for $j \neq k$.) For $k \in \N$, define $\calA_k = \alpha(p_1, \ldots, p_k)$. Then $\calA_1 \subseteq \calA_2 \subseteq \ldots$ is a non-decreasing sequence of finite subalgebras. Moreover, $\calA = \alpha(\calC) = \cup_{k=1}^{\infty} \calA_k$, since each element of $\calA$ may be expressed in terms of elementary operations applied to a finite subset of $\calC$. 

$(3 \implies 4)$ Since $\calA_k$ is finite, the atoms of $\calA_k$ form a finite partition $\Delta_k$ of 1. Then each element of $\calA_k$ is a finite disjunction of elements of $\Delta_k$, so $\alpha(\cup_{k=1}^{\infty} \Delta_k) = \cup_{k=1}^{\infty} \calA_k$. Also, since $\calA_k \subseteq \calA_{k+1}$, each element of $\Delta_k$ is a disjunction of atoms of $\calA_{k+1}$, so $\Delta_{k+1} \leq \Delta_k$. 

$(4 \implies 5)$ For each $k \in \N$, $\Delta_k := \{ \Delta_{k,1}, \ldots, \Delta_{k,n_k} \}$, with $\{ n_k \}_{k=1}^{\infty}$ a non-decreasing sequence. Without loss of generality, suppose each partition is ordered so that the elements of $\Delta_{k+1}$ that form $\Delta_{k,1}$ by disjunction are enumerated first, then the elements of $\Delta_{k+1}$ that form $\Delta_{k,2}$ by disjunction, and so on. Then the set
\[
\calT := \{ \vee_{j=1}^i \Delta_{k,j} : k \in \N, i \in \{ 1, \ldots, n_k \} \}
\]
forms a countable chain in $\calA$. Moreover, $\Delta_{k,j}$ is a difference of elements of $\calT$ for every $k \in \N$ and $j \in \{ 1, \ldots , n_k \}$, and every $A \in \calT$ is a finite disjunction of elements of $\cup_{k=1}^{\infty} \Delta_k$. Thus $\alpha(\cup_{k=1}^{\infty} \Delta_k)= \alpha(\calT)$.

$(5 \implies 1)$ Each element of $\calA$ is formed by performing basic operations on a finite subset of $\calT$, and thus the elements of $\calA$ can be systematically enumerated.
\qed 


Theorem~\ref{chain_generated_field} can now be proved as follows.

\noindent {\bf Proof of Theorem~\ref{chain_generated_field}: }
$(1 \implies 2)$ Given $\calT \subseteq \calA$, one may obtain a countable subchain $\calC \subseteq \calT$ such that for any $A \in \calT$ and $\epsilon > 0$ there are $B, C \in \calC$ with $[B] \leq [A] \leq [C]$ and $\nu(C) - \nu(B) < \epsilon$, as in the proof of Theorem~\ref{nullmodification_thm}. Clearly $\sigma[\alpha(\calC)] \subseteq \sigma[\alpha(\calT)]$. Now consider $A \in \calT$, and define $B_0 := \bigvee \{ [B] : B \in \calC, [B] \leq [A] \}$ and $C_0 := \bigwedge \{ [C] : C \in \calC, [A] \leq [C] \}$, where $\bigvee$ and $\bigwedge$ represent supremum and infimum respectively. Note the supremum and infimum exist and are elements of $[\calF]$ by Theorem~\ref{calF_quotient_complete}. Moreover, $B_0 \leq [A] \leq C_0$ and $\nu(C_0 - B_0) < \epsilon$ for any $\epsilon > 0$. But this implies $B_0 = [A] = C_0$, by Proposition~\ref{quotient_properties}(5). Hence $[A] \in \sigma[\alpha(\calC)]$, $[\calT] \subseteq \sigma[\alpha(\calC)]$ and $[\calA] = \sigma[\alpha(\calT)] = \sigma[\alpha(\calC)]$.

$(2 \implies 1)$ Given a countable set $\calC \subseteq \calA$, Lemma~\ref{countablefield} implies there is a countable chain $\calT^{\prime} \subseteq [\calA]$ such that $\alpha(\calT^{\prime}) = \alpha[\calC]$. But then $\sigma(\alpha(\calT^{\prime})) = \sigma(\alpha[\calC])$. Suppose $\calT^{\prime} := \{ [A_k^{\prime}] : k \in \N \}$ for some $\{ A_k^{\prime} \}_{k=1}^{\infty} \subseteq \calA$. (This implicitly invokes the axiom of countable choice.) Set $A_1 := A_1^{\prime}$ so that trivially $[A_1] = [A_1^{\prime}]$ and $\{ A_1 \}$ forms a chain in $\calA$. Inductively define for each $k \geq 2$ the sets $B_k := \bigcup \{ A_j : j < k, [A_j] \leq [A_k^{\prime}] \}$ and $C_k := \N \cap \bigcap \{ A_j : j < k, [A_j] \geq [A_k^{\prime}] \}$. That is, $B_k$ is the largest element of the finite sequence $A_1, \ldots, A_{k-1}$ with $[A_j] \leq [A_k^{\prime}]$ or $\emptyset$ if no such set exists, and $C_k$ is the smallest element of that finite sequence with $[A_j] \geq [A_k^{\prime}]$ or $\N$ if no such set exists. Set $A_k := (A_k^{\prime} \cup B_k) \cap C_k$ and note $A_k \setminus A_k^{\prime} \subseteq B_k \setminus A_k^{\prime}$ and $A_k^{\prime} \setminus A_k \subseteq A_k^{\prime} \setminus C_k$. But then since $[B_k] \leq [A_k^{\prime}] \leq [C_k]$, one must have $B_k \setminus A_k^{\prime} \in \calN$, $A_k^{\prime} \setminus C_k \in \calN$ and $A_k \triangle A_k^{\prime} \in \calN$. Thus $[A_k] = [A_k^{\prime}]$ and $\{ A_1, \ldots, A_k \}$ forms a chain in $\calA$. It follows that $\calT := \{ A_k \}_{k=1}^{\infty} \subseteq \calA$ is a chain with $[\calT] = \calT^{\prime}$. Moreover, $[\calA] = \sigma[\alpha(\calC)] = \sigma[\alpha(\calT)]$, since $[\alpha(\calC)] = \alpha[\calC] = \alpha(\calT^{\prime}) = \alpha[\calT] = [\alpha(\calT)]$.

If either statement holds, then Lemma~\ref{sigma_approximation} and Corollary~\ref{countableadditivity} together imply that for any $A \in \calA$ and $\epsilon > 0$ there is $A_{\epsilon} \in \alpha(\calC)$ with $\nu([A] + [A_{\epsilon}]) < \epsilon$. But then $\nu(A \triangle A_{\epsilon}) < \epsilon$, implying $\calA$ is separable, since $\alpha(\calC)$ is countable.
\qed

\begin{corollary} \label{Lp_Polish}
If $\calA \subset \calF$ is a field of sets and $\calT \subseteq \calA$ is a chain such that $[\calA] = \sigma[\alpha(\calT)]$, then $\calL_p(\calA)$ is a Polish space for $p \in \{0\} \cup [1,\infty)$. 
\end{corollary}


\noindent {\bf Proof: }
The function space $\calL_1(\calA)$ is complete by Corollary~\ref{countableadditivity} and Theorem~4.2(2) of~\cite{keith2022}, noting $[\calA] = \alpha(\calA \cup \calN) / \calN$. Hence $L_p(\calA)$ is complete by Theorem~3.4 of~\cite{basile2000}. Moreover, $\calA$ is separable by Theorem~\ref{chain_generated_field}, and then $\calL_p(\calA)$ is separable by Proposition~\ref{separable_Lp}.
\qed 

Note for any chain $\calT \subset \calF$, the above corollary applies to the field of sets $\calA := \xi^{-1}(\sigma[\alpha(\calT)])$.

Since any $T_1$-measurable function $f$ induces a chain of inverse images of the form $f^{-1}(a, \infty)$, where $a \in \mathbb{R}$, the above theorem opens the door to exploring complete, separable $\calK_p$ and $\calL_p$ spaces induced by $T_1$-measurable functions. Of particular interest are the $\calK_p$ and $\calL_p$ spaces generated by the maximal chains in $\calF$ discussed in Corollary~\ref{maximal_chains}, since all $\calL_p$ spaces generated by chains in $\calF$ are closed subspaces of the $\calL_p$ spaces generated by such maximal chains. Exploring $\calK_p$ and $\calL_p$ spaces generated by chains in this manner may form the subject of a future paper.

\section{Conclusion}\label{sec13}

A common strategy in mathematics is to study objects of interest in aggregrate, to uncover structural properties of the class that shed light on the inter-relationships between members. This has been a fruitful approach in probability, much of which is concerned with random variables and their various modes of convergence. Surprisingly, sequences with a Ces\`aro limit have not been studied in aggregrate, despite their obvious relevance to many topics in probability and analysis, particularly the study of ergodic processes. A discrete-time, real-valued ergodic process may be regarded as a map from a probability space into the space of real sequences with a Ces\`aro limit, with probability one. It is therefore appropriate to investigate the nature of that latter space.

This paper represents a first step in that direction. The space $\calF$ of binary sequences with a Ces\`aro limit has been studied in Sections~\ref{F_and_nu} to~\ref{quotient_space}. The main structural property to emerge from this part of the paper is Corollary~\ref{countableadditivity}, which establishes that $\calF$ can be factored to produce a monotone class $[\calF]$. This finding is significant in light of the monotone class theorem for Boolean algebras (Theorem~\ref{monotone_for_Boolean}): it implies that algebras embedded in $[\calF]$ can be expanded to countably complete algebras in $[\calF]$.

Spaces of more general sequences with a Ces\`aro limit are studied in Sections~\ref{cesaro_integrals} to~\ref{separability_section}. These spaces, dubbed $K_p(\calA)$ spaces, are comprised of sequences with a kind of measurability property (Condition~1 of Definition~\ref{Kp}) with respect to a field of sets $\calA \subset \calF$, and also a kind of integrability property (Condition~2 of Definition~\ref{Kp}). One of the main results in this second part of the paper is that $K_p(\calA)$ spaces are vector spaces with a pseudonorm defined in terms of the Ces\`aro limit, and an isometric embedding into a $L_p(\calA)$ space. In particular, for sequences in these spaces, the Ces\`aro limit is equal to an integral. Conditions under which this embedding (or the induced isometric embedding of the normed vector space $\calK_p(\calA)$ into $\calL_p(\calA)$) is an isomorphism are identified in Theorems~\ref{Kp_Lp_isomorphism} and~\ref{Kp_null_sets_add_no_functions}.

Sections~\ref{Lp_countably_additive}, \ref{complete_calF} and~\ref{separability_section} identify conditions under which $L_p$ spaces and the $K_p$ spaces isomorhpic to them are complete and separable. Section~\ref{Lp_countably_additive} clarifies why $L_p(\calA)$ and $K_p(A)$ are appropriately studied using finitely additive measure theory: the cases in which these reduce to (countably additive) measure spaces are of a highly specific form, but are nevertheless complete and separable. More generally, Theorems~\ref{isomorphism} and~\ref{completeness} provide complementary characterisations of complete $L_p(\calA)$ spaces, the first in terms of an embedding in a conventional Lebesgue function space (constructed on a countably additive measure space), and the second in terms of an expansion of the underlying field $\calA$ to a field $\lambda(\calA)$, based on Corollary~\ref{countableadditivity}. Section~\ref{separability_section} proposes a strategy for generating $\calL_p(\calA)$ spaces that are Polish spaces (complete, separable metric spaces), using a chain of sets in $\calF$, such as that consisting of the inverse images of a $T_1$-measurable function. This is flagged as a potential topic for further research. 

In the process of proving the results contained in this paper, an analytical tool dubbed {\em null modification} has been developed (a form of this technique is also described in the companion paper~\cite{keith2022}). The technique involves manipulating sets (Proposition~\ref{nullmodification}), chains (Theorem~\ref{nullmodification_thm}), sequences (Corollary~\ref{nullmodification_sequence}) or functions (Theorem~\ref{nullmodification_fn}) by adding and/or deleting null sets in such a manner as to imbue the modified object with desirable properties. This technique may be useful in other contexts.

The theory developed herein pertains to certain spaces of sequences with Ces\`aro limits, defined in terms of some field of sets $\calA \subset \calF$. An interesting direction for future research is to study in aggregate the space of all sequences with Ces\`aro limits, or at least those that satisfy Definition~\ref{Kp} with the field $\calA$ replaced by the additive class $\calF$. Spaces of this latter form are appropriately dubbed $K_p(\calF)$ spaces, and may form the subject of a future paper.

\backmatter

\bmhead{Supplementary information}

All proofs are provided in the preprint version of this paper, together with additional lemmas on which these proofs depend. There is no data associated with this manuscript.

\bmhead{Acknowledgments}

The authors are grateful to the Australian Research Council Centre of Excellence for Mathematical and Statistical Frontiers (CE140100049) for their (non-financial) support.

\section*{Declarations}

The authors have no competing interests to declare that are relevant to the content of this article.

\end{document}